
\magnification=1200

\font\tengoth=eufm10
\font\sevengoth=eufm7
\newfam\gothfam
\textfont\gothfam=\tengoth
\scriptfont\gothfam=\sevengoth
\def\goth{\fam\gothfam\tengoth}

\font\tenbboard=msbm10
\font\sevenbboard=msbm7
\newfam\bboardfam
\textfont\bboardfam=\tenbboard
\scriptfont\bboardfam=\sevenbboard
\def\bboard{\fam\bboardfam\tenbboard}

\newif\ifpagetitre
\newtoks\auteurcourant \auteurcourant={\hfill}
\newtoks\titrecourant \titrecourant={\hfill}

\pretolerance=500 \tolerance=1000 \brokenpenalty=5000
\newdimen\hmargehaute \hmargehaute=0cm
\newdimen\lpage \lpage=14.3cm
\newdimen\hpage \hpage=20cm
\newdimen\lmargeext \lmargeext=1cm
\hsize=11.25cm
\vsize=18cm
\parskip=0cm
\parindent=12pt

\def\margehaute{\vbox to \hmargehaute{\vss}}
\def\margebasse{\vss}

\output{\shipout\vbox to \hpage{\margehaute\nointerlineskip
  \corpsdepage\margebasse}
  \advancepageno \global\pagetitrefalse
  \ifnum\outputpenalty>-2000 \else\dosupereject\fi}

\def\corpsdepage{\hbox to \lpage{\hss\pagetexte\hskip\lmargeext}}
\def\pagetexte{\vbox{\makeheadline\pagebody\makefootline}}
\headline={\ifpagetitre\titleheadline \else
  \ifodd\pageno\rightheadline \else \leftheadline\fi\fi}
\def\leftheadline{\hfil\the\auteurcourant\hfil}
\def\rightheadline{\hfil\the\titrecourant\hfil}
\def\titleheadline{\hfill}
\pagetitretrue

\font\petcap=cmcsc10
\def\pc#1#2|{{\tenrm#1\sevenrm#2}}
\def\pd#1 {{\pc#1|}\ }

\def\pointir{\discretionary{.}{}{.\kern.35em---\kern.7em}\nobreak\hskip 0em
 plus .3em minus .4em }

\def\titre#1|{\message{#1}
             \par\vskip 30pt plus 24pt minus 3pt\penalty -1000
             \vskip 0pt plus -24pt minus 3pt\penalty -1000
             \centerline{\bf #1}
             \vskip 5pt
             \penalty 10000 }
\def\section#1|
               {\par\vskip .3cm
               {\bf #1}\pointir}
\def\ssection#1|
             {\par\vskip .2cm
             {\it #1}\pointir}

\long\def\th#1|#2\finth{\par\vskip 5pt
              {\petcap #1\pointir}{\it #2}\par\vskip 3pt}
\long\def\tha#1|#2\fintha{\par\vskip 5pt
               {\petcap #1.}\par\nobreak{\it #2}\par\vskip 3pt}

\def\rem#1|
{\par\vskip 5pt
                {{\it #1}\pointir}}
\def\rema#1|
{\par\vskip 5pt
                {{\it #1.}\par\nobreak }}

\def\proof{\noindent {\it Proof. }}

\def\qed{\quad\raise -2pt \hbox{\vrule\vbox to 10pt
{\hrule width 4pt
\vfill\hrule}\vrule}}

\def\cqfd{\ifmmode
\unkern\quad\hfill\qed
\else
\unskip\quad\hfill\qed\bigskip
\fi}

\newcount\n
\def\exo{\advance\n by 1 \par \vskip .3cm {\bf \the \n}. }

\def\ad{\mathop{\rm ad}\nolimits}

\def\Ad{\mathop{\rm Ad}\nolimits}

\def\demi{{\textstyle{1\over 2}}\,}

\font\msbmten=msbm10
\def\ltimes{\mathbin{\hbox{\msbmten\char'156}}}


\hfill

\vskip 20 pt

\centerline{MINIMAL REPRESENTATIONS OF SIMPLE REAL LIE GROUPS  }

\vskip 2pt

\centerline{ OF NON HERMITIAN TYPE  }

\vskip 20 pt
 
\centerline{Dehbia Achab}

\vskip 4 pt

\centerline{Institut de Math\'ematiques de Jussieu}

\vskip 4 pt

\centerline{Universit\'e Pierre et Marie Curie}

\vskip 4 pt

\centerline{4 place Jussieu, case 247, 75252 Paris cedex 05}

\vskip 4 pt

\centerline{\tt dehbia.achab@imj-prg.fr}

\vskip 10 pt

\vskip 15pt

\noindent
{\bf Abstract} 
\it In the recent paper {\rm [AF12]}, we introduced an analysis of the Brylinski-Kostant model for  spherical minimal representations for simple real Lie groups of non Hermitian type.  We generalize here that analysis and give a unified geometric realization to  a family of unitary irreducible representations of such groups.  \rm
\vskip 10 pt
\noindent
{\it Key words}: Minimal representation, Lie algebra, Jordan algebra, Conformal group, Bernstein identity,  hypergeometric function, Meijer $G$-function.
\vskip 10 pt
\noindent
{\it Mathematics Subject Index 2010}:17C36, 22E46, 32M15, 33C80.
\vskip 20 pt
\noindent
Introduction
\vskip 4 pt
\noindent
1.  Construction process of  complex simple Lie algebras and   real  forms of non Hermitian type.
\vskip 4 pt
\noindent
2.  $\prod\limits_{i=1}^s(K_i)_{\bboard R}$-invariant Hilbert subspaces 
of ${{\cal O}}(\prod\limits_{i=1}^s\Xi_i)$.
\vskip 4 pt
\noindent
3.  Irreducible representations  of the Lie algebra.
\vskip 4 pt
\noindent
4. Unitary  representations of the   corresponding real Lie group.
\vskip 4 pt
\noindent
5. The ${\goth sl}(3,{\bboard R})$-case.
\vskip 4 pt
\noindent
References

\vskip 15cm
\vfill \eject
\section Introduction|

  In the papers [A11] and [AF12] the following principles were established:
  \vskip 1 pt

  1) Given a finite dimensional complex vector space  $V$ and  a homogeneous polynomial $Q$ on $V$ such that  the associated structure group  ${\rm Str}(V,Q)$  has an open symmetric orbit, one can equip $V$ with the structure of a  semisimple complex Jordan algebra (see Theorem 1.1 in [AF12]).
  \vskip 1 pt

  2)  There exists a finite covering $K$ of the conformal group ${\rm Conf}(V,Q)$ and  a cocycle $\mu:K\times V \rightarrow \bboard C$ such that  the corresponding cocycle representation of $K$ on the polynomial functions on $V$ leaves the space ${\cal W}$, spanned by $Q$ and its translates $z\mapsto Q(z-a)$ with $a\in V$, invariant, producing an irreducible representation $\kappa$ of $K$ on ${\cal W}$ (see Proposition 1.1 and Corollary 1.2 in [A11]).
  \vskip 1 pt

    3) Assuming  the polynomial $Q$ of degree 4, ${\rm Lie}(K)\oplus {\cal W}$ carries the structure of a complex simple Lie algebra $\goth g$
  (see Theorem 3.1 and Theorem 3.3 in [A11]). One constructs  a non compact real form ${\goth g}_{\bboard R}$ of  $\goth g$,  of non Hermitian type, starting  with a Euclidean  real form $V_{\bboard R}$ of $V$.  The construction also yields a compact real form $K_{\bboard R}$ of the group $K$. It turns out that we obtain in this way  (see table 1 in [A11]) all the simple  real Lie algebras of non Hermitian type which possess a strongly minimal real nilpotent orbit, i.e. such that  ${\cal O}_{min}\cap  \goth g_{\bboard R}\ne\emptyset$, where ${{\cal O}}_{min}$ is the minimal nilpotent adjoint orbit of $\goth g$ (see tables 4.7 and 4.8 in [B98]). Recall that from a point of view of representation theory, the condition that ${{\cal O}}_{min}\cap  \goth g_{\bboard R}\ne\emptyset$ is natural, since it is a necessary condition for a simple real Lie group with Lie algebra $\goth g_{\bboard R}$ to admit an irreducible unitary representation  with associated complex nilpotent orbit ${\cal O}_{min}$ (see Theorem 8.4 in [V91]).
\vskip 1 pt

4)  The $K$-orbit $\Xi=K.Q \subset {\cal W}$ carries the structure of a complex line bundle over the conformal compactification  $\bar V$ of $V$. Its restriction  $\Xi_0$ to $V$ is trivial and 
the $K$-action on $\Xi$ yields a left regular representation $\pi$ of $K$ on the space ${\cal O}(\Xi)$ of holomorphic functions  on $\Xi$.
\vskip 1 pt

5)  The subspaces ${\cal O}_m(\Xi)$ of fixed degree of homogeneity $m$ in the fiber direction  are irreducible $K$-modules whose elements restrict to polynomials on $\Xi_0$. Homogeneity  in the fiber direction allows to identify ${\cal O}_m(\Xi)$ with  a space $\tilde{\cal O}_m(V)$ of polynomials on $V$.
\vskip 1 pt
6) There exists a $K_{\bboard R}$-invariant  inner product on  $\tilde{\cal O}_m(V)$ which has a reproducing kernel  $H^m$,  where $H$ is obtained from $Q$ by polarization.
Adding these kernels with arbitrarily chosen non-negative weights yields a multiplicity free
unitary  $K_{\bboard R}$-representation on a Hilbert subspace of  ${\cal O}(\Xi)$.
\vskip 1 pt
7)  The polynomials $p\in {\cal W}$ act  on ${\cal O}_{\rm fin}(\Xi)=\sum\limits_{m\in \bboard N}{\cal O}_m(\Xi)$ by multiplication and differentiation  such that the corresponding maps ${\cal M},  {\cal D} : {\cal W} \rightarrow  $End$({\cal O}_{\rm fin}(\Xi))$ are $K$-equivariant.
\vskip 1 pt
8)  $E=Q$ and $F=1$ in ${\cal W}$ can be complemented by $H=[E,F] \in {\goth g}$ to an 
$sl_2$-triple which allows to determine a specific condition  (T) depending on the pair $(V,Q)$ under which, there exists a unique sequence  $\delta=(\delta_m)_{m\in \bboard N}$ acting diagonally on  ${\cal O}_{\rm fin}(\Xi)=\sum\limits_{m\in \bboard N}{\cal O}_m(\Xi)$  such that the map 
$\rho:  {\cal W} \rightarrow  $End$({\cal O}_{\rm fin}(\Xi)), p\mapsto {\cal M}(p)-\delta\circ{\cal D}(p)$ complements $d\pi$ to a representation $d\pi+\rho$  of ${\goth g}={\rm Lie}(K)+{\cal W}$ on  ${\cal O}_{\rm fin}(\Xi)$.
\vskip 1pt
Since ${\cal O}_0(\Xi)$ is reduced to the constants and then consists of $K$-fixed vectors, this representation is spherical.
\vskip 1 pt
A list of the cases satisfying condition (T) is given and that determine the non Hermitian simple real Lie groups  ${\goth g}_{\bboard R}$  which admit spherical minimal representations.
\vskip 1 pt
9)  When the condition (T) is satisfied, given the sequence $\delta$,  it is possible to  determine a set of weights such that  the restriction of d$\pi+\rho$ to ${\goth g}_{\bboard R}$ is infinitesimally unitary. Moreover, Nelson's criterion can be used to  show that this representation integrates to  the simply connected Lie group ${G_{\bboard R}}$ with Lie algebra ${\goth g}_{\bboard R}$.  
\vskip 6 pt
As the  geometric setting in [AF12]  allowed   to  give  realizations of only spherical minimal representations, we consider in this paper a more general geometric setting which allows to realize a large class of  minimal representations. We resolve  this problem by applying principles  1)-7), unless 3), separately for the simple summands $V_i$ of the Jordan algebra $V$ which means that $Q$ gets factored  as a product of powers of Jordan determinants. Then,  (neglecting coverings) K is a product of  the analogous groups $K_i$ for the simple summands $V_i$ and on the level  of $K$ representations, one has  a tensor product situation.  
The variety $\Xi$ is a quotient of   $\widetilde\Xi:=\prod\limits_{i=1}^s\Xi_i$, where  $s$ is the number of simple summands,  $\Xi_i=K_i\cdot \Delta_i$ is the $K_i$-orbit  of  the Jordan determinant $\Delta_i$.  This results in an ${\bboard N}^s$-grading rather than an ${\bboard N}$-grading, a fact which only slightly complicates the calculations. The technical  heart of this paper then is to make principle 8) work.   To this end,  one has to  consider holomorphic functions on  $\widetilde\Xi$ instead of $\Xi$,   and  replace ${\cal O}_{\rm fin}(\Xi)$ by spaces 
${\cal O}_{q,\rm fin}(\widetilde\Xi)=\sum\limits_{m\in {\bboard N}}=\otimes_{i=1}^s{\cal O}_{k_im+q_i}(\Xi_i)$ where  $(k_1,\ldots,k_s)$ are the multiplicities given by
$Q(z)=\prod\limits_{i=1}^s\Delta_i^{k_i}(z_i)$ and $q=(q_1,\ldots,q_s)\in {\bboard N}^s$.  Then the strategy of proof given in [AF12] can be used to prove principle 8) also in the product setup. The extra freedom of the parameter $q$ turns out to be  enough to produce  all the representations 
constructed by Brylinski-Kostant. The representation  so obtained    is realized on a Hilbert space ${\cal F}_q(\widetilde\Xi)$  of holomorphic functions on $\widetilde\Xi $.  
\hfill \eject
There is an explicit formula for the  reproducing kernel of ${\cal F}_q(\widetilde\Xi)$
involving a hypergeometric function  ${}_1F_2$.
\vskip 2 pt
The space ${\cal F}_q(\widetilde\Xi)$ is a weighted Bergman space. The norm is given by an integral on $\prod\limits_{i=1}^s({\bboard C}\times V_i)$,  with a weight
taking in general both positive and negative values, involving  a Meijer $G$-function   $G_{1,3}^{3,0}$. This  result permits to get an analogy with the classical Fock space. The paper [AF12] corresponds to the case $q=0.$
\vskip 1 pt
A similar  theory can be developped for   Fock models of minimal representations of simple real Lie groups of Hermitian type. It  is the subject of  another paper. 
\vskip 1 pt
The recent work by T. Kobayashi, B. Orsted, J. Hilgert and J. M\"ollers (see  [HKMO12]),  consists in constructing Schr\"odinger models for minimal representations of  real  Lie groups  of Hermitian type,   which arise as   conformal groups of simple real  Jordan algebras.  
Also, in a  series of papers,  A.  Dvorsky and S. Sahi  obtained explicit Hilbert spaces for  unipotent representations  of conformal groups of simple real Jordan  algebras (see [S92], [DS99], [DS03]).

\vskip 30 cm
\section 1. Construction process of  complex simple Lie algebras and of   simple real Lie algebras of non Hermitian type|

Let $V$ be a finite dimensional complex vector space and $Q$ a 
homogeneous polynomial on $V$. Define
$$ L={\rm Str}(V,Q)=\{g\in {\rm GL}(V) \mid \exists \gamma=\gamma(g),
Q(g\cdot x)=\gamma(g)Q(x)\}.$$
Assume that there exists $e \in V$ such that 
\vskip 1pt
\quad (1)  The symmetric bilinear form $\langle x,y \rangle =-D_xD_y\log Q(e),$
is non-degenerate.
\vskip 1pt
\quad (2) The orbit $\Omega= L\cdot e$ is open.
\vskip 1pt
\quad (3) The orbit $\Omega=L\cdot e$ is symmetric, i.e.  the pair
$( L, L_0)$, with $ L_0=\{g \in  L \mid g\cdot e=e\},$
is symmetric, which means that there is an involutive 
automorphism $\nu $ of $L$ 
such that $L_0$ is open in $\{g\in L \mid \nu (g)=g\}$.
\vskip 2 pt
\noindent
The  vector space $V$ is then equipped naturally with a product which makes it  a semisimple Jordan algebra  (see [AF12] Theorem 1.1 and [FK94]  for   Jordan algebras theory).
\vskip 2 pt
 
The conformal group ${\rm Conf}(V,Q)$ is the group of rational transformations $g$ of $V$ generated by: the translations $z\mapsto z+a$ ($a\in V$), the dilations $z\mapsto \ell \cdot z$
($\ell \in L$), and the inversion $j:z\mapsto -z^{-1}$ ( see [M78]). A transformation $g\in {\rm Conf}(V,Q)$
is conformal in the sense that the differential $Dg(z)$ belongs to $ {\rm Str}(V,Q)$ at any point $z$ where $g$ is defined.

Let $\cal W$ be the space of polynomials on $V$ 
generated by the translated $Q(z-a)$ of $Q$, with $a\in V$.  Let $\kappa$ be the cocycle representation of  ${\rm Conf}(V,Q)$ or of a covering of order two of it on $\cal W$,  defined in [A11] and [AF12] as follows:
\noindent
\vskip 1 pt
{\it Case 1.} In case there exists a character $\chi$ of Str$(V,Q)$ such that
$\chi^2=\gamma$, then let $K={\rm Conf}(V,Q)$. Define the cocycle 
$$\mu(g,z)=\chi((Dg(z)^{-1})\quad (g\in K,\ z\in V),$$
and the representation $\kappa$ of $K$ on $\cal W$,
$$(\kappa(g)p)(z)=\mu(g^{-1},z)p(g^{-1}\cdot z).$$
The function $\kappa(g)p$ belongs actually to $\cal W$ (see [FG96], Proposition 6.2).
The cocycle $\mu(g,z)$ is a polynomial in $z$ of degree 
$\leq \, {\rm deg}\, Q$ and
$$\eqalign{
(\kappa(\tau_a)p)(z)&=p(z-a) \quad (a\in V),\cr
(\kappa(\ell )p)(z)&=\chi(\ell )p(\ell^{-1}\cdot z) \quad (\ell\in L),\cr
(\kappa(j)p)(z)&=Q(z)p(-z^{-1}).\cr}$$

\medskip

\noindent
{\it Case 2.} Otherwise the group $K$ is defined as the set of pairs
$(g,\mu)$ with $g\in  {\rm Conf}(V,Q)$, and $\mu$ is a rational
function on $V$ such that 
$$\mu(z)^2=\gamma(Dg(z))^{-1}.$$
We consider on $K$ the product  $(g_1,\mu_1)(g_2,\mu_2)=(g_1g_2,\mu_3)$
with $\mu_3(z)=\mu_1(g_2\cdot z)\mu_2(z)$.
For $\tilde g=(g,\mu) \in K$, define $\mu(\tilde g,z):=\mu(z)$.
Then $\mu(\tilde g,z)$ is a cocycle:
$\mu(\tilde g_1\tilde g_2,z)=\mu(\tilde g_1,\tilde
g_2\cdot z)\mu(\tilde g_2,z),$
where $\tilde g\cdot z=g\cdot z$ by definition.
\vskip 2 pt
Recall  from [AF12], section 1, that the representation $\kappa$ of $K$ on $\cal W$ is irreducible and 
$$\bigl(\kappa(\tilde g)p\bigr)(z)=\mu({\tilde g}^{-1},z)p({\tilde g}^{-1}\cdot z).$$
\vskip 2 pt
 Observe that the inverse in $K$ of $\sigma=(j,Q(z))$ is
$\sigma^{-1}=(j,Q(-z))$ and then equals $\sigma$ if the degree of $Q$ is even.
From now on, assume that the polynomial $Q$ has degree 4.
Then there exists $H\in{\goth z}({\goth l})$, where ${\goth l}=$Lie$(L)$ which defines  gradings of  ${\goth k}=$Lie$(K)$ and of ${\goth p}={\cal W}$:
$${\goth k}={\goth k}_{-1}+{\goth k}_0+{\goth k}_1,\quad {\goth p}={\goth p}_{-2}+{\goth p}_{-1}+{\goth p}_0+{\goth p}_1+{\goth p}_2,$$
with 

\centerline{${\goth k}_j=\{X\in {\goth k}\mid \ad (H)X=jX\}, \quad (j=-1,0,1),$}
 
\centerline{${\goth k}_{-1}\simeq V,\quad {\goth k}_0=$Lie$(L),\quad
\Ad (\sigma): {\goth k}_j\rightarrow {\goth k}_{-j},$}

\vskip 2 pt
\noindent
and where   
$${\goth p}_j=\{p\in {\goth p} \mid {\rm d}\kappa(H)p=jp\}$$
is the set of polynomials in $\goth p$,  homogeneous  of degree 
$j+2$.  Furthermore 
$\kappa(\sigma) : {\goth p}_j \rightarrow
{\goth p}_{-j}$, and
$${\goth p}_{-2}={\bboard C}, \quad {\goth p}_2={\bboard C}\, Q, 
\quad {\goth p}_{-1} \simeq V,\quad {\goth p}_1\simeq V.$$
Then for   ${\goth g}={\goth k}\oplus {\goth p}$, $E=Q$, $F=1$,  we established (see Theorem 3.1 in [A11]) the existence on ${\goth g}$ of a unique simple Lie algebra structure such that:
$$\eqalign{
&(i)  \quad [X,X']=[X,X']_{\goth k} \quad (X,X' \in {\goth k}), \cr
&(ii) \quad [X,p]={\rm d}\kappa(X)p \quad (X\in {\goth k}, p\in {\goth p}), \cr
&(iii) \quad [E,F]=H. \cr}$$

\medskip
Recall  now the  real form ${\goth g}_{\bboard R}$ of $\goth g$ which will be considered in the sequel. It has been introduced in [A11] and [AF12] .  We fix a  Euclidean real form $V_{\bboard R}$ of the complex Jordan algebra $V$, denote by $z\mapsto
\bar z$ the conjugation of $V$ with respect to $V_{\bboard R}$, and then 
consider  the involution $g\mapsto \bar g$ of ${\rm Conf}(V,Q)$ given  by: 
$\bar g\cdot z=\overline{g\cdot \bar z}$.  
For $(g,\mu)\in K$ define
$$\overline{(g,\mu)}=(\bar g,\bar\mu),\  {\rm where}\ 
\bar\mu(z)=\overline{\mu(\bar z)}.$$
The involution
$\alpha$ defined by $\alpha (g)=\sigma\circ\bar g\circ\sigma^{-1}$ is a Cartan
involution of $K$ (see Proposition 1.1. in  [P02]),  and  $K_{\bboard R}=\{g\in K\mid \alpha(g)=g\}$
is a compact real form of $K$.
 \vskip 1 pt
 \noindent
Let $\goth u$ be the  compact real form of $\goth g$ such that 
${\goth k}\cap {\goth u}={\goth k}_{\bboard R}$, the Lie algebra of $K_{\bboard R}$.  Denote by $\goth p_{\bboard R}=\goth p\cap(i{\goth u})$. Then, the real Lie algebra  defined by
\vskip 2 pt

\centerline{$\goth g_{\bboard R}=\goth k_{\bboard R}+\goth p_{\bboard R}\quad \quad (*)$}

\vskip 2 pt
\noindent
 is a real form of $\goth g$ and the decomposition $(*)$ is its Cartan decomposition.
 Looking at the subalgebra ${\goth g}^0$ isomorphic to ${\goth sl}(2,{\bboard C})$ generated by the triple $(E,F,H)$, one sees that  $H\in i{\goth k}_{\bboard R}, E+F\in {\goth p}_{\bboard R}$ and $i(E-F)\in {\goth p}_{\bboard R}$. It follows that  ${\goth p}_{\bboard R}={\cal U}({\goth k}_{\bboard R})(E+F)+{\cal U}({\goth k}_{\bboard R})(i(E-F)).$
  \vskip 1 pt
  
Since the complexification of the Cartan decomposition of ${\goth g}_{\bboard R}$ is $\goth g=\goth k+\goth p$ and since $\goth p$ is a simple $\goth k$-module, it follows that the simple real Lie algebra ${\goth g}_{\bboard  R}$ is of non Hermitian type.
\vskip 2 pt
In [AF12], we have established that 
${\goth p}_{\bboard R}=\{p\in {\goth p}\mid \beta (p)=p\}$, 
where  we defined for a polynomial  $p\in {\goth p}$, $\bar p=\overline{p(\bar z)},$ and 
 considered the antilinear involution $\beta $ of $\goth p$ given by
$\beta (p)=\kappa (\sigma )\bar p.$
 \vskip 6 pt
\noindent
 Observe that $E+F$ and $i(E-F)$ belong to ${\cal W}_{\bboard R}$, and that 
 
  \vskip 4pt
 
 \centerline{${\cal W}_{\bboard R}={\rm d}\kappa({\cal U}({\goth k}_{\bboard R}))(E+F)+{\rm d}\kappa({\cal U}({\goth k}_{\bboard R}))(i(E-F))$,}

 \vskip 4 pt
 \noindent
Since $iH$ belongs to ${\goth k}_{\bboard R}$ and 
 
\vskip 4 pt

\centerline{$[E+F,i(E-F)]=-2i[E,F]=-2iH$,}

\vskip 4pt

\centerline{$[iH,E+F]=2i(E-F), \quad [iH,i(E-F)]=-2(E+F), $}

\vskip 4pt
\noindent
then  the real Lie subalgebra of ${\goth g}_{\bboard R}$   generated by $iH, E+F, i(E-F)$ is isomorphic to ${\goth su}(1,1)$.
\vskip 6 pt
\noindent
 
Another real form of the Lie algebra ${\goth g}$ can be obtained naturally by considering the real analogous of ${\goth g}$, given by 
$\tilde{\goth g}_{\bboard R}=\tilde{\goth k}_{\bboard R}+\tilde{\goth p}_{\bboard R} $, where $\tilde{\goth k}_{\bboard R}={\rm Lie}({\rm Str}(V_{\bboard R}))$,   $\tilde{\goth p}_{\bboard R}=\tilde{\cal W}_{\bboard R}$ is the real subspace of ${\goth p}$ generated by the polynomials  $Q(x-a)$ on $V_{\bboard R}$  with $a\in V_{\bboard R}$ in such a way that for 

\vskip 2 pt

\centerline{$\tilde{\goth p}_j=\{p\in {\tilde{\goth p}}_{\bboard R} \mid {\rm d}\kappa(H)p=jp\}$,}

\vskip 1 pt

\centerline{$\tilde{\goth p}_{\bboard R}=\tilde{\goth p}_{-2}\oplus\tilde{\goth p}_{-1}\oplus\ldots\oplus\tilde{\goth p}_0\oplus\ldots\oplus\tilde{\goth p}_{1}\oplus\tilde{\goth p}_{2}$}

\vskip 2 pt
\noindent
 is  the  eigen-space decomposition of $\tilde{\goth p}_{\bboard R}$ under ${\rm ad}(\tilde H)=d\kappa(\tilde H)$, where we notice that the element $H$ of ${\goth z}({\goth k})$ belongs to ${\goth z}(\tilde{\goth k}_{\bboard R})$ too.

\vskip 6 pt
\noindent
Observe that $E$ and $F$ belong to $\tilde{\cal W}_{\bboard R}$, and that 

\centerline{$\tilde{\cal W}_{\bboard R}={\rm d}\kappa({\cal U}(\tilde{\goth k}_{\bboard R}))E+{\rm d}\kappa({\cal U}(\tilde{\goth k}_{\bboard R}))F$}

\vskip 4 pt
\noindent
Since $H$ belongs to $\tilde {\goth k}_{\bboard R}$  and since 
 
\vskip 2 pt

\centerline{$[E,F]=H, [H,E]=2E, [H,F]=-2F, $}

\vskip 2 pt
\noindent
then  the real Lie subalgebra of $\tilde{\goth g}_{\bboard R}$   generated by $H, E, F$  is isomorphic to ${\goth sl}(2,\bboard R)$.

 \vskip 4 pt
 \noindent
The table in next page gives the classification of the simple real Lie algebras ${\goth g}_{\bboard R}$ obtained in this way.  It turns out that they are exactly the simple real Lie algebras of non Hermitian type which  satisfy the condition ${\cal O}_{min}\cap{ \goth g}_{\bboard}\ne\emptyset$ (see tables 4.7 and 4.8 in  [B98] ).  It also gives the real Lie algebras $\tilde{\goth g}_{\bboard R}$. 
\vskip 1 pt
We have used the notation: 
$$\varphi _p(z)=z_1^2+\cdots +z_p^2,\quad (z\in {\bboard C}^p).$$
In case of an exceptional Lie algebra $\goth g$, the real form ${\goth g}_{\bboard R}$ has been identified by computing the Cartan signature.

\vfill \eject

\centerline{\bf Table 1} \rm

\hsize=14cm
\def\tvi{\vrule height 12pt depth 1pt width 0pt}

\def\traithorizontal{\noalign{\hrule}}
$$\eqalignno{&\vbox{\offinterlineskip \halign{
\tvi # & # & # & # & #&# \cr
 \qquad  \qquad \quad  $V$ & $Q$ & $\goth k$ & $\goth g$ & ${\goth g}_{\bboard R}$ & $\tilde{\goth g}_{\bboard R}$  \quad \cr
\tvi \cr
\traithorizontal 
\qquad  \qquad \quad${\bboard C}$ & $z^4$ & ${\goth sl}(2,{\bboard C})$ \hfill 
& ${\goth sl}(3,{\bboard C})$ 
& ${\goth sl}(3,{\bboard R})$ & ${\goth sl}(3,{\bboard R})$ \quad\cr
\tvi \cr
\traithorizontal 
\qquad  \qquad \quad${\bboard C}^p$ & $\varphi _p(z)^2$ & ${\goth so}(p+2,{\bboard C})$ \hfill 
& ${\goth sl}(p+2,{\bboard C})$ 
& ${\goth sl}(p+2,{\bboard R})$ & ${\goth sl}(p+2,{\bboard R})$ \quad \cr
\tvi \cr
\traithorizontal
\qquad  \qquad \quad${\bboard C}\oplus {\bboard C}$ \hfill & $z_1^2z_2^2$
& ${\goth so}(3,{\bboard C})\oplus {\goth so}(3,{\bboard C})$
& ${\goth so}(6,{\bboard C})$ & ${\goth so}(3,3)$& ${\goth so}(6)$ \quad  \cr
\tvi \cr
\traithorizontal
\qquad   \qquad \quad ${\bboard C}\oplus {\bboard C}\oplus{\bboard C}$ \hfill & $z_1^2z_2z_3$
& ${\goth sl}(2,{\bboard C})^{\oplus 3}$
& ${\goth so}(7,{\bboard C})$ & ${\goth so}(3,4)$& ${\goth so}(7)$ \quad \cr
\tvi \cr
\traithorizontal
\qquad  \qquad \quad ${\bboard C}\oplus {\bboard C}\oplus{\bboard C}\oplus{\bboard C}$ \hfill & $z_1z_2z_3z_4$
& ${\goth sl}(2,{\bboard C})^{\oplus 4}$
& ${\goth so}(8,{\bboard C})$ & ${\goth so}(4,4)$ &${\goth so}(8)$ \quad \cr
\tvi \cr
\traithorizontal
\qquad   \qquad \quad ${\bboard C}^p\oplus {\bboard C}$ \hfill & $\varphi _p(z){z'}^2$
& ${\goth so}(p+2,{\bboard C})\oplus {\goth so}(3,{\bboard C})$
& ${\goth so}(p+5,{\bboard C})$ & ${\goth so}(p+2,3)$ & ${\goth so}(p+5)$ \quad\cr
\tvi \cr
\traithorizontal
\qquad   \qquad \quad ${\bboard C}^p\oplus {\bboard C}\oplus{\bboard C}$ \hfill & $\varphi _p(z)z'z''$
& ${\goth so}(p+2,{\bboard C})\oplus {\goth sl}(2,{\bboard C}^{\oplus 2})$
& ${\goth so}(p+6,{\bboard C})$ & ${\goth so}(p+2,4)$ & ${\goth so}(p+6)$ \quad\cr
\tvi \cr
\traithorizontal
\qquad  \qquad \quad ${\bboard C}^{p_1}\oplus {\bboard C}^{p_2}$ \hfill & $\varphi _{p_1}(z)\varphi _{p_2}(z')$
& $\oplus_{i=1}^2{\goth so}(p_i+2,{\bboard C})$
& ${\goth so}(p_1+p_2+4,{\bboard C})$ & ${\goth so}(p_1+2,p_2+2)$ & ${\goth so}(p_1+p_2+4)$  \quad \cr
\tvi \cr
\traithorizontal
\qquad  \qquad  \quad Sym$(4,{\bboard C})$ & $\det z$ & ${\goth sp}(8,{\bboard C})$ & ${\goth e}_6$
& ${\goth e}_{6(6)}$ & ${\goth e}_{6(6)}$  \quad\cr
\tvi \cr
\qquad  \qquad \quad ${\rm M}(4,{\bboard C})$ & $\det z$ & ${\goth sl}(8,{\bboard C})$ 
& ${\goth e}_7$ & ${\goth e}_{7(7)}$ & ${\goth e}_{7(7)}$ \quad\cr
\tvi \cr
\qquad   \qquad \quad Skew$(8,{\bboard C})$ & ${\rm Pfaff}(z)$ & ${\goth so}(16,{\bboard C})$ & ${\goth e}_8$
& ${\goth e}_{8(8)}$& ${\goth e}_{8(8)}$ \quad \cr
\tvi \cr
\traithorizontal
\qquad   \qquad \quad Sym$(3,{\bboard C})\oplus {\bboard C}$ & $\det z \cdot z'$ & 
${\goth sp}(6,{\bboard C})\oplus {\goth sl}(2,{\bboard C})$ \hfill 
& ${\goth f}_4$ & ${\goth f}_{4(4)}$ & ${\goth f}_{4(4)}$  \quad \cr
\tvi \cr
\qquad   \qquad \quad $M(3,{\bboard C})\oplus {\bboard C}$ \hfill & $\det z\cdot z' $ & 
${\goth sl}(6,{\bboard C})\oplus {\goth sl}(2,{\bboard C})$ \hfill &
${\goth e}_6$ & ${\goth e}_{6(2)}$ & ${\goth e}_{6(6)}$ \quad \cr
\tvi \cr
\qquad   \qquad \quad Skew$(6,{\bboard C})\oplus {\bboard C}$ & ${\rm Pfaff}( z) \cdot z' $ 
& ${\goth so}(12,{\bboard C})\oplus {\goth sl}(2,{\bboard C})$ \hfill & 
${\goth e}_7$ & ${\goth e}_{7(-5)}$ & ${\goth e}_{7(7)}$ \quad\cr
\tvi \cr
\qquad  \qquad \quad  Herm$(3,{\bboard O})_{\bboard C}\oplus {\bboard C}$ & $\det z\cdot z'$ &
${\goth e}_7\oplus {\goth sl}(2,{\bboard C})$ \hfill & ${\goth e}_8$ & ${\goth e}_{8(-24)}$ & ${\goth e}_{8(8)}$ \quad \cr
\tvi \cr
\traithorizontal
\qquad   \qquad \quad ${\bboard C}\oplus {\bboard C}$ \hfill & $z^3\cdot z'$ & 
${\goth sl}(2,{\bboard C})\oplus {\goth sl}(2,{\bboard C})$ \hfill & ${\goth g}_2$ & ${\goth g}_{2(2)}$ & ${\goth g}_{2(2)}$  \quad \cr
\tvi \cr
\traithorizontal
}}\cr}$$

\vfill \eject

\hsize=11,25cm

The Jordan algebra $V$ is a direct sum of simple ideals: 
$$V=\bigoplus_{i=1}^sV_i,$$
and
$$Q(z)=\prod _{i=1}^s\Delta_i(z_i)^{k_i}\quad
(z=(z_1,\ldots,z_s)),$$
where $\Delta_i$ is the
determinant polynomial of the simple Jordan algebra 
$V_i$ and the $k_i$ are positive integers. The degree of $Q$ is equal
to $\sum _{i=1}^sk_ir_i$, where $r_i$ is the rank of $V_i$.
Define
$$L_i={\rm Str}(V_i,\Delta_i)=\{g\in GL(V_i) \mid \exists\gamma_i(g)\in {\bboard C}^*,
\Delta_i(g\cdot x)=\gamma_i(g)\Delta_i(x)\}.$$

In a similar manner, we consider the conformal groups ${\rm Conf}(V_i,\Delta_i)$ ,  the corresponding groups  $K_i$ analogous of $K$,  the spaces  
${\goth p}_i$  of polynomials on $V_i$ 
generated by the polynomials  $\Delta_i(z_i-a_i)$ with $a_i\in V_i$  and the cocycle representations  
$\kappa_i$ of $K_i$, analogous of $\kappa$,  on the space  ${\goth p}_i$  given by $\bigl(\kappa_i(g_i)p\bigr)(z_i)=\mu_i(g_i^{-1},z_i)p(g_i^{-1}\cdot z_i).$
More precisely, $K_i={\rm Conf}(V_i,\Delta_i)$ when there exists a character $\chi_i$ of $L_i$ satisfying $\chi_i^{2}=\gamma_i$ and $K_i=\{\tilde g_i=(g_i,\mu_i) \mid g_i\in {\rm Conf}(V_i,\Delta_i), \mu_i(z_i)^{2}=\gamma_i(Dg_i(z_i))^{-1}\}$ if the character $\chi_i$ doesn't exist.
\vskip 1 pt
\noindent
Since  the differential  of the inversion $j : V_i\rightarrow V_i, z_i\mapsto -z_i^{-1}$ is given by $Dj(z_i)=P_i(z_i)^{-1}$ where $P_i$ is the quadratic representation of  $V_i$ and since  $P_i(z_i)$ belongs to $L_i$ and satisfies $\Delta_i(P_i(z_i)x_i)=\Delta_i^2(z_i)\Delta_i(x_i)$ then $\gamma_i(P_i(z_i))=\Delta_i^{2}(z_i)$ and the element $\sigma_i=(j,\Delta_i(z_i))$ belongs to $K_i$. The inverse in $ K_i$ of $\sigma_i=(j,\Delta_i(z_i))$ is
$\sigma_i^{-1}=(j,\Delta_i(-z_i))=(j,(-1)^{r_i}\Delta_i(z_i))$.
\vskip 2 pt
\noindent
Furthermore, when $\chi_i$ exists, the representation $\kappa_i$ of $K_i$, is given by

\centerline{$(\kappa_i(g_i)p)(z_i)=\mu_i(g_i^{-1},z_i)p(g_i^{-1}\cdot z_i)$}

\vskip 2 pt
\noindent
where $\mu_i(g_i,z_i)=\chi_i((Dg_i(z_i)^{-1})$. In particular 
$$\eqalign{
(\kappa_i(\tau_{a_i})p)(z_i)&=p(z_i-a_i) \quad (a_i\in V_i),\cr
(\kappa_i(\ell_i )p)(z_i)&=\chi_i(\ell_i )p(\ell_i^{-1}\cdot z_i) \quad (\ell_i\in L_i),\cr
(\kappa_i(j)p)(z_i)&=\Delta_i(z_i)p(-z_i^{-1}).\cr}.$$

\medskip
\noindent
And, when $\chi_i$ doesn't exist,  the representation $\kappa_i$ is given  by

\centerline{$\bigl(\kappa_i(\tilde g_i)p\bigr)(z_i)=\mu_i({\tilde g_i}^{-1},z_i)p({\tilde g_i}^{-1}\cdot z_i)$}

\vskip 2 pt
\noindent
where $\tilde g_i=(g_i,\mu_i(z_i))$ and $\mu_i(\tilde g_i,z_i)=\mu_i(z_i)$. In particular

\centerline{$(\kappa_i(\sigma_i)p)(z_i)=\Delta_i(z_i)p(-z_i^{-1})$ and $(\kappa_i(\sigma_i^{-1})p)(z_i)=(-1)^{r_i}\Delta_i(z_i)p(-z_i^{-1})$.}

\bigskip
Let ${\cal W}_i^{(k_i)}$ be  the  vector  space generated by the polynomials  $\Delta_i^{k_i}(z_i-a_i)$ for $a_i\in V_i$. Then, the  group $K_i$ acts on ${\cal W}_i^{(k_i)}$ by the representation  $\kappa_i^{(k_i)}$ given by
$$(\kappa_i^{(k_i)}(g_i)p)(z_i)=\mu_i(g_i,z_i)^{k_i}p(g_i^{-1}\cdot z_i).$$

\vskip 6 pt
\noindent
Observe that   $\prod\limits_{i=1}^s L_i $ acting on $V$ by  $g\cdot z=(g_i\cdot z_i)$, for  $g=(g_i)\in \prod\limits_{i=1}^sL_i$, and $z=(z_i) \in V$, is a subgroup of $L$,  that they have  the  same Lie algebra  $\goth l=\sum\limits_{i=1}^s{\goth l}_i$ with ${\goth l}_i=$Lie$( L_i)$,  and that 

\centerline{$\gamma(g)=\prod\limits_{i=1}^s\gamma_i^{k_i}(g_i)$}

\vskip 2 pt
\noindent
Also,  the groups $\prod\limits_{i=1}^sK_i $  and  $K$ have the same Lie algebra  $\goth k=\sum\limits_{i=1}^s{\goth k}_i$ with 

\vskip 2 pt

\centerline{${\goth k}_i=$Lie$(K_i)\simeq V_i\oplus{\goth l}_i\oplus V_i$.}

\vskip 1 pt

\th Proposition 1.1| The space ${\goth p}$ is  the tensor product of the ${\cal W}_i^{(k_i)}$.
\finth
\vskip 2 pt
\proof In fact,  let  $\psi$ be the multilinear map 
 ${\goth p}_1\times\ldots\times{\goth p}_s \rightarrow {\goth p},  (p_1,\ldots,p_s)\mapsto p$, with 
$p(z)=p_1(z_1)\ldots p_s(z_s)$  for $z=\sum\limits_{i=1}^sz_i$. 
Since $\kappa(\tau_a)Q(z)=\prod\limits_{i=1}^s(\kappa_i^{(k_i)}(\tau_{a_i})\Delta_i^{k_i})(z_i)$,  then the complex vector space  ${\goth p}$ is generated by $\psi({\cal W}_1^{(k_1)}\times\ldots\times{\cal W}_s^{(k_s)})$. 
It remains to prove that for a multilinear map $f : {\cal W}_1^{(k_1)}\times\ldots\times{\cal W}_s^{(k_s)} \rightarrow U$, there is a linear map $\tilde f : {\goth p}\rightarrow U$ such that $\tilde f\circ\psi=f$. 
\vskip 2 pt
\noindent
Recall that  ${\goth p}={\goth p}_{-2}+{\goth p}_{-1}+{\goth p}_0+{\goth p}_1+{\goth p}_2$,   ${\goth k}={\goth k}_{-1}+{\goth k}_0+{\goth k}_1$  with
${\goth k}_{-1} \simeq V\simeq {\goth k}_1=\kappa(\sigma)({\goth k}_{-1})$,  ${\goth p}_{-2}={\bboard C}\cdot 1, {\goth p}_2={\bboard C}\cdot Q$  and that  for  every $ p\in {\goth p}_1, q\in {\goth p}_{-1}, $ there is  a unique $X \in {\goth k}_{-1}$ and a unique   $Y\in 
{\goth k}_{-1}$ such that  
$p=d\kappa(X)Q$ and  $q=\kappa(\sigma)d\kappa(Y)Q$ (see  [A11], Lemma1.1 ). 
\vskip 2 pt

\vskip 2 pt
\noindent
Let's  define the linear map $\tilde f$ on ${\goth p}$  as follows:
\vskip 2 pt
\noindent
- On ${\goth p}_{2}$ and ${\goth p}_{-2}$:  $\tilde f(Q)=f(\Delta_1^{k_1},\ldots,\Delta_s^{k_s})$ ; $\tilde f(1)=f(1,\ldots,1)$.

\vskip 2 pt
\noindent
- On ${\goth p}_{1}$:  for $p\in {\goth p}_{1}$, let $X\in {\goth k}_{-1}$ be such that $p=[X,Q]$, then put 

\vskip 2 pt

\centerline{$\tilde f(p)=\tilde f(d\kappa(X)Q)=f(d\kappa_1^{(k_1)}(X_1)\Delta_1^{k_1},\ldots,d\kappa_s^{(k_s)}(X_s)\Delta_s^{k_s})$}

\vskip 2 pt
\noindent
where exp$(X)=\tau_a$,    $a=\sum\limits_{i=1}^sa_i$ and  $X_i \in {\goth k}_i$ such that exp$(X_i)=\tau_{a_i}$. 
\vskip 2 pt
\noindent
- On ${\goth p}_{-1}$:  for $q\in  {\goth p}_{-1}$, let $Y\in {\goth k}_{-1}$ such that $q=\kappa(\sigma)[Y,Q]$., then put 

\vskip 2 pt

\centerline{$\tilde f(q)=f(\kappa_1^{(k_1)}(\sigma_1){\rm d}\kappa_1^{(k_1)}(Y_1)\Delta_1^{k_1},\ldots,\kappa_s^{(k_s)}(\sigma_s){\rm d} \kappa_s^{(k_s)}(Y_s)\Delta_s^{k_s})$.}

\vskip 2 pt
\noindent
where   exp$(Y)=\tau_b$,  $b=\sum\limits_{i=1}^sb_i$ and $Y_i\in {\goth k}_i$ such that exp$(Y_i)=\tau_{b_i}$. 
\vskip 2 pt
\noindent
- On ${\goth p}_{0}$: for $\tilde p\in {\goth p}_0$, there is $(X,p)\in {\goth k}_{-1}\times{\goth p}_1$ such that $\tilde p=[X,p]$, and, since there is   $Y\in {\goth k}_{-1}$ such that $p=[Y,Q]$ then $\tilde p=[X,[Y,Q]]$.  Let $\{p_1,\ldots,p_{n_0}\}$ be a basis of the space ${\goth p}_0$ where $p_j=[X^{(j)},[Y^{(j)},Q]]$ with $X^{(j)}, Y^{(j)} \in {\goth k}_{-1}$.  Denote by 
$a^{(j)}=\sum\limits_{i=1}^sa_i^{(j)},   b^{(j)}=\sum\limits_{i=1}^sb_i^{(j)}$  the elements such that exp$(X^{(j)})=\tau_{a^{(j)}}$ and exp$(Y^{(j)})=\tau_{b^{(j)}}$ and let  $X_i^{(j)}, Y_i^{(j)} \in {\goth k}_i$  be  such that  exp$(X_i^{(j)})=\tau_{a_i^{(j)}}$ and exp$(Y_i^{(j)})=\tau_{b_i^{(j)}}$. 
\vskip 6 pt
\noindent
We define $\tilde f$ on this basis by 

\vskip 6 pt

\centerline{$\tilde f(p_j)=f([X_1^{(j)},[Y_1^{(j)},\Delta_1^{k_1}]]\prod\limits_{i=2}^s\Delta_i^{k_i},\ldots,[X_s^{(j)},[Y_s^{(j)},\Delta_s^{k_s}]]\prod\limits_{i=1}^{s-1}\Delta_i^{k_i})$.}

\vskip 6 pt
\noindent
Since 
$$({\rm d}\kappa(X)Q)(z)=({\rm d}\kappa_1^{(k_1)}(X_1)\Delta_1^{k_1})(z_1)\ldots ({\rm d}\kappa_s^{(k_s)}(X_s)\Delta_s^{k_s})(z_s),$$
$$(\kappa(\sigma){\rm d}\kappa(Y)Q)(z)=$$
$$(\kappa_1^{(k_1)}(\sigma_1){\rm d}\kappa_1^{(k_1)}(Y_1)\Delta_1^{k_1})(z_1)\ldots (\kappa_s^{(k_s)}(\sigma_s){\rm d}\kappa_s^{(k_s)}(Y_s)\Delta_s^{k_s})(z_s),$$

\vskip 2 pt
\noindent
and 
$$([X^{(j)},[Y^{(j)},Q]])(z)=$$
$$\sum\limits_{l=1}^s([X_l^{(j)},[Y_l^{(j)},\Delta_l^{k_l}]])(z_l)\prod\limits_{i\ne l}^s\Delta_i^{k_i}(z_i),$$

\vskip 2 pt
\noindent
then  

\centerline{$\tilde f\circ\psi=f$.}

\hfill \qed

\vskip 4 pt
Let  $\tilde\kappa=\kappa_1^{(k_1)}\otimes\ldots\otimes\kappa_s^{(k_s)}$ be  the  tensor product of the representations $\kappa_i^{(k_i)}$. 
The  infinitesimal representations $d\tilde\kappa$ and $d\kappa$ of ${\goth k}$ on ${\goth p}$ are equal 
and  the structure of simple Lie algebra  described above on ${\goth g}={\goth k}+{\goth p}$ can be obtained by considering the action $\tilde\kappa$ of the direct product group  $\prod\limits_{i=1}^sK_i $  (instead of the action $\kappa$ of $K$)  on  the tensor product space ${\goth p}=\otimes_{i=1}^s{\cal W}_i^{(k_i)}$.
\vskip 6 pt
The orbit $\Xi$ of $Q$ under the action of the group  $\widetilde K$, has dimension ${\rm dim}(V)+1$ and  is a conical variety : for $\xi\in \Xi, \lambda\in {\bboard C}^*$, $\lambda \xi$ belongs to $\Xi$.

\vskip 2 pt
\hfill \eject

Let $\Xi_i$  be the
$K_i$-orbit  of $\Delta_i$ in ${\cal W}_i$ under $\kappa_i$:
$$\Xi_i=\{\kappa_i(g_i)\Delta_i\mid g_i\in  K_i\}.$$
Then $\Xi_i$ is a conical variety.   It can  be seen as a line bundle
over the conformal compactification of $V_i$. A polynomial $ \xi_i\in {\cal W}_i$ can be written
$$\xi_i(v)=w_i\Delta_i(v)+\ {\rm terms \  of \ degree}\ < r_i \quad (w_i\in {\bboard C}),$$
and $w_i=w_i(\nu_i)$ is a linear form on ${\cal W}_i$ which is semi-invariant under the preimage in $K_i$ of the  maximal 
parabolic subgroup $P_{\rm max}^{(i)}=L_i\ltimes N_i$, where $N_i$ is the group of translations $z_i\in V_i\mapsto z_i+a_i$, for $a_i\in V_i$.
The set $\Xi_{i,0}=\{\xi_i \in \Xi_i \mid w_i(\xi_i)\ne 0\}$
 is open and dense in $\Xi_i$.
 A polynomial $\xi_i\in \Xi_{i,0}$  can be written $\xi_i(v_i)=w_i\Delta_i(v_i-z_i) \quad
(w_i \in {\bboard C}^*, z_i\in V_i).$
Hence  we get a coordinate system $(w_i,z_i)\in {\bboard C}^*\times V_i$ for $\Xi_{i,0}$.
\vskip 4 pt
\noindent
In this coordinate  system, the cocycle action 
of $ K_i$ is given by 
$$\kappa_i(g_i) : (w_i,z_i) \mapsto \bigl(\mu_i(g_i,z_i)w_i,g_i\cdot z_i\bigr).$$

\vskip 4 pt
Denote by  $\Xi_i^{(k_i)}$   the
$K_i$-orbit  of $\Delta_i^{k_i}$ in ${\cal W}_i^{(k_i)}$ under $\kappa_i^{(k_i)}$:
$$\Xi_i^{(k_i)}=\{\kappa_i^{(k_i)}(g_i)\Delta_i^{k_i}\mid g_i\in  K_i\}.$$
Since $\kappa_i^{(k_i)}(g_i)\Delta_i^{k_i}=(\kappa_i(g_i)\Delta_i)^{k_i}$, then 
$\Xi_i^{(k_i)}=\{\xi_i^{k_i} \mid \xi_i\in \Xi_i\}.$
Furthermore,  a polynomial $\xi_i^{(k_i)} \in {\cal W}_i^{(k_i)}$ can be  written 
$$\xi_i^{(k_i)}(v)=\tilde w_i\Delta_i^{k_i}(v)+\ {\rm terms \  of \ degree}\ < k_ir_i \quad (\tilde w_i\in {\bboard C}),$$
 the set $\Xi_{i,0}^{(k_i)}=\{\xi_i^{(k_i)} \in \Xi_i^{(k_i)}\mid \tilde w_i(\xi_i^{(k_i)})\ne 0\}$
 is open and dense in $\Xi_i^{(k_i)}$
and a polynomial $\xi_i^{(k_i)}\in \Xi_{i,0}^{(k_i)}$  can be written $\xi_i^{(k_i)}(v_i)=\tilde w_i\Delta_i^{k_i}(v_i-z_i) \quad
(\tilde w_i \in {\bboard C}^*, z_i\in V_i).$
Hence  we get a coordinate system $(\tilde w_i,z_i)\in {\bboard C}^*\times V_i$ for $\Xi_{i,0}^{(k_i)}$.
\vskip 4 pt
\noindent
In this coordinate  system, the cocycle action 
of $ K_i$  is given by 
$$\kappa_i^{(k_i)}(g_i) : (\tilde w_i,z_i) \mapsto \bigl(\mu_i^{k_i}(g_i,z_i)\tilde w_i,g_i\cdot z_i\bigr).$$

\vskip 4 pt
Define on $\widetilde\Xi=\prod\limits_{i=1}^s\Xi_i$ the equivalence relation 
\vskip 4 pt

\centerline{$(\xi_1,\ldots,\xi_s) \sim (\lambda_1\xi_1,\ldots,\lambda_s\xi_s)$,  for $ (\lambda_1,\ldots,\lambda_s)\in {\bboard C}^s$,  $\prod_{i=1}^s\lambda_i^{k_i}=1$.}

\vskip 4 pt
\noindent
\th Proposition 1.2|  The  map   $\widetilde\Xi \rightarrow \Xi,  (\xi_1,\ldots,\xi_s) \mapsto \xi_1^{k_1}\ldots\xi_s^{k_s}$ is a covering.

\bigskip
\noindent
The group $K_i$  acts on the space  ${{\cal O}}(\Xi_i)$  of holomorphic functions on $\Xi_i$  by: 

\vskip 2 pt

\centerline{$\bigl(\pi_i(g_i) f_i\bigr)({\xi_i})=f_i\bigl(\kappa_i(g_i^{-1})\xi _i\bigr)$.}
\vskip 2 pt
\noindent
If  $\xi_i(v_i)=w_i\Delta_i(v_i-z_i)$,  and $f_i\in {{\cal O}}(\Xi_i)$, we denote by
$ f_i(\xi_i)=\phi_i(w_i,z_i)$  the restriction of $f_i$ to $\Xi_{i,0}$. Then  $\phi_i\in {\cal O}({\bboard C}^*\times V_i)$. 
In the coordinates $(w_i,z_i)$, the  representation $\pi_i$ is given by
$$\eqalignno{(\pi_i(g_i)\phi_i)(w_i,z_i)&=\phi_i(\mu_i(g_i^{-1},z_i)w_i,g_i^{-1}\cdot z_i)\cr}$$

\vskip 2 pt
\noindent
 The group $K_i$  acts on the space  ${{\cal O}}(\Xi_i^{(k_i)})$   by: 
 \vskip 2 pt
 
 \centerline{$\bigl(\pi_i^{(k_i)}(g_i) f_i^{(k_i)}\bigr)({\xi_i^{(k_i)}})=f_i^{(k_i)}\bigl(\kappa_i^{(k_i)}(g_i^{-1})\xi _i^{(k_i)}\bigr)$.}

\vskip 4 pt
\noindent
If  $\xi_i^{(k_i)}(v_i)=\tilde w_i\Delta_i^{k_i}(v_i-z_i)$,  and $f_i\in {{\cal O}}(\Xi_i^{(k_i)})$, we denote by
$ f_i^{(k_i)}(\xi_i^{(k_i)})=\phi_i^{(k_i)}(\tilde w_i,z_i)$  the restriction of $f_i^{(k_i)}$ to $\Xi_{i,0}^{(k_i)}$. Then $\phi_i^{(k_i)}\in {\cal O}({\bboard C}^*\times V_i)$. 
In the coordinates $(\tilde w_i,z_i)$, the  representation $\pi_i^{(k_i)}$ is given by
$$\eqalignno{(\pi_i^{(k_i)}(g_i)\phi_i^{(k_i)})(\tilde w_i,z_i)&=\phi_i^{(k_i)}(\mu_i^{k_i}(g_i^{-1},z_i)\tilde w_i,g_i^{-1}\cdot z_i)\cr}$$
\vskip 2 pt
\noindent
Observe that for every $f_i^{(k_i)}\in {{\cal O}}(\Xi_i^{(k_i)})$, one can associate a unique $f_i\in {{\cal O}}(\Xi_i)$ such that $f_i(\xi_i)=f_i^{(k_i)}({\xi_i}^{k_i})$ in such a way that 
$$\bigl(\pi_i^{(k_i)}(g_i) f_i^{(k_i)}\bigr)({\xi_i}^{k_i})=\bigl(\pi_i(g_i)f_i\bigr)(\xi_i).$$
\vskip 4 pt
\noindent
The group $\widetilde K$  acts on the tensor product  space $\otimes_{i=1}^s{{\cal O}}(\Xi_i)$   by the tensor product representation $\pi=\otimes_{i=1}^s\pi_i$  given by:
$$\bigl(\pi(g)(f_1\otimes\ldots\otimes f_s)\bigr)(\xi_1,\ldots,\xi_s)=\bigl(\pi_1(g_1)f_1\bigr)\otimes\ldots\otimes\bigl(\pi_s(g_s)f_s\bigr)(\xi_1,\ldots,\xi_s)$$
$$\qquad \qquad \qquad \qquad \qquad  \quad  \quad =f_1(\kappa_1(g_1^{-1})\xi _1)\ldots f_s(\kappa_s(g_s^{-1})\xi _s)\quad (g=(g_i))$$
In the coordinates $(w_i,z_i)$, the  representation $\pi$ is given by
$$\eqalignno{(\pi(g)(\phi_1\otimes\ldots\otimes\phi_s))((w_1,z_1),\ldots,(w_s,z_s))&=\prod\limits_{i=1}^s\phi_i(\mu_i(g_i^{-1},z_i)w_i,g_i^{-1}\cdot z_i).\cr}$$

\vskip 4pt

\finth
\hfill
\eject

\vskip 4 pt
\noindent
\section 2. $\prod\limits_{i=1}^s({K_i})_{\bboard R}$-invariant Hilbert subspaces 
of $\otimes_{i=1}^s{{\cal O}}(\Xi_i)$|
  \vskip 2 pt
For  every  $m_i, p_i \in{\bboard Z}$,  we  denote by ${{\cal O}}_{p_i}(\Xi_i^{(k_i)})$ and ${{\cal O}}_{m_i}(\Xi_i)$ the spaces of holomorphic functions
$f_i^{(k_i)}$ on $\Xi_i^{(k_i)}$ and $f_i$ on $\Xi_i$  such that for every $\lambda\in {\bboard C}^*$,

\vskip 2 pt
\centerline{$f_i^{(k_i)}(\lambda\xi_i^{(k_i)})=\lambda^{p_i}f_i(\xi_i^{(k_i)})$ and $f_i(\lambda\xi_i)=\lambda^{m_i}f_i(\xi_i).$}

\vskip 2 pt

The space ${{\cal O}}_{m_i}(\Xi_i)$ (resp.  ${{\cal O}}_{p_i}(\Xi_i^{(k_i)})$)  is invariant under the representation
$\pi_i$ (resp. $\pi_i^{(k_i)}$).
\vskip 4 pt
The embedding $f_i^{(k_i)}\in {{\cal O}}_{p_i}(\Xi_i^{(k_i)})\mapsto  f_i\in  {{\cal O}}_{k_ip_i}(\Xi_i) $, where $f_i(\xi_i)=f_i^{(k_i)}(\xi_i^{k_i})$ intertwines the  restrictions $\pi_{i,p_i}^{(k_i)}$ and $\pi_{i,k_ip_i}$ of the representations $\pi_i^{(k_i)}$ and $\pi_i$ to ${{\cal O}}_{p_i}(\Xi_i^{(k_i)})$ and ${{\cal O}}_{k_ip_i}(\Xi_i)$.
\bigskip
From now on,  we essentially work with  the spaces ${{\cal O}}_{m_i}(\Xi_i)$. The representations 
$(\pi_{i,p_i}^{(k_i)}, {{\cal O}}_{p_i}(\Xi_i^{(k_i)}))$ and $(\pi_{i,k_ip_i}, {{\cal O}}_{k_ip_i}(\Xi_i))$ of $K_i$ are irreducible and equivalent.
\bigskip
If $f_i\in  {{\cal O}}_{m_i}(\Xi_i)$, then its 
restriction $\phi_i$ to $\Xi_{i,0}$ can be written 
$\phi_i(w_i,z_i)=w_i^{m_i}\psi_i(z_i)$
where $\psi_i$ is a holomorphic function on $V_i$. We will write 
$\widetilde {{\cal O}}_{m_i}(V_i)$ for the space of the functions $\psi_i$ corresponding to the
functions $f_i\in {{\cal O}}_{m_i}(\Xi_i)$, and denote by $\tilde\pi_{i,m_i}$ the representation of
$K_i$  on $\widetilde {{\cal O}}_{m_i}(V_i)$  corresponding to the restriction $\pi_{i,m_i}$ of
$\pi$ to ${{\cal O}}_{m_i}(\Xi_i)$. It is given by
$$(\tilde\pi_{i,m_i}(g_i)\psi_i)(z_i)=\mu_i(g_i^{-1},z_i)^{m_i}\psi_i(g_i^{-1}\cdot z_i).$$
In particular,
$$(\tilde\pi_{i,m_i}(\sigma_i)\psi_i)(z_i)=\Delta_i^{m_i}(z_i)\psi_i(-z_i^{-1}).$$

\bigskip

 \th Theorem 2.1|
\vskip 1 pt
 \noindent
{\rm (i)} ${{\cal O}}_{m_i}(\Xi_i)=\{0\}$ if $m_i<0$.
\vskip 1 pt
 \noindent
{\rm (ii)} The space ${{\cal O}}_{m_i}(\Xi_i) $ is  finite dimensional, and the
representation $\pi_{i,m_i}$ is  irreducible.
\vskip 1 pt
 \noindent
{\rm (iii)} The functions $\psi_i$ in $\widetilde{{\cal O}}_{m_i}(V_i)$  are polynomials with degree $\leq m_ir_i$.
\vskip 1 pt
 \noindent
{\rm (iv)}  $\widetilde{{\cal O}}_{m_i}(V_i)=\{\psi_i\in {\cal O}(V_i) \mid  \tilde\pi_{i,m_i}(\sigma_i)\psi_i \in {\cal O}(V_i)\}$.
\vskip 1 pt
 \noindent
{\rm (v)}  $ \widetilde{\cal O}_{m_i}(V_i) \subset  \widetilde{\cal O}_{m_i+1}(V_i)$. 

\finth

\vskip 1 pt
\noindent
The proof  of {\rm (i)}, {\rm (ii)}, {\rm (iii)} is similar to that of   Theorem 2.2 in [AF12] and is omitted. It remains to prove {\rm (iv)} and {\rm (v)}.
\vskip 1 pt
\noindent
{\rm (iv)} The inclusion $\widetilde{{\cal O}}_{m_i}(V_i) \subset \{\psi_i\in {\cal O}(V_i) \mid  \tilde\pi_{i,m_i}(\sigma_i)\psi_i \in {\cal O}(V_i)\}$ is obvious. Let $\psi_i\in {\cal O}(V_i)$  be such that   $\pi_{i,m_i}(\sigma_i)\psi_i \in {\cal O}(V_i)$. Then, $\psi_i$ is a polynomial.
\vskip 1 pt
\hfill
\eject
 Since $(\pi_{i,m_i}(\sigma_i)\psi_i)(z_i)=\Delta_i^{m_i}(z_i)\psi_i(-z_i^{-1})$, then the degree of $\psi_i$ is $\leq$ to the degree of $\Delta_i^{m_i}$ which is $m_ir_i$.  Moreover, since $\pi_{i,m_i}(\tau_{a_i})\psi_i$ and $\pi_{i,m_i}(l_i)\psi_i$ are both holomorphic on $V_i$ for every $a_i\in V_i$ and $l_i\in K_{i,0}$, and since the elements $\sigma_i, \tau_{a_i}$ and $l_i$ generate the group $K_i$,  it follows that  $\pi_{i,m_i}(g_i)\psi_i$ is holomorphic on $V_i$ for every $g_i\in K_i$. Then 
 the function $f_i$ defined on $\Xi_i$ by    $f_i(\xi_i)=w_i^{m_i}\psi_i(z_i)$  for $\xi_i(v)=w_i\Delta_i(v-z_i)$, belongs to ${\cal O}_{m_i}(\Xi_i)$ and then the function $\psi_i$ belongs to   $\widetilde{{\cal O}}_{m_i}(V_i)$.
\vskip 1 pt
\noindent
{\rm (v)} Let $\psi_i$ be an element of $ \widetilde{{\cal O}}_{m_i}(V_i)$. The  function $\pi_{i,m_i}(\sigma_i)\psi_i$ given by  $(\pi_{i,m_i}(\sigma_i)\psi_i)(z_i)=\Delta_i^{m_i}(z_i)\psi_i(-z_i^{-1})$ is  holomorphic  on $V_i$.  Then the function 
$\pi_{i,m_i+1}(\sigma_i)\psi_i( z_i)=\Delta_i(z_i)(\pi_{i,m_i}(\sigma_i)\psi_i )(z_i)$ is holomorphic on $V_i$, i.e. $\psi_i$ belongs to  the space $\widetilde{\cal O}_{m_i+1}(V_i)$.\quad \hfill \qed

\bigskip

\noindent
For    $(m_1,\ldots,m_s)\in {\bboard N}^s$,  we  denote by ${{\cal O}}_{(m_1,\ldots,m_s)}(\widetilde\Xi)=\otimes_{i=1}^s{\cal O}_{m_i}(\Xi_i)$.    It  is invariant under the tensor product  representation $\pi=\otimes_{i=1}^s\pi_i$. 
 \vskip 4 pt
 \noindent
For  $f=f_1\otimes\ldots\otimes f_s\in {{\cal O}}_{(m_1,\ldots,m_s)}(\widetilde\Xi)$,  and if $\phi_i(w_i,z_i)$ is the restriction of $f_i$ to $\Xi_{i,0}$,  then   $\phi_i(w_i,z_i)=w_i^{m_i}\psi_i(z_i)$ and

\vskip 6 pt

\centerline{$(\phi_1\otimes\ldots\otimes\phi_s)((w_1,z_1),\ldots,(w_s,z_s))=w_1^{m_1}\ldots w_s^{m_s}\psi_1(z_1)\ldots\psi_s(z_s).$}

\vskip 6 pt
\noindent
We write  $\widetilde {{\cal O}}_{(m_1,\ldots,m_s)}(V)=\otimes_{i=1}^s\widetilde{{\cal O}}_{m_i}(V_i)$ and  $\tilde\pi _{(m_1,\ldots,m_s)}=\otimes_{i=1}^s\tilde\pi_{m_i}$.  Then 
\vskip 4 pt
\noindent

\centerline{$\bigl(\tilde\pi_{(m_1,\ldots,m_s)}(g)(\psi_1\otimes\ldots\otimes\psi_s) \bigr)(z_1,\ldots,z_s)=\prod\limits_{i=1}^s{\mu_i(g_i^{-1},z_i)}^{m_i}\psi_i(g_i^{-1}\cdot z_i).$}

\vskip 2 pt
\noindent
For $\sigma=(\sigma_i)$, \quad $(\tilde\pi_{(m_1,\ldots,m_s)}(\sigma)1)(z_1,\ldots,z_s)=\Delta_1^{m_1}(z_1)\ldots. \Delta_s^{m_s}(z_s).$
\bigskip

We now  consider the  compact real forms $(K_i)_{\bboard R}$  analogous of $K_{\bboard R}$ for the groups $K_i$ and 
define on  the space ${{\cal O}}_{(m_1,\ldots,m_s)}(\widetilde\Xi)$ a $\prod\limits_{i=1}^s(K_i)_{{\bboard R}}$-invariant inner product. 
Define the subgroup $K_{i,0}$ of $K_i$ as $K_{i,0}=L_i$ in Case 1, and the preimage
of $L_i$ in Case 2, relatively to the covering map $K_i\to {\rm Conf}(V_i,\Delta_i)$,
and also $(K_{i,0})_{\bboard R}=K_{i,0} \cap (K_i)_{\bboard R}$. 
The coset space $M_i=(K_i)_{\bboard R}/(K{i,0})_{\bboard R}$, 
is a compact  Hermitian space and is the
conformal compactification of $V_i$. 
There is on $M_i$ a $(K_i)_{\bboard R}$-invariant probability measure, 
for which $M_i\setminus V_i$ has measure 0. Its restriction $m_{i,0}$ to $V_i$ is a probability measure with a density .
\vskip 4 pt
Let $H_i(z_i,z_i')$ be the polynomial on $V_i\times V_i$, holomorphic in $z_i$,
anti-holomorphic in $z_i'$ such that $H_i(x_i,x_i)=\Delta_i(e_i+x_i^2)$ for $x_i\in(V_i)_{\bboard R}.$ Put $H_i(z_i)=H_i(z_i,z_i)$. If $z_i$ is invertible, then
$H_i(z_i)=\Delta_i(\bar z_i)\Delta_i(\bar z_i^{-1}+z_i)$.
\vskip 4 pt

\th Proposition 2.2|
For $g_i\in (K_i)_{\bboard R}$,
\vskip 4 pt

\centerline{$H_i(g_i\cdot z_i,g_i\cdot z_i)'\mu_i(g_i,z_i)\overline{\mu_i(g_i,z_i')}=H_i(z_i,z_i'),$}

\vskip 2 pt

\centerline{$H_i(g_i\cdot z_i)\vert\mu_i(g_i,z_i)\vert^2=H_i(z_i).$}

\finth
\vskip 4 pt
We define the norms of  $\psi_i\in \widetilde{{\cal O}}_{m_i}(V_i)$  and  of $\psi\in   \widetilde{{\cal O}}_{(m_1,\ldots,m_s)}(V)$ by
$$\Vert \psi_i\Vert_{i,m_i}^2={1\over a_{i,m_i}}
\int_{V_i}\vert \psi_i(z_i)\vert^2H_i(z_i)^{-m_i}m_{i,0}(dz_i),$$
$$\Vert \psi\Vert_{(m_1,\ldots,m_s)}^2={1\over a_{(m_1,\ldots,m_s)}}
\int_{\prod\limits_{i=1}^sV_i}\vert \psi(z_1,\ldots,z_s)\vert^2\prod\limits_{i=1}^sH_i(z_i)^{-m_i}m_{i,0}(dz_i),$$
$$a_{i,m_i}=\int_{V_i}H_i(z_i)^{-m_i}m_{i,0}(dz_i), \quad a_{(m_1,\ldots,m_s)}=\prod\limits_{i=1}^sa_{i,m_i}.$$

\th Proposition 2.3| 
{ \rm (i)}  The  norm $\Vert\cdot\Vert_{i,m_i}$  is $(K_i)_{{\bboard R}}$-invariant. 
Hence,  $\widetilde{{\cal O}}_{m_i}(V_i)$ is a Hilbert subspace of
${{\cal O}}(V_i)$.  The reproducing kernel of $\widetilde{{\cal O}}_{m_i}(V_i)$ is given by
$$\tilde{\cal K}_{i,m_i}(z_i,z_i')= H_i(z_i,z_i')^{m_i} .$$
\vskip 1 pt
\noindent
{\rm (ii)} The norm $\Vert\cdot\Vert_{(m_1,\ldots,m_s)}$ is $\prod\limits_{i=1}^s(K_i)_{{\bboard R}}$-invariant.  Hence,  $\widetilde{{\cal O}}_{(m_1,\ldots,m_s)}(V)$ is a Hilbert subspace of
${{\cal O}}(\prod\limits_{i=1}^sV_i)$ with  reproducing kernel  
$$\tilde{\cal K}_{(m_1,\ldots,m_s)}((z_1,\ldots,z_s),(z_1',\ldots,z_s'))=\prod\limits_{i=1}^s\tilde{\cal K}_{i,m_i}(z_i,z_i')=\prod\limits_{i=1}^sH_i(z_i,z_i')^{m_i} .$$
\finth
\noindent
Since ${{\cal O}}_{m_i}(\Xi_i)$  is isomorphic to $\widetilde{{\cal O}}_{m_i}(V_i)$, then 
 ${{\cal O}}_{m_i}(\Xi_i)$ becomes an invariant Hilbert subspace of ${{\cal O}}(\Xi_i)$, with
reproducing kernel given by ${\cal K}_{i,m_i}(\xi_i,\xi'_i)=\Phi_i (\xi_i,\xi_i')^{m_i},$
where  $\Phi_i (\xi_i,\xi'_i)=H_i(z_i,z'_i)w_i\overline{w_i'}$. 
Hence  ${{\cal O}}_{(m_1,\ldots,m_s)}(\widetilde\Xi)$ is an invariant Hilbert subspace of ${{\cal O}}(\prod\limits_{i=1}^s\Xi_i)$ with
reproducing kernel 

\vskip 1 pt

\centerline{${\cal K}_{(m_1,\ldots,m_s)}((\xi_1,\ldots,\xi_s),(\xi_1',\ldots,\xi_s'))=\prod\limits_{i=1}^s{\cal K}_{i,m_i}(\xi_i,\xi'_i)=\prod\limits_{i=1}^s\Phi_i (\xi_i,\xi_i')^{m_i}.$}
\vskip 1 pt
\hfill
\eject

\th Theorem 2.4|
\vskip 2 pt
\noindent
{\rm (i)} The group $(K_i)_{{\bboard R}}$  acts
multiplicity free on the space ${{\cal O}}(\Xi_i)$.  
The irreducible $(K_i)_{{\bboard R}}$-invariant  subspaces of ${{\cal O}}(\Xi_i)$ 
are the spaces ${{\cal O}}_{m_i}(\Xi_i)$ with  $m_i\in {\bboard N}$. 
If ${\cal H}_i \subset {{\cal O}}(\Xi_i)$ is a 
$(K_i)_{{\bboard R}}$-invariant Hilbert subspace, the reproducing kernel
of ${\cal H}_i$ can be written 
$${\cal K}_i(\xi_i,\xi_i')=\sum _{m_i\in {\bboard N}}c_{i,m_i}\Phi_i (\xi_i,\xi_i')^{m_i},$$
where $(c_{i,m_i})$  is a sequence of positive numbers  such that   the series  $\sum\limits _{m_i\in {\bboard N}} c_{i,m_i}\Phi_i (\xi_i,\xi_i')^{m_i}$ converges uniformly on compact subsets in $\Xi_i$.
\vskip 2 pt
\noindent
{\rm (ii)} 
The group $\prod\limits_{i=1}^s(K_i)_{{\bboard R}}$  acts
multiplicity free on the space $\otimes_{i=1}^s{{\cal O}}(\Xi_i)$.  
The irreducible $\prod\limits_{i=1}^s(K_i)_{{\bboard R}}$-invariant  subspaces of $\otimes_{i=1}^s{{\cal O}}(\Xi_i)$ 
are the spaces ${{\cal O}}_{(m_1,\ldots,m_s)}(\widetilde\Xi)$ with  $(m_1,\ldots,m_s)\in {\bboard N}^s.$
The tensor product  space ${\cal H}=\otimes_{i=1}^s{\cal H}_i\subset \otimes_{i=1}^s{{\cal O}}(\Xi_i)$  of $(K_i)_{{\bboard R}}$-invariant Hilbert spaces  is a 
$\prod\limits_{i=1}^s(K_i)_{{\bboard R}}$-invariant Hilbert subspace and the reproducing kernel
of $\cal H$ can be written 
$${\cal K}((\xi_1,\ldots,\xi_s),(\xi_1',\ldots,\xi_s'))=\prod\limits_{i=1}^s{\cal K}_i(\xi_i,\xi_i').$$

\finth
\vskip 1 pt
\noindent
The proof  of {\rm (i)} is similar to that of   Theorem 2.5 in [AF12] and is omitted. The {\rm (ii)}  is then deduced from {(i)} for the tensor product situation.
\bigskip
In case of a weighted Bergman space there is an integral formula for the numbers  $ c_{i,p_i}$. 
For a positive function $\nu_i(\xi_i)$ on $\Xi_i$, 
consider the subspace ${\cal H}_i \subset {{\cal O}}(\Xi_i)$ 
of functions $\phi_i$ such that 
$$\Vert\phi_i\Vert_i^2=\int_{{\bboard C}\times V_i}\vert\phi_i(w_i,z_i)\vert^2\nu_i(w_i,z_i)m(dw_i)m_{i,0}(dz_i)\ <\infty,$$

\vskip 4pt
\noindent
and   then consider the subspace ${\cal H}=\otimes{\cal H}_i \subset \otimes_{i=1}^s{{\cal O}}(\Xi_i)$ 
of   functions $\phi$ such that  
$$\Vert\phi\Vert^2=$$
$$\int_{\prod\limits_{i=1}^s{\bboard C}\times V_i}\vert\phi((w_1,z_1),\ldots,(w_s,z_s))\vert^2\prod\limits_{i=1}^s\nu_i(w_i,z_i)m(dw_i)m_{i,0}(dz_i)\ <\infty.$$

\hfill
\eject

\vskip 1 pt
\th Theorem 2.5|
Let $F_i$ be a positive function on $[0,\infty[$, and define 

\vskip 2 pt

\centerline{$\nu_i(w_i,z_i)=F_i(H_i(z_i)\vert w_i\vert^{2})H_i(z_i).$}

\vskip 2 pt
\noindent
{\rm (i)} Then ${\cal H}_i$  is $(K_i)_{{\bboard R}}$-invariant.  

\vskip 3 pt
\noindent
{\rm (ii)} If $\phi_i(w_i,z_i)=\sum\limits _{m_i\in{\bboard N}}w_i^{m_i}\psi_{i,m_i}(z_i)$, then
\vskip 2 pt

\centerline{$\Vert
\phi_i\Vert_i^2=\sum\limits_{m_i \in{\bboard N}}{1\over c_{i,m_i}}\Vert\psi_{i,m_i}\Vert_{i,m_i}^2$}

\vskip 2 pt
\noindent
with 
$${1\over c_{i,m_i}}=\pi a_{i,m_i} \int_{[0,+\infty[}F_i(u_i)u_i^{m_i}du_i.$$
{\rm (iii)} The reproducing kernel of ${\cal H}_i$ is given by 
$$\eqalignno{&{\cal K}_i(\xi_i,\xi_i')=\sum\limits_{m_i\in{\bboard N}}c_{i,m_i}\Phi_i(\xi_i,\xi_i')^{m_i}.\cr}$$
\noindent
{\rm (iv)} The space ${\cal H}$  is $\prod\limits_{i=1}^s(K_i)_{{\bboard R}}$-invariant. If $\phi=\phi_1\otimes\ldots\otimes\phi_s$, then 
$$\Vert
\phi\Vert^2=\sum _{(m_1,\ldots,m_s)\in{\bboard N}^s}{1\over c_{(m_1,\ldots,m_s)}}\prod\limits_{i=1}^s\Vert\psi_{i,m_i}\Vert_{i,m_i}^2,$$
with $c_{(m_1,\ldots,m_s)}=\prod\limits_{i=1}^sc_{i,m_i}.$  The 
reproducing kernel  is given by 
$${\cal K}((\xi_1,\ldots,\xi_s),(\xi_1',\ldots,\xi_s'))=\prod\limits_{i=1}^s{\cal K}_i(\xi_i,\xi_i').$$
\finth
\vskip 2 pt
\noindent
The proof  is similar to that of  Theorem 2.6. in [AF12] and  {\rm (iv)} follows  by applying {\rm (i)}, {\rm (ii)} and {\rm (iii)}  to   the tensor product situation.
\vskip 2 pt
\th  Proposition 2.6|
\vskip 3 pt
\noindent
 {\rm (i)}The polynomial $\Delta_i^{k_i}$ satisfies the following    Bernstein identity 
$$\Delta_i^{k_i}\Bigl({\partial \over \partial z_i}\Bigr) \Delta_i^{k_i\alpha } =B_i(\alpha )\Delta_i^{k_i\alpha -k_i} 
\quad (z_i\in V_i),$$
\noindent
with
$$B_i(\alpha )=b_i(k_i\alpha )b_i(k_i\alpha -1)\ldots b_i(k_i\alpha -k_i+1),$$
\noindent
where  $b_i$ is the Bernstein polynomial relative to  the determinant $\Delta _i$.
\vskip 1 pt
\noindent
{\rm (ii)} Furthermore
$$\Delta_i^{k_i}\Bigl({\partial \over \partial z_i}\Bigr) H_i(z)^{k_i\alpha }
=B_i(\alpha )\overline{\Delta_i^{k_i}(z_i)}H_i(z_i)^{k_i(\alpha-1)}.$$
\finth
\vskip 2 pt
The measure $m_{i,0}$ has a density with respect to the Lebesgue measure $m(dz_i)$ on $V_i $: $m_{i,0}(dz_i)={1\over C_{i,0}}H_{i,0}(z_i) m(dz_i),$
with 
$$H_{i,0}(z_i)= H_i(z_i)^{-2{n_i\over r_i}}, \quad 
C_{i,0}=\int _{V_i} H_{i,0}(z_i)m(dz_i). $$
The Lebesgue measure $m(dz_i)$ will be chosen such that $C_{i,0}=1$.
\vskip 6 pt
\noindent
Recall that we have introduced the numbers
$$a_{i,m_i}=\int _{V_i}H_i(z_i)^{-m_i}m_{i,0}(dz_i),  \quad a_{(m_1,\ldots,m_s)}=\prod\limits_{i=1}^sa_{i,m_i}.$$
\th Proposition 2.7|
$$a_{i,m_i}= {\Gamma _{\Omega _i}(2{n_i\over r_i})
\over \Gamma _{\Omega _i}({n_i\over r_i})}
\prod _{i=1}^s {\Gamma _{\Omega _i}(m_i+{n_i\over r_i})
\over \Gamma _{\Omega _i}(m_i+2{n_i\over r_i})},$$
where $\Gamma _{\Omega _ i}$ is the Gindikin gamma function of the symmetric cone $\Omega _i$ in the Euclidean Jordan algebra $(V_i)_{\bboard R}$.
\finth
\vskip 1 pt
The proof is similar to that of  Proposition 3.2 in [AF12],  in the situation of $\Delta_i$ instead of $Q$.
\bigskip

\section 3. Irreducible representations of the Lie algebra|
\vskip 3 pt
In the sequel, we construct  irreducible   representations  $\rho_q$  of the Lie algebra $\goth g$, which will be the infinitesimal  versions of  minimal representations.  Recall that the group $\widetilde K=\prod\limits_{i=1}^sK_i$  acts on $\otimes_{i=1}^s{\cal O}(\Xi_i)$  by $\pi=\otimes_{i=1}^s\pi_i$, 
$$\bigl(\pi (g)f\bigr)(\xi_1,\ldots,\xi_s)=\prod\limits_{i=1}^sf_i\bigl(\kappa_i(g_i)^{-1}\xi _i\bigr), \quad  (g=(g_i), f=f_1\otimes\ldots\otimes f_s)$$
which is given in   the coordinates $(w,z)=((w_1,z_1),\ldots,(w_s,z_s))$ by
$$\eqalignno{\pi(g)\phi((w_1,z_1),\ldots,(w_s,z_s))&=\prod\limits_{i=1}^s\phi_i(\mu_i(g_i^{-1},z_i)w_i,g_i^{-1}\cdot z_i)\cr}$$
 This leads to a representation d$\pi$ of the Lie algebra $\goth k$ in  $\otimes_{i=1}^s{\cal O}(\Xi_i)$.
\vskip 4 pt
\hfill
\eject
Following the method of R. Brylinski and B. Kostant, we will
construct a representation
$\rho$ of
${\goth g}={\goth k}+{\goth  p}$ on  some subspace  of the space of finite sums 
$${{\cal O}}_{\rm fin}(\widetilde\Xi)=\sum _{(m_1,\ldots,m_s)\in  {\bboard N}^s}{{\cal O}}_{(m_1,\ldots,m_s)}(\widetilde\Xi),$$
such that,  

\centerline{$\forall X\in {\goth k}$, $\rho (X)=$d$\pi (X)$.}

\vskip 6 pt

We define first a representation  $\rho$ of the subalgebra generated by $E,F, H$,
isomorphic to ${\goth sl}(2,{\bboard C})$. In particular
$$\rho (H)={\rm d}\pi (H)={{\rm d}\over {\rm d}t}\Big| _{t=0}\pi (\exp tH).$$
Hence, for $\phi\in {{\cal O}}_{(m_1,\ldots,m_s)}(\widetilde\Xi)$,
$$\rho(H)\phi=({\cal E}-{(m_1r_1+\ldots+m_sr_s)\over 2})\phi$$
where $\cal E$ is the Euler operator
$${\cal E}\phi ((w_1,z_1),\ldots,(w_s,z_s))={{\rm d}\over {\rm d}t}\Big|_{t=0}\phi((w_1,e^tz_1),\ldots,(w_s,e^tz_s)).$$
The space  ${\cal O}_{(m_1,\ldots,m_s)}(\widetilde\Xi)$ can be identified  to the  space of functions  $\phi_{(m_1,\ldots,m_s)}((w_1,z_1),\ldots,(w_s,z_s))=\prod\limits_{i=1}^sw_i^{m_i}\psi_i(z_i)$ with $\psi_i\in
 \widetilde{{\cal O}}_{m_i}(V_i)$. Then, each $\phi\in {\cal O}_{\rm fin}(\widetilde\Xi)$ can be  written
$$\phi((w_1,z_1),\ldots,(w_s,z_s))=\sum\limits_{(m_1,\ldots,m_s)}\phi_{(m_1,\ldots,m_s)}((w_1,z_1),\ldots,(w_s,z_s)).$$
\bigskip
 Let ${\bboard D}(V)^L$ denote the algebra of $L$-invariant holomorphic differential operators on the open set of invertible elements of $V$.  One can show that this  algebra is isomorphic  to the algebra ${\bboard D}(\Omega)^G$ of $G$-invariant differential operators on $\Omega$, the symmetric cone  associated to the Euclidean real form $V_{\bboard R}$ of $V$, where $G$ is the connected component of the identity of the group $G(\Omega)=\{g\in GL(V) \mid g\Omega=\Omega\}$.  In [AF12] we established the following results:
 \hfill
 \eject
\vskip 1pt
1)  If $V$ is simple with rank $r$,   degree $d$ and dimension $n$,  and $Q=\Delta $, the determinant polynomial, then ${\bboard D}(V)^L$
is isomorphic to the algebra ${\cal P}({\bboard C}^r)^{{\goth S} _r}$ of symmetric polynomials in $r$ variables. The map
$$D\mapsto \tilde\gamma (D),\quad {\bboard D}(V)^L\to {\cal P}({\bboard C}^r)^{{\goth S}_r},$$
is the Harish-Chandra isomorphism , with  $\tilde\gamma (D)$ determined by 
$$  \tilde\gamma (D_{\alpha})(\lambda)=\prod\limits_{j=1}^r(\lambda_j-\alpha+{d\over 4}(r-1))$$ 
where  $D_{\alpha}=\Delta(z)^{1+\alpha}\Delta({\partial\over\partial z})\Delta(z)^{-\alpha}.$ (Recall that ${n\over r}=1+(r-1){d\over 2}$).
\vskip 1pt 
2) If  $V=\sum\limits_{i=1}^sV_i$ is semisimple and $Q=\prod\limits_{i=1}^s\Delta_i^{k_i}$,  ${\bboard D}(V)^L$ is isomorphic to the algebra
$\prod _{i=1}^s {\cal P}({\bboard C}^{r_i})^{{\goth S}_{r_i}}.$
The isomorphism is given by
$$D\mapsto \tilde\gamma (D)=\bigl(\tilde\gamma _1(D),\ldots ,\tilde\gamma _s(D)\bigr),$$
where $\tilde\gamma _i$ is the Harish-Chandra isomorphism relative to the algebra ${\bboard D}(V_i)^{L_i}$.
\vskip 2 pt
And, for $D\in {\bboard D}(V)^L$, we defined the adjoint $D^*$ by
$D^*=J\circ D\circ J$, where $Jf(z)=f\circ j (z)=f(-z^{-1})$  and established that $\tilde\gamma (D^*)(\lambda )=\tilde\gamma (D)(-\lambda )$.
\vskip 4 pt
In the present setting we define the Maass differential operator ${\bf D}_{\alpha }^{(i)}$ as

\centerline{${\bf D}_{\alpha }^{(i)}=\Delta_i(z_i) ^{k_i+\alpha }\Delta _i^{k_i}\Bigl({\partial \over \partial z_i}\Bigr)
\Delta _i(z_i)^{-\alpha}$}

\centerline{$
=\prod _{j=1}^{k_i} \Delta_i (z_i)^{\alpha +k_i-j+1}\Delta_i \Bigl({\partial \over \partial z_i}\Bigr)
\Delta_i (z_i)^{-(\alpha +k_i-j)}.$}

It is a $L_i$-invariant holomorphic differential operator on $V_i$. We write
$$\gamma _{\alpha }^{(i)} (\lambda^{(i)} )=\tilde\gamma_i({\bf D}_{\alpha }^{(i)})(\lambda^{(i)} ).$$
Then
$$\gamma _{\alpha }^{(i)} (\lambda^{(i)} )=\prod _{j=1}^r\big[\lambda _j^{(i)}-\alpha +\demi ({n_i\over r_i}-1)\big]_{k_i}$$
where $n_i$ is the dimension of $V_i$ and where we  have used the Pochhammer symbol $[a]_k=a(a-1)\ldots (a-k+1)$.
\vskip 4 pt
\noindent
Observe that the spaces $\widetilde{\cal O}_{m_i}(V_i)$ of  polynomial functions,  considered in Theorem 2.1, are stable under  the action of the Maass operator  ${\bf D}_{\alpha }^{(i)}$. 

 \th Lemma 3.1| For $\psi_i\in  \widetilde{\cal O}_{m_i}(V_i)$, the polynomial  $\Delta_i\Bigl({\partial\over \partial z_i}\Bigr)\psi_i$ belongs to  $ \widetilde{\cal O}_{m_i-1}(V_i)$.
 \finth
\vskip 1 pt
\proof
Let $\psi_i$ be an element of $ \widetilde{{\cal O}}_{m_i}(V_i)$.  We have to show that the function $(\pi_{i,m_i-1}(\sigma_i)(\Delta_i\Bigl({\partial\over \partial z_i}\Bigr)\psi_i))$ is holomorphic on $V_i$. In fact , this function is given by  
$$(\pi_{i,m_i-1}(\sigma_i)(\Delta_i\Bigl({\partial\over \partial z_i}\Bigr)\psi_i))(z_i)=\Delta_i^{(m_i-1)}(z_i)(\Delta_i\Bigl({\partial\over \partial z_i}\Bigr)\psi_i)(-z_i^{-1})=$$
$$\Delta_i^{m_i}(z_i)(\Delta_i\Delta_i\Bigl({\partial\over \partial z_i}\Bigr)\psi_i)(-z_i^{-1})=\Delta_i^{m_i}(z_i)({\bf D}_0^{(i)}\psi_i)(-z_i^{-1})$$
$$=\pi_{i,m_i}(\sigma_i)({\bf D}_0^{(i)}\psi_i)(z_i).$$ 
and since the  polynomial function ${\bf D}_0^{(i)}\psi_i$ belongs to  $ \widetilde{{\cal O}}_{m_i}(V_i)$, it follows that  $\pi_{i,m_i}(\sigma_i)({\bf D}_0^{(i)}\psi_i)(z_i)$ is holomorphic on $V_i$.

\hfill \qed
 
 \bigskip
 We introduce a  multiplication operator  ${\cal M}$ and a differential operator ${\cal D}$   on the space ${\cal O}_{\rm fin}(\widetilde\Xi)$ defined as follows : for  $\phi_{(m_1,\ldots,m_s)}$ in  ${\cal O}_{(m_1,\ldots,m_s)}(\widetilde\Xi)$, given by $\phi_{(m_1,\ldots,m_s)}((w_1,z_1),\ldots,(w_s,z_s))=\prod\limits_{i=1}^sw_i^{m_i}\psi_i(z_i)$,
$$({\cal M}\phi_{(m_1,\ldots,m_s)} )((w_1,z_1),\ldots,(w_s,z_s))=\prod\limits_{i=1}^sw_i^{m_i+k_i}\psi_i(z_i),$$
\noindent
and
$$({\cal D}\phi_{(m_1,\ldots,m_s)})((w_1,z_1),\ldots,(w_s,z_s))=\prod\limits_{i=1}^sw_i^{m_i-k_i}\biggl(\Delta_i^{k_i}\Bigl({\partial\over \partial z_i}\Bigr)\psi_i\biggr)(z_i).$$
\noindent
\vskip 2 pt
\noindent
Using  the item (v)  of Theorem 2.1,  and  Lemma 3.1,  one can see  that,  the operator ${\cal M}$ maps the space ${\cal O}_{(m_1,\ldots,m_s)}(\widetilde\Xi)$ into ${\cal O}_{(m_1+k_1,\ldots,m_s+k_s)}(\widetilde\Xi)$ and  that  the operator 
${\cal D}$ maps  ${\cal O}_{(m_1,\ldots,m_s)}(\widetilde\Xi)$ into ${\cal O}_{(m_1-k_1,\ldots,m_s-k_s)}(\widetilde\Xi)$ if every $m_i$ is  $\geq k_i$ and into $\{0\}$   if there is  $m_i$ such that   $m_i<k_i$. 

\bigskip
 We  denote by  ${\cal M}^{\sigma}$ and ${\cal D}^{\sigma}$ the conjugate operators: 
$$
{\cal M}^{\sigma}=\pi(\sigma){\cal M}\pi(\sigma)^{-1}, \quad
{\cal D}^{\sigma}=\pi(\sigma){\cal D}\pi(\sigma)^{-1}.$$
\noindent
Then
\vskip 1 pt
\hfill
\eject
\noindent
$$\eqalign{(&{\cal M}^{\sigma}\phi_{(m_1,\ldots,m_s)})(w,z)=(\pi(\sigma){\cal M}\pi(\sigma)\phi_{(m_1,\ldots,m_s)})(w,z)\cr
&=({\cal M}\pi(\sigma)\phi_{(m_1,\ldots,m_s)})((\Delta_1(z_1)w_1,-z_1^{-1}),\ldots,(\Delta_s(z_s)w_s,-z_s^{-1})),}$$
\noindent
and 
 $$\eqalign{&(\pi(\sigma)\phi_{(m_1,\ldots,m_s)})(w,z)=\cr
 &\phi_{(m_1,\ldots,m_s)}((\Delta_1(z_1)w_1,-z_1^{-1}),\ldots,(\Delta_s(z_s)w_s,-z_s^{-1}))\cr
&=\prod\limits_{i=1}^sw_i^{m_i}\Delta_i^{m_i}(z_i)\psi_i(-z_i^{-1})=\prod\limits_{i=1}^sw_i^{m_i}\Delta_i^{m_i}(z_i)(\psi_i\circ j_i)(z_i).}$$
\vskip 2 pt
\noindent
Then  
$$({\cal M}\pi(\sigma)\phi_{(m_1,\ldots,m_s)})(w,z)=\prod\limits_{i=1}^sw_i^{m_i+k_i}\Delta_i^{m_i}(z_i)(\psi_i\circ j_i)(z_i)$$
\noindent
and 
$$\eqalignno{&({\cal M}\pi(\sigma)\phi_{(m_1,\ldots,m_s)})((\Delta_1(z_1)w_1,-z_1^{-1}),\ldots,(\Delta_s(z_s)w_s,-z_s^{-1}))\cr}$$
$$\eqalignno{&=\Delta_i^{(m_i+k_i)}(z_i)\Delta_i^{m_i}(-z_i^{-1})(\psi_i\circ j_i)(-z_i^{-1}).\cr}$$
It follows that 
$$({\cal M}^{\sigma}\phi_{(m_1,\ldots,m_s)})((w_1,z_1),\ldots,(w_s,z_s))=\prod\limits_{i=1}^sw_i^{m_i+k_i}\Delta_i^{k_i}(z_i)\psi_i(z_i).$$

Furthermore, if there is $m_i$ such that $p_i<k_i$, then $({\cal D}^{\sigma}\phi_{(m_1,\ldots,m_s)})=0$, and if every $m_i$ is $\geq k_i$, then
\vskip 1pt
$({\cal D}^{\sigma}\phi_{(m_1,\ldots,m_s)})(w,z)=(\pi(\sigma){\cal D}\pi(\sigma)\phi_{(m_1,\ldots,m_s)})(w,z)$
$$=({\cal D}\pi(\sigma)\phi_{(m_1,\ldots,m_s)})((\Delta_1(z_1)w_1,-z_1^{-1}),\ldots,(\Delta_s(z_s)w_s,-z_s^{-1})),$$
$$\eqalign{&(\pi(\sigma)\phi_{(m_1,\ldots,m_s)})(w,z)=\cr
&\phi_{(m_1,\ldots,m_s)}((\Delta_1(z_1)w_1,-z_1^{-1}),\ldots,(\Delta_s(z_s)w_s,-z_s^{-1}))\cr
&=\prod\limits_{i=1}^sw_i^{m_i}\Delta_i^{m_i}(z_i)\psi_i(-z_i^{-1})=\prod\limits_{i=1}^sw_i^{m_i}\Delta_i^{m_i}(z_i)(\psi_i\circ j_i)(z_i),}$$
$$({\cal D}\pi(\sigma)\phi_{(m_1,\ldots,m_s)})(w,z)=\prod\limits_{i=1}^sw_i^{m_i-k_i}\Delta_i^{k_i}\Bigl({\partial\over \partial z_i}\Bigr)(\Delta_i^{m_i}(\psi_i\circ j_i))(z_i).$$
Then 
$$\eqalign{&({\cal D}^{\sigma}\phi_{(m_1,\ldots,m_s)})(w,z)=\cr
&\prod\limits_{i=1}^sw_i^{m_i-k_i}(\Delta_i^{k_i-m_i}\Delta_i^{k_i}\Bigl({\partial\over \partial z_i}\Bigr)\Delta_i^{m_i}(\psi_i\circ j_i))\circ j_i(z_i)\cr
&=\prod\limits_{i=1}^sw_i^{m_i-k_i}(({{\bf D}_{-m_i}^{(i)}})^*\psi_i)(z_i).}$$
\vskip 2 pt

Given a sequence $(\delta_{(m_1,\ldots,m_s)})_{(m_1,\ldots,m_s)\in {\bboard N}^s}$,  one defines  the diagonal operator 
$\delta$ on ${{\cal O}}_{\rm fin}(\widetilde\Xi)$   by  
$$\delta(\sum _{(m_1,\ldots,m_s)\in{\bboard N}^s} \phi_{(m_1,\ldots,m_s)})=\sum _{(m_1,\ldots,m_s)\in {\bboard N}^s} \delta_{(m_1,\ldots,m_s)}\phi_{(m_1,\ldots,m_s)},$$
$$\eqalignno{&
\rho(F)={\cal M}-\delta\circ {\cal D}, \quad 
\rho(E)=\pi(\sigma)\rho(F)\pi(\sigma)^{-1}
={\cal M}^{\sigma}-\delta\circ {\cal D}^{\sigma}.\cr}$$
\vskip 2 pt
\th Lemma 3.2| 
$$[\rho(H),\rho(E)]=2\rho(E),$$
$$[\rho(H),\rho(F)]=-2\rho(F).$$
\finth
\vskip 1 pt
\proof
Since
$\rho(H) {\cal M}$ and ${\cal M}\rho(H)$ map  $\prod\limits_{i=1}^sw_i^{m_i}\psi_i(z_i)$   to respectively  

\centerline{$ ({\cal E}-\sum\limits_{i=1}^s{(m_i+k_i)r_i\over 2})\prod\limits_{i=1}^sw_i^{m_i+k_i}\psi_i(z_i),$ and $ ({\cal E}-\sum\limits_{i=1}^s{m_ir_i\over 2})\prod\limits_{i=1}^sw_i^{m_i+k_i}\psi_i(z_i),$}

\noindent
one obtains $[\rho(H),{\cal M}]=-{\sum\limits_{i=1}^sk_ir_i\over 2}{\cal M}=-2{\cal M}.$
\vskip 4 pt
\noindent
- If  $m_i\geq k_i$, then
$\rho(H)\delta\circ {\cal D} $ and $\delta\circ{\cal D}\rho(H)$ map  $\prod\limits_{i=1}^sw_i^{m_i}\psi_i(z_i)$ to respectively

\centerline{$\delta_{(m_1-k_1,\ldots,m_s-k_s)}({\cal E}-\sum\limits_{i=1}^s{(m_i-k_i)r_i\over 2})(\prod\limits_{i=1}^s(\Delta_i^{k_i}\Bigl({\partial\over\partial z_i}\Bigr)\psi_i(z_i)w_i^{m_i-k_i}),$}

 \centerline{$\delta_{(m_1-k_1,\ldots,m_s-k_s)}(\prod\limits_{i=1}^s\Delta_i^{k_i}\Bigl({\partial\over \partial z_i}\Bigr))({\cal E}-\sum\limits_{i=1}^s{m_ir_i\over 2})(\prod\limits_{i=1}^s\psi_i(z_i)w_i^{m_i-k_i}),$}

\vskip 6pt
 \noindent
then, by using the identity 

\vskip 6 pt

\centerline{$[\prod\limits_{i=1}^s\Delta_i^{k_i}\Bigl({\partial\over \partial z_i}\Bigr),{\cal E}]=4\prod\limits_{i=1}^s\Delta_i^{k_i}\Bigl({\partial\over \partial z_i}\Bigr),$}

\vskip 2 pt
\noindent
one obtains
\vskip 4 pt

\centerline{$[\rho (H),\delta\circ{\cal D}]( \prod\limits_{i=1}^sw_i^{m_i}\psi_i(z_i))=
-2\delta_{(m_1-k_1,\ldots,m_s-k_s)}\prod\limits_{i=1}^s\Delta_i^{k_i}\Bigl({\partial \over \partial z_i}\Bigr)\psi_i(z_i)w_i^{m_i-k_i}$}

\vskip 2 pt
\noindent
i.e. 
\vskip 4 pt

\centerline{$[\rho(H),\rho(F)]\phi_{(m_1,\ldots,m_s)}=-2\rho(F)\phi_{(m_1,\ldots,m_s)}.$}

\vskip 2 pt
\noindent
- If there is $i$ such that $m_i< k_i$, then 
$$\rho(H)\delta\circ {\cal D} \phi_{(m_1,\ldots,m_s)}=0 \quad , \quad \delta\circ{\cal D}\rho(H)\phi_{(m_1,\ldots,m_s)}=0,$$

\vskip 2 pt

\centerline{$[\rho(H),\rho(F)]\phi_{(m_1,\ldots,m_s)}=[\rho(H),{\cal M}]\phi_{(m_1,\ldots,m_s)}=-2\rho(F)\phi_{(m_1,\ldots,m_s)}$.}

\vskip 2pt
\noindent
Finally, one gets 
$$[\rho(H),\rho(F)]=-2\rho(F).$$
Furthermore,  the operator $\delta$ commutes with $\pi(\sigma)$, and 
$\pi(\sigma)\rho(H)\pi(\sigma)^{-1}=-\rho(H)$,
we get also 
$[\rho(H),\rho(E)]=2\rho(E)$.
\hfill \qed
\vskip 4 pt
\noindent
\th Proposition 3.3| \rm 
The subspaces ${{\cal O}}_{(m_1,\ldots,m_s)}(\widetilde\Xi)$ are invariant under $[\rho (E),\rho (F)]$
and the restriction of  $[\rho(E),\rho(F)]$ to ${{\cal O}}_{(m_1,\ldots,m_s)}(\widetilde\Xi)$ 
commutes with the $(\prod\limits_{i=1}^sL_i)$-action:
\vskip 2 pt

\centerline{$[\rho(E),\rho(F)]: 
\psi(z)\prod\limits_{i=1}^sw_i^{m_i} \mapsto (P_{(m_1,\ldots,m_s)}\psi)(z)\prod\limits_{i=1}^sw_i^{m_i},$}

\vskip 2 pt
\noindent
 where $\psi(z)=\prod\limits_{i=1}^s\psi_i(z_i)$ and $ P_{(m_1,\ldots,m_s)}$ is an $\prod\limits_{i=1}^sL_i$-invariant holomorphic differential operator
on $V$  given as follows:
\vskip 1 pt
\noindent
- If every $m_i$ is $\geq k_i$, then
$$\eqalign{&P_{(m_1,\ldots,m_s)}=\cr
&\delta_{(m_1,\ldots,m_s)}(\prod\limits_{i=1}^s{\bf D}_{-1}^{(i)}-\prod\limits_{i=1}^s({\bf D}_{-m_i-k_i}^{(i)})^*)\cr
&+\delta_{(m_1-k_1,\ldots,m_s-k_s)}(\prod\limits_{i=1}^s({\bf D}_{-m_i}^{(i)})^*-\prod\limits_{i=1}^s{\bf D}_{0}^{(i)})}$$
\vskip 2 pt
\noindent
- If there is $i$ such that $m_i<k_i$, then
$$P_{(m_1,\ldots,m_s)}=\delta_{(m_1,\ldots,m_s)}((\prod\limits_{i=1}^s{\bf D}_{-1}^{(i)}-\prod\limits_{i=1}^s({\bf D}_{-m_i-k_i}^{(i)})^*).$$
\finth
\hfill
\eject
\vskip 1 pt
\proof
\vskip 1 pt
\noindent
-  If every $m_i$ is $\geq k_i$, then the restrictions to 
 ${{\cal O}}_{(m_1,\ldots,m_s)}(\widetilde\Xi)$ are given by 
$${\cal M}^{\sigma}{\cal D}=\prod\limits_{i=1}^s{\bf D}_0^{(i)},
{\cal D}{\cal M}^{\sigma}=\prod\limits_{i=1}^s{\bf  D}_{-1}^{(i)}, $$
$${\cal M}{\cal D}^{\sigma}=\prod\limits_{i=1}^s({{\bf D}_{-m_i}^{(i)}})^*, 
{\cal D}^{\sigma}{\cal M}=\prod\limits_{i=1}^s({{\bf D}_{-m_i-k_i}^{(i)}})^*.$$
\noindent
Then the   restriction of  $[\rho(E),\rho(F)]$ 
to ${{\cal O}}_{(m_1,\ldots,m_s)}(\widetilde\Xi)$ is  given by 
$$\eqalign{
&[\rho(E),\rho(F)]=\cr
&[{\cal M}^{\sigma }
-\delta\circ{\cal D}^{\sigma},{\cal M}-\delta\circ {\cal D} ]\cr
&=[{\cal M},\delta \circ {\cal D}^{\sigma}]+[\delta \circ {\cal D},{\cal M}^{\sigma}] \cr
&={\cal M}\delta {\cal D}^{\sigma}-\delta {\cal D}^{\sigma}{\cal M}
+\delta {\cal D} {\cal M}^{\sigma}-{\cal M}^{\sigma} \delta \circ {\cal D} \cr
&=\delta_{(m_1,\ldots,m_s)}({\cal D}{\cal M}^{\sigma}-{\cal D}^{\sigma}{\cal M})
+\delta_{(m_1-k_1,\ldots,m_s-k_s)}({\cal M} {\cal D}^{\sigma}-{\cal M}^{\sigma}{\cal D}) \cr
&=\delta_{(m_1,\ldots,m_s)}(\prod\limits_{i=1}^s{\bf D}_{-1}^{(i)}-\prod\limits_{i=1}^s({{\bf D}_{-m_i-k_i}^{(i)}})^*)\cr
&+\delta_{(m_1-k_1,\ldots,m_s-k_s)}(\prod\limits_{i=1}^s({{\bf D}_{-m_i}^{(i)}})^*-\prod\limits_{i=1}^s{\bf D}_0^{(i)}).}$$
\noindent
Then, by the Harish-Chandra isomorphism the operator $P_{(m_1,\ldots,m_s)}$ corresponds to the polynomial
$p_{(m_1,\ldots,m_s)}=\tilde\gamma (P_{(m_1,\ldots,m_s)})$,

$$p_{(m_1,\ldots,m_s)}(\lambda )=\delta_{(m_1,\ldots,m_s)}\bigl(\prod\limits_{i=1}^s\gamma_{-1}^{(i)}(\lambda^{(i)})-\prod\limits_{i=1}^s\gamma_{-m_i-k_i}^{(i)}(-\lambda^{(i)})\bigr)$$
$$
+\delta_{(m_1-k_1,\ldots,m_s-k_s)}\bigl(\prod\limits_{i=1}^s\gamma_{-m_i}^{(i)}(-\lambda^{(i)})-\prod\limits_{i=1}^s\gamma_0^{(i)}(\lambda^{(i)})\bigr)$$

\vskip 2 pt
\noindent
- If  there is $i$ such that $m_i < k_i$, then the restrictions to 
 ${{\cal O}}_{(m_1,\ldots,m_s)}(\widetilde\Xi)$ are : 
$$\eqalignno{{\cal M}^{\sigma}{\cal D}&=0,
{\cal D}{\cal M}^{\sigma}=\prod\limits_{i=1}^s{\bf  D}_{-1}^{(i)}, 
{\cal M}{\cal D}^{\sigma}=0, 
{\cal D}^{\sigma}{\cal M}=\prod\limits_{i=1}^s({{\bf D}_{-m_i-k_i}^{(i)}})^*.\cr}$$
\vskip 2 pt
\noindent
It follows that the restriction of $[\rho(E),\rho(F)]$ to ${{\cal O}}_{(m_1,\ldots,m_s)}(\widetilde\Xi)$ is  given by
$$\eqalignno{&[\rho(E),\rho(F)]=\delta\circ{\cal D}{\cal M}^{\sigma}-\delta\circ{\cal D}^{\sigma}{\cal M}\cr
&=\delta_{(m_1,\ldots,m_s)}[\prod\limits_{i=1}^s{\bf  D}_{-1}^{(i)}-\prod\limits_{i=1}^s({{\bf D}_{-m_i-k_i}^{(i)}})^*].\cr}$$
Then, by the Harish-Chandra isomorphism the operator $P_{(m_1,\ldots,m_s)}$ corresponds to the polynomial
$p_{(m_1,\ldots,m_s)}=\tilde\gamma (P_{(m_1,\ldots,m_s)})$,
$$p_{(m_1,\ldots,m_s)}(\lambda )=\delta_{(m_1,\ldots,m_s)}\bigl(\prod\limits_{i=1}^s\gamma_{-1}^{(i)}(\lambda^{(i)})-\prod\limits_{i=1}^s\gamma_{-m_i-k_i}^{(i)}(-\lambda^{(i)})\bigr)$$

\vskip 2 pt
\noindent
where 
\vskip 2 pt

\centerline{($\lambda=(\lambda^{(1)},\ldots,\lambda^{(s)})$ with
$\lambda^{(i)}=(\lambda_1^{(i)},\ldots,\lambda_{r_i}^{(i)})\in {\bboard C}^{r_i}$).}

\bigskip
The question is now whether it is possible to choose the sequence $(\delta_{(p_1,\ldots,p_s)})$
in such a way that $[\rho(E),\rho(F)]=\rho(H)$. 
\vskip 4 pt
\noindent
 Recall that,  restricted
to ${{\cal O}}_{(m_1,\ldots,m_s)}(\widetilde\Xi)$, 
$$\rho(H)={\cal E}-{(m_1r_1+\ldots+m_sr_s)\over 2},$$
where $\cal E$ is the Euler operator
$$({\cal E}\phi_{(m_1,\ldots,m_s)})((w_1,z_1),\ldots,(w_s,z_s))=(\prod\limits_{i=1}^sw_i^{m_i}){{\rm d}\over {\rm d}t} \big|_{t=0}(\prod\limits_{i=1}^s\psi_i(e^tz_i)).$$
Then  it amounts to checking that,
$$p_{(m_1,\ldots,m_s)}(\lambda)=\tilde\gamma ({\cal E})(\lambda )-{(m_1r_1+\ldots+m_sr_s)\over 2}\qquad (S)$$
\vskip 10pt
From now on, we introduce the subspaces   ${{\cal O}}_{q,\rm fin}(\widetilde\Xi)$ of ${{\cal O}}_{\rm fin}(\widetilde\Xi)$, for  $q=(q_1,\ldots,q_s)\in {\bboard N}^s$  defined by 
$${{\cal O}}_{q,\rm fin}(\widetilde\Xi)=\sum\limits_{m\in {\bboard N}}{{\cal O}}_{(k_1m+q_1,\ldots,k_sm+q_s)}(\widetilde\Xi)\qquad  (m_i=k_im+q_i)$$
\vskip 2 pt
\noindent
The  subspaces  ${{\cal O}}_{q,\rm fin}(\widetilde\Xi)$  which are obviously stable under $\rho(H)$,   are  stable  under  the operators $\rho(E)$ and $\rho(F)$ iff  $\exists i,  q_i<k_i$.
 \vskip 2 pt
\noindent
 In fact,   for  every $q$, the space ${{\cal O}}_{q,\rm fin}(\widetilde\Xi)$ is stable under the operators ${\cal M}$ and $\pi(\sigma)$.  But,  stability under the operator ${\cal D}$ (and ${\cal D}^{\sigma}$) requires the condition  ${\cal D}({{\cal O}}_{(k_1m+q_1,\ldots,k_sm+q_s)})=\{0\}$ for $m=0$,  which requires the existence of $i$ such that $\Delta^{k_i}({\partial\over\partial z_i})(\tilde{{\cal O}}_{q_i}(V_i))=0$. Using  (iii) of Theorem 2.1, one gets $q_i<k_i$ for some $i$. 

\vskip 4 pt
\noindent
From now on, assume that there is  $q_i$ such that $q_i < k_i$ and denote by 
 $\rho_q(H)$, $\rho_q(E)$ and $\rho_q(F)$  the   restrictions of $\rho(H)$, $\rho(E)$ and $\rho(F)$  to ${{\cal O}}_{q,\rm fin}(\widetilde\Xi)$.
 \vskip 1 pt
 \noindent
 \th Theorem 3.4| 
\vskip 1 pt
\noindent
(i)
 For every  $m\in {\bboard N}^*$, one can  choose  
$(\delta_{(k_1m+q_1,\ldots,k_sm+q_s)})$ such that 
$$\eqalignno{[\rho_q(H),\rho_q(E)]&=2\rho_q(E),  [\rho_q(H),\rho_q(F)]=-2\rho_q(F),
 [\rho_q(E),\rho_q(F)]=\rho_q(H)\cr}$$
if and only if  there is  a constant $\eta_q$  such that  for all $i$, ${q_i\over k_i}+{n_i\over k_ir_i}=\eta_q$, 
 and then,  for  $A=\prod _{i=1}^sk_i^{k_ir_i}$,
$$\delta_{(k_1m+q_1,\ldots,k_sm+q_s)}={1\over A(m+\eta_q )(m+\eta_q +1)}.$$ 
\vskip 2 pt
\noindent
(ii)  For $m=0$, one can  choose  
$(\delta_{(q_1,\ldots,q_s)})$ such that 
$$\eqalignno{[\rho_q(H),\rho_q(E)]&=2\rho_q(E),  [\rho_q(H),\rho_q(F)]=-2\rho_q(F),
 [\rho_q(E),\rho_q(F)]=\rho_q(H)\cr}$$
if and only if $q$   is   as in table 2 and $\delta_{q}$ as obtained below.
\finth

\vskip 2 pt
\proof    
(i) Assume $m\ne 0$. First,  recall (from Theorem 2.1 ) that  for  $m_i=k_im+q_i$,   the space ${{\cal O}}_{(m_1,\ldots,m_s)}(\widetilde\Xi)={{\cal O}}_{(k_1m+q_1,\ldots,k_sm+q_s)}(\widetilde\Xi)$,  consists in polynomial functions  $\prod\limits_{i=1}^s\psi_{m,i}(z_i)w_i^{k_im+q_i}$, where the polynomials  $\psi_{m,i}$ have degree $\leq k_ir_im+q_ir_i$.

$$\eqalign{
&\prod\limits_{i=1}^s\gamma _{\alpha }^{(i)} (\lambda^{(i)} )
=\prod _{i=1}^s \prod _{j=1}^{r_i} \big[\lambda _j^{(i)}-\alpha +\demi ({n_i\over r_i}-1)\big]_{k_i}\cr
&=\prod _{i=1}^s\prod _{j=1}^{r_i}\prod _{k=1}^{k_i}
\Bigl(\lambda _j^{(i)}-\alpha +\demi \bigl({n_i\over r_i}-1\Bigr)-(k-1)\Bigr)\cr
&=A\prod _{i=1}^s\prod _{j=1}^{r_i}\prod _{k=1}^{k_i}
\Bigl({\lambda _j^{(i)}\over k_i}-{\alpha\over k_i} +{1\over 2k_i}\bigl({n_i\over r_i}-1\bigr)-{k-1\over k_i}\Bigr).\cr}$$
\noindent
Denote by 
$$
X_{jk}^{(i)}={\lambda _j^{(i)}\over k_i}+{1\over 2k_i}\bigl({n_i\over r_i}-1\bigr)-{k-1\over k_i},$$
\noindent
and
$$b_{k_im+q_i}^{(i)}=m+{q_i\over k_i}+{n_i\over k_ir_i}-1.$$
\noindent
 Then 
$$\eqalignno{
&p_{(k_1m+q_1,\ldots,k_sm+q_s)}(\lambda )=A \delta _{(k_1m+q_1,\ldots,k_sm+q_s)}\times\cr
&\Bigl(\prod _{i=1}^s\prod _{j=1}^{r_i}\prod _{k=1}^{k_i} (X_{jk}^{(i)}+1)
-\prod _{i=1}^s\prod _{j=1}^{r_i}\prod _{k=1}^{k_i} (X_{jk}^{(i)}-b_{k_im+q_i}^{(i)}-1) \Bigr)\cr
&+A\delta _{(k_1(m-1)+q_1,\ldots,k_s(m-1)+q_s)}\times\cr
&\Bigl(\prod _{i=1}^s\prod _{j=1}^{r_i}\prod _{k=1}^{k_i}(X_{jk}^{(i)}-b_{k_im+q_i}^{(i)})
-\prod _{i=1}^s\prod _{j=1}^{r_i}\prod _{k=1}^{k_i}(X_{jk}^{(i)})\Bigr), \cr}$$
$$\tilde\gamma ({\cal E})(\lambda )=\sum _{i=1}^s\sum _{j=1}^{r_i}\sum _{k=1}^{k_i}X_{jk}^{(i)}
-\demi \sum _{i=1}^s\sum _{j=1}^{r_i}\sum _{k=1}^{k_i}b_{k_im+q_i}^{(i)}.$$
\vskip 1 pt
\noindent
The restriction of $\rho(H)$ to ${{\cal O}}_{(k_1m+q_1,\ldots,k_sm+q_s)}(\widetilde\Xi)$, is  given by   
$$\rho_q(H)={\cal E}-\sum\limits_{i=1}^s{k_ir_im+q_ir_i\over 2}={\cal E}-2m-\sum\limits_{i=1}^s{q_ir_i\over 2}.$$

\noindent
The identity (S) is then proved  by using the following  lemma:

\vskip 1 pt
\th Lemma 3.5|
To a partition $p=(p_1,\ldots ,p_{\ell })$ of 4 and length $\ell $:
$p_1+\dots +p_{\ell }=4,$
and  numbers $\gamma _{ij}$ ($1\leq i\leq \ell $, $1\leq j\leq p_i-1$), 
one associates
the polynomial $F$ in the $\ell $ variables $T_1,\ldots ,T_{\ell }$:
$F(T_1,\ldots ,T_{\ell })=\prod _{i=1}^{\ell } T_i\prod _{j=1}^{p_i-1}(T_i+\gamma _{ij}).$
Given $\alpha ,\beta,c\in {\bboard R}$, and $b_1,\ldots b_{\ell }\in {\bboard R}$, then

\centerline{$
\alpha \bigl(F(T_1+1,\ldots ,T_{\ell }+1) 
-F(T_1-b_1-1,\ldots ,T_{\ell }-b_{\ell }-1)\bigr)$}

\centerline{$+\beta \bigl(F(T_1-b_1,\ldots ,T_{\ell }-b_{\ell })-F(T_1,\ldots ,T_{\ell }\bigr) 
=\sum _{i=1}^{\ell }T_i+c $}

\noindent
is an identity in the variables $T_1,\ldots ,T_{\ell }$ if and only if
there exists $b\ne 0$ such that  $b_1=\ldots =b_{\ell }=b,\ \alpha ={1\over (b+1)(b+2)},\ \beta ={1\over b(b+1)}$ and $$ c=\sum _{i=1}^{\ell }\sum _{j=1}^{k_i-1}\gamma _{ij}-2b.$$
\finth

\vskip 1 pt
\noindent
By this lemma, we see in particular  that for every $m$,  the  $b_{k_im+q_i}^{(i)}$ must be equal  for every $i$,  which means that there is $\eta_q={q_i\over k_i}+{n_i\over k_ir_i} \quad  ( \forall i)$
\vskip 4pt
\noindent
(ii) Now,  consider the situation of $m=0$ and determine the cases where the constant $\eta_0$ and the  constants $\delta_{(q_1,\ldots,q_s)}$ exist.
\vskip 4 pt
\noindent
Assume $m=0$. The $\psi_{0,i}$ have degree $\leq q_ir_i$.  
 Let's consider the restriction of the operators ${\cal D}{\cal M}^{\sigma}, {\cal D}^{\sigma}{\cal M}, {\cal M}{\cal D}^{\sigma}, {\cal M}^{\sigma}{\cal D}$ and $\rho(H)$ to the space ${{\cal O}}_{(q_1,\ldots,q_s)}(\widetilde\Xi)$. 
The restriction of $\rho(H)$ to ${{\cal O}}_{(q_1,\ldots,q_s)}(\widetilde\Xi)$, is  given by   
$$\rho_q(H)={\cal E}-\sum\limits_{i=1}^s{q_ir_i\over 2}.$$
\vskip 1 pt
\noindent
 By Proposition 3.3, the restriction  of $[\rho(E),\rho(F)]$ 
to ${{\cal O}}_{(q_1,\ldots,q_s)}(\widetilde\Xi)$,   is given by

$$\delta_{(q_1,\ldots,q_s)}[(\prod\limits_{i=1}^s{\bf D}_{-1}^{(i)}-\prod\limits_{i=1}^s({{\bf D}_{-p_i-k_i}^{(i)}})^*)$$
\noindent
then
$$p_{(q_1,\ldots,q_s)}(\lambda )$$
$$=
A \delta _{(q_1,\ldots,q_s)}\Bigl(\prod _{i=1}^s\prod _{j=1}^{r_i}\prod _{k=1}^{k_i} (X_{jk}^{(i)}+1)
-\prod _{i=1}^s\prod _{j=1}^{r_i}\prod _{k=1}^{k_i} (X_{jk}^{(i)}-b_{q_i}^{(i)}-1) \Bigr).$$
\noindent
Furthermore, 
$$\tilde\gamma ({\cal E})(\lambda )=\sum _{i=1}^s\sum _{j=1}^{r_i}\sum _{k=1}^{k_i}X_{jk}^{(i)}
-\demi \sum _{i=1}^s\sum _{j=1}^{r_i}\sum _{k=1}^{k_i}b_{q_i}^{(i)}.$$

\noindent
Let's  check the identity (S) case by case.

\vskip 2pt
\noindent
(1)  $\widetilde{\rm SL}(3,{\bboard R})$ :  $V={\bboard C},  Q(z)=z^4$,   $r=1, k=4$, $q < 4$, $\eta_q={q\over 4}+{1\over 4}$.
The functions $\phi_{q}$  are  given by $\phi_{q}(w,z)=z^{\alpha}w^q$ with $\alpha \leq q$.  The identity   (S)  becomes
\vskip 2 pt
\noindent

\centerline{$\delta_{q}\Bigl(\prod\limits_{j=1}^4(\alpha+j)-\prod\limits_{j=1}^4(q-\alpha+j)\Bigr)=\alpha-{q\over 2}.$}

\vskip 6 pt

\noindent
One has to determine $\delta_{q}$ such that this identity is available for every $\alpha \leq  q$ such that $\alpha \geq q-\alpha$, i.e. for ${q\over 2} \leq \alpha \leq q$. It follows that 
for $q=0$, the identity is trivial and $\delta_{0}$ is arbitrarily,  for  $q=1$,   $\alpha=1$ and $\delta_{1}={1\over 2(120-24)}$,  for $q=2$, $\alpha=2$ and $\delta_{2}$ is arbitrarily and for $q=3$, $\alpha=2$ or $\alpha=3$ then (S) is impossible.

\vskip 2 pt
\noindent
(2) $\widetilde{\rm SL}(p+2,{\bboard R})$ : $V={\bboard C}^p,  Q(z)=(z_1^2+\ldots+z_p^2)^2$, $r=2, k=2$,  $q\in \{0,1\}$, $\eta_q={q\over 2}+{p\over 4}$.
The functions $\phi_{q}$  are  given by $\phi_{q}(w,z)=\Delta(z)^{\alpha}w^q$ with $\alpha  \leq q $.  The identity (S)  becomes 
\vskip 2 pt

\centerline{$\delta_{q}\Bigr((\alpha+2)(\alpha +1)-(q-\alpha+2)(q-\alpha+1)\Bigl)=2\alpha-q.$}

\vskip 6 pt

\noindent
One has to determine $\delta_{q}$ such  that  (S) occurs  for every $\alpha$ such that  ${q\over 2} \leq  \alpha \leq q$.  If $q=0$, the identity is trivial and $\delta_0$ is arbitrarily, and  if $q=1$, it  becomes $\delta_{1}(8\alpha-4)=2\alpha-1$, then $\delta_1={1\over 4}$.
\hfill
\eject
\vskip 2 pt
\noindent
(3) $\widetilde{\rm SO}(3,3)$ : $V={\bboard C}\oplus{\bboard C}, Q(z)=z_1^2z_2^2$,  $r_1=r_2=1, k_1=k_2=2$, $q_1=q_2=:\tilde q\in\{0,1], \eta_q={\tilde q\over 2}+{1\over 2}$.
The functions $\phi_{q}$  are  given by $\phi_{(q_1,q_2)}((w_1,z_1),(w_2,z_2))=z_1^{\alpha_1}z_2^{\alpha_2}w_1^{q_1}w_2^{q_2}$ with  $\alpha_1\leq 1$ and  $\alpha_2\leq 1$.
The identity (S)  becomes 
\vskip 2 pt
\noindent

\centerline{$\delta_{(q_1,q_2)}(\alpha_1+2)(\alpha_1 +1)(\alpha_2+2)(\alpha_2+1)$}

\vskip 2 pt

\centerline{$-\delta_{(q_1,q_2)}(q_1-\alpha_1+2)(q_1-\alpha_1+1)(q_2-\alpha_2+2)(q_2-\alpha_2+1)=\alpha_1+\alpha_2-\tilde q.$}

\vskip 6 pt
\noindent
One has to determine $\delta_q=\delta_{(q_1,q_2)}$ such  that this identity is available for every $\alpha_i$ such that  ${q_i\over 2} \leq  \alpha_i \leq q_i$. It follows that, if $\tilde q=0$, then  $\alpha_i=0$ and $\delta_{(0,0)}$ is arbitrarily, and if $\tilde q=1$, then $\alpha_1=\alpha_2=1$ and $\delta_{(1,1)}={1\over 36-4}$.
\vskip 4pt
\noindent
(4) $\widetilde{\rm SO}(4,4)$ : $V=\oplus^4{\bboard C}$, $Q(z)=z_1z_2z_3z_4$, $r_i=1=k_i$,  $q_i=\tilde q=0$,  $\eta_q=\tilde q+1=1$. 
Then  $\phi_{q}=\phi_0$  are   constants,   (S) is trivial.

\vskip 2 pt
\noindent
(5) $\widetilde{\rm SO}(p+2,3)$ :  $V={\bboard C}^p\oplus{\bboard C}, Q(z)=(z_1^2+\ldots+z_p^2){z'}^2$, $p$ odd, $r_1=2, r_2=1, k_1=1, k_2=2$, $q_1=0$,  $q_2={p-1\over 2}$, $\eta_q={p\over 2}$.
The functions $\phi_{q}$  are  given by $\phi_{(q_1,q_2)}((w,z),(w',z'))=(z')^{\alpha_2}(w')^{q_2}$ with $\alpha_2\leq q_2 $.
The identity (S)  becomes 
\vskip 2 pt

\centerline{$\delta_{(0,q_2)}[(\alpha_2+2)(\alpha_2+1)-(q_2-\alpha_2+2)(q_2-\alpha_2+1)]=\alpha_2-{p-1\over 4}$}

\vskip 2 pt
\noindent
i.e.,
\vskip 2 pt

\centerline{$\delta_{(0,q_2)}(p+5)(\alpha_2 -{p-1\over 4})=\alpha_2-{p-1\over 4}.$}

\vskip 6 pt
\noindent
One has to determine $\delta_{(0,q_2)}$ such  that this identity is available for every  $\alpha_2$ such that   ${q_2\over 2} \leq \alpha_2 \leq q_2$. It follows that   $\delta_{(0,q_2)}={1\over p+5}$. 

\vskip 2 pt
\noindent
(6) $\widetilde{\rm SO}(p+2,4)$ :  $V={\bboard C}^p\oplus {\bboard C}\oplus{\bboard C}$,  $Q(z)=(z_1^2+\ldots+z_p^2)z'z''$, $p$ even,  $r_1=2, r_2=r_3=1, k_i=1$, $(q_1,q_2,q_3)=(0,{p\over 2}-1,{p\over 2}-1)$,  $\eta_q={p\over 2}$.
The functions $\phi_{q}$  are  given by 
\vskip 2 pt

\centerline{$\phi_{(q_1,q_2,q_3)}((w,z),(w',z'),(w'',z''))=(z')^{\alpha_2}(z'')^{\alpha_3}(w')^{q_2}(w'')^{q_3}.$}

\vskip 2 pt
\noindent
 with $\alpha_2\leq q_2,  \alpha_3 \leq q_3 $.
The identity (S)  becomes 
\vskip 2 pt

\centerline{$\delta_{(0,q_2,q_3)}[(\alpha_2+1)(\alpha_3+1)-(q_2-\alpha_2+1)(q_3-\alpha_3+1)]=\alpha_2+\alpha_3-({p\over 2}-1), $.}

\vskip 2pt
\noindent
i.e.

\vskip 2 pt

\centerline{$\delta_{(0,q_2,q_3)}({p\over 2}+1)[\alpha_2+\alpha_3-({p\over 2}-1)]=\alpha_2+\alpha_3-({p\over 2}-1).$}

\vskip 6 pt
\noindent
One has to determine $\delta_{(0,q_2,q_3)}$ such  that this identity is available for every  $\alpha_2$ and $\alpha_3$ such that  ${q_i\over 2} \leq  \alpha_i \leq q_i$. It follows that     $\delta_{(0,q_2,q_3)}={1\over{p\over 2}+1}$.  
\hfill
\eject
\vskip 2 pt
\noindent
(7) $\widetilde{\rm SO}(p_1+2,p_2+2)$ : 
$V={\bboard C}^{p_1}\oplus{\bboard C}^{p_2}, Q=(z_1^2+\ldots+z_{p_1}^2)({z_1'}^2+\ldots+{z_{p_2}'}^2)$,
\vskip 1 pt
\noindent
$r_1=r_2=2, k_1=k_2=1$,  $p_1-p_2\geq 0$ even,    $q_1=0$,  $q_2={p_1-p_2\over 2}$,   $\eta_q={p_1\over 2}={p_1\over 2}$.
The functions $\phi_{q}$  are  given by 
\vskip 2 pt

\centerline{$\phi_{(q_1,q_2)}((w_1,z),(w_2,z'))=\Delta(z')^{\alpha_2}w_1^{q_1}w_2^{q_2}$}

\vskip 2 pt
\noindent
 with $\alpha_2\leq q_2$.
The identity (S)  becomes 

\vskip 2 pt

\centerline{$\delta_{(0,q_2)}[(\alpha_2+1)-(q_2-\alpha_2+1)]=2\alpha_2-{p_1-p_2\over 2}$}

\vskip 2 pt
\noindent
i.e.,
\vskip 2 pt

\centerline{$\delta_{(0,q_2)}(2\alpha_2-{p_1-p_2\over 2})=2\alpha_2-{p_1-p_2\over 2}.$}

\vskip 6 pt
\noindent
One has to determine $\delta_{(0,q_2)}$ such  that this identity is available for every  $\alpha_2$  such that  ${q_2\over 2} \leq  \alpha_2 \leq q_2$.  It follows that  if  $p_1=p_2$, then $q=0$ and $\delta_{0}$ is arbitrarily, and if  $p_1\ne p_2$, then    $\delta_{(0,{p_1-p_2\over 2})}=1$.  
\vskip 2 pt
\noindent
(8) $E_{6(6)}, E_{7(7)},E_{8(8)}$ : $(V,Q(z))$ is respectively the pair

\vskip 1 pt

\centerline{$({\rm Sym}(4,{\bboard C}), \det z)$ ;  $(M(4,{\bboard C}), \det z)$ ;  $({\rm Skew}(8,{\bboard C}),{\rm Pfaff}(z))$,}
\noindent
 $r=4, k=1$, $q=0$,  $\eta_q=q+{n\over 4}=1+{3d\over 2}$, with $d=1, 2$ or $4$ respectively.
 \vskip 1pt 
\noindent
The functions $\phi_{0}$  are  given by $\phi_{q}(w,z)=c$, 
the identity (S) is  trivial.

\vskip 2 pt
\noindent
(9)   $F_{4(4)}, E_{6(2)},E_{7(-5)}, E_{8(-24)}$ :  $(V,Q(z))$ is respectively the pair

\vskip 1 pt

\centerline{$({\rm Sym}(3,{\bboard C})\oplus {\bboard C},\det(z_1)\cdot z_2), \quad (M(3,{\bboard C})\oplus {\bboard C}, \det(z_1)\cdot z_2)$}

\vskip 1 pt

\centerline{(${\rm Skew}(6,{\bboard C})\oplus {\bboard C},{\rm Pfaff}(z_1)\cdot z_2, \quad ({\rm Herm}(3,{\bboard C})\oplus {\bboard C},\det(z_1)\cdot z_2)$,}
\noindent
$r_1=3, r_2=1, k_1=k_2=1$,  $q_1=0$,    $q_2={n_1\over 3}-1$,  $\eta_q=q_1+{n_1\over 3}=1+d_1$.
\vskip 1pt 
\noindent
The functions $\phi_{q}$  are  given by $\phi_{(q_1,q_2)}((w_1,z_1),(w_2,z_2))=z_2^{\alpha_2}w_2^{q_2}$ with $\alpha_2\leq q_2$.
The identity (S)  becomes 

\vskip 2 pt

\centerline{$\delta_{(0,q_2)}[(\alpha_2+1)-(q_2-\alpha_2+1)]=\alpha_2-{q_2\over 2}$}
\vskip 2 pt
\noindent
i.e.
\vskip 2 pt

\centerline{
$\delta_{(0,q_2)}(2\alpha_2-q_2)=\alpha_2-{q_2\over 2}.$}

\vskip 6 pt
\noindent
One has to determine $\delta_{(0,q_2)}$ such  that this identity is available for every  $\alpha_2$  such that  ${q_2\over 2} \leq  \alpha_2 \leq q_2$. Then  $\delta_{(0,q_2)}={1\over 2}$.

\vskip 2 pt
\noindent
(10) $G_{2(2)} : V={\bboard C}\oplus{\bboard C}, Q=z_1^3z_2,  r_1=r_2=1, k_1=3, k_2=1, q_1=2, q_2=0, \eta_q={q_1\over 3}+{1\over 3}=1$.
The functions $\phi_{q}$  are  given by $\phi_{(2,0)}((w_1,z_1),(w_2,z_2))=z_1^{\alpha_1}w_1^{q_1}$ with $\alpha_1\leq 2$.
The identity (S)  becomes 
\vskip 2 pt

\centerline{$\delta_{(2,0)}[(\alpha_1+1)-(q_1-\alpha_1+1)]=\alpha_1-{q_1\over 2}$}

\vskip 2 pt
\noindent
i.e.

\vskip 2 pt

\centerline{$\delta_{(2,0)}(2\alpha_1-2)=\alpha_1-1.$}

\vskip 6 pt
\noindent
One has to determine $\delta_{(2,0)}$ such  that this identity is available for every  $\alpha_1$  such that  ${q_1\over 2}=1 \leq  \alpha_1 \leq q_1=2$. Then  $\delta_{(2,0)}={1\over 2}$.

\bigskip
\noindent
\hfill
\eject
Now, we extend $\rho$ to a linear map from ${\goth p}$ to   the space End(${\cal O}_{\rm fin})$.
\vskip 4 pt
\th Proposition 3.6|
\vskip 1 pt
\noindent
{\rm (i)} There is a unique linear  map 
${\goth p}\to {\rm End}\bigl({{\cal O}}_{\rm fin}(\widetilde\Xi)\bigr),\quad p\mapsto {\cal M}(p),$
such that ${\cal M}(1)={\cal M}$,  ${\cal M}(Q)={\cal M}^{\sigma}$ and, for $X\in{\goth k}$,
${\cal M}\bigl([X,p]\bigr)=[{\rm d}\pi(X),{\cal M}(p)].$
\vskip 4 pt

\noindent
{\rm (ii)} 
There is a unique linear map 
${\goth p}\to {\rm End}\bigl({{\cal O}}_{\rm fin }(\widetilde\Xi)\bigr),\quad p\mapsto {\cal D}(p),$
such that ${\cal D}(1)={\cal D}$,   ${\cal D}(Q)={\cal D}^{\sigma}$ and,   for $X\in{\goth k}$,
${\cal D}\bigl([X,p]\bigr)=[{\rm d}\pi(X),{\cal D}(p)].$
\vskip 4 pt

\noindent
{\rm (iii)}  ${\cal M}\bigl(\tilde\kappa (g)p\bigr)=\pi (g){\cal M}(p)\pi (g^{-1})$ and  ${\cal D}\bigl(\tilde\kappa (g)p\bigr)=\pi (g){\cal D}(p)\pi (g^{-1}).$

\finth
\vskip 1 pt
\proof Recall  that  ${\goth p}={\goth p}_{-2}+{\goth p}_{-1}+{\goth p}_0+{\goth p}_1+{\goth p}_2, {\goth k}={\goth k}_{-1}+{\goth k}_0+{\goth k}_1$  with
${\goth k}_{-1} \simeq V\simeq {\goth k}_1=\kappa(\sigma)({\goth k}_{-1})$,  ${\goth p}_{-2}={\bboard C}\cdot 1, {\goth p}_2={\bboard C}\cdot Q$  and that  for  every $ p_1\in {\goth p}_1, p_{-1}\in {\goth p}_{-1}, $ there is  a unique $X_{-1} \in {\goth k}_{-1}$ and a unique   $X_1\in 
{\goth k}_1$ such that  $p_1={\rm d}\kappa(X_{-1})Q=[X_{-1},Q]$ and  $p_{-1}={\rm d}\kappa(X_1)1=[X_1,1]$. 
\vskip 2 pt
\noindent
 Furthermore,  since ${\goth p}={\cal U}({\goth k})Q={\cal U}({\goth k}_{-1}+{\goth k}_0)Q$,  it follows that the space ${\goth p}_0$ is given by ${\goth p}_0=[{\goth k}_{-1},{\goth p}_1]$.

 \vskip 4 pt
 \noindent
{\rm (i)}  Let ${\cal M} : {\goth p} \rightarrow $End$({\cal O}_{\rm fin})$ be the linear map defined as follows:
\vskip 4 pt
\noindent
- On ${\goth p}_{-2}$ and ${\goth p}_2$: ${\cal M}(1)={\cal  M} , {\cal M}(Q)={\cal M}^{\sigma}$.

\vskip 4 pt

\noindent
- On ${\goth p}_1$:  for  $p_1\in {\goth p}_1$, let $X_{-1}\in {\goth k}_{-1}$ such that $p_1=[X_{-1},Q]$, then put 
\vskip 4 pt

\centerline{${\cal M}(p_1)={\cal M}({\rm d}\kappa(X_{-1})Q)=[{\rm d}\pi(X_{-1}),{\cal M}(Q)]$.}

\vskip 4 pt
\noindent
- On ${\goth p}_{-1}$:  for  $p\in {\goth p}_{-1}$, let $X_1\in {\goth k}_{1}$ such that $p_{-1}=[X_1,1]$, then put 
\vskip 4 pt

\centerline{${\cal M}(p_{-1})={\cal M}({\rm d}\kappa(X_1)1)=[{\rm d}\pi(X_1),{\cal M}(1)]$.}

\vskip 4 pt
\noindent
- On ${\goth p}_0$ :   for  $p_0\in {\goth p}_0$ there is $(X_{-1},p_1)\in {\goth k}_{-1}\times{\goth p}_1$ such that $p_0=[X_{-1},p_1]$.  Observe that if there is another $(Y_{-1},q_1)\in {\goth k}_{-1}\times{\goth p}_1$ such that  $[X_{-1},p_1]=[Y_{-1},q_1]$  then  $[{\rm d}\pi(X_{-1}),{\cal M}(p_1)]=[{\rm d}\pi(Y_{-1}),{\cal M}(q_1)].$ 
\vskip 4pt
\noindent
 In fact, $p_1=[X_{-1}',Q]$ and $q_1=[Y_{-1}',Q]$, with $X_{-1}', Y_{-1}' \in {\goth k}_{-1}$. Let 
exp$(X_{-1})=\tau_a$,  exp$(X_{-1}')=\tau_{a'}$,   exp$(Y_{-1})=\tau_b$ and  exp$(Y_{-1}')=\tau_{b'}$, where 
$a=\sum\limits_{i=1}^sa_i,  \quad a'=\sum\limits_{i=1}^sa_i', \quad b=\sum\limits_{i=1}^sb_i, \quad b'=\sum\limits_{i=1}^sb_i'$.
Furthermore, for $\phi(w,z)=\prod\limits_{i=1}^sw_i^{m_i}\psi_i(z_i)$ and ${\omega}=\prod\limits_{i=1}^sw_i^{m_i+k_i}$,  
\vskip 2 pt

\centerline{$({\cal M}(Q)\phi)(w,z)=({\cal M}^{\sigma}\phi)(w,z)=\prod\limits_{i=1}^sw_i^{m_i+k_i}\Delta_i^{k_i}(z_i)\psi_i(z_i),$}

\vskip 2 pt
\noindent
 $$\eqalign{&([{\rm d}\pi(X_{-1}),{\cal M}(p_1)]\phi)(w,z)=([{\rm d}\pi(X_{-1}),[{\rm d}\pi(X_{-1}'),{\cal M}(Q)]]\phi)(w,z)\cr
&=({\rm d}\pi(X_{-1}){\rm d}\pi(X_{-1}'){\cal M}(Q)\phi)(w,z)-({\rm d}\pi(X_{-1}){\cal M}(Q){\rm d}\pi(X_{-1}')\phi)(w,z)\cr
 &-({\rm d}\pi(X_{-1}'){\cal M}(Q){\rm d}\pi(X_{-1})\phi)(w,z)+({\cal M}(Q){\rm d}\pi(X_{-1}'){\rm d}\pi(X_{-1})\phi)(w,z)\cr
&={\omega}{\rm d}\pi(X_{-1})({{\rm d}\over {\rm d}t}_{\mid_{t=0}}((\prod\limits_{i=1}^s\Delta_i^{k_i}(z_i+ta_i'))(\prod\limits_{i=1}^s\psi_i(z_i+ta_i'))) \cr
&-{\omega}{\rm d}\pi(X_{-1})((\prod\limits_{i=1}^s\Delta_i^{k_i}(z_i)){{\rm d}\over {\rm d}t}_{\mid_{t=0}}(\prod\limits_{i=1}^s\psi_i(z_i+ta_i'))\cr
&-{\omega}{\rm d}\pi(X_{-1}')( (\prod\limits_{i=1}^s\Delta_i^{k_i}(z_i)){{\rm d}\over {\rm d}t}_{\mid_{t=0}}(\prod\limits_{i=1}^s\psi_i(z_i+ta_i))    \cr
&+{\omega}(\prod\limits_{i=1}^s\Delta_i^{k_i}(z_i)){\rm d}\pi(X_{-1}')({{\rm d}\over {\rm d}t}_{\mid_{t=0}}(\prod\limits_{i=1}^s\psi_i(z_i+ta_i)))   \cr
&={\omega}{\rm d}\pi(X_{-1})(({{\rm d}\over {\rm d}t}_{\mid_{t=0}}(\prod\limits_{i=1}^s\Delta_i^{k_i}(z_i+ta_i')))(\prod\limits_{i=1}^s\psi_i(z_i))   \cr
&-{\omega}({{\rm d}\over {\rm d}r}_{\mid_{r=0}}(\prod\limits_{i=1}^s\Delta_i^{k_i}(z_i+ra_i')))({{\rm d}\over {\rm d}t}_{\mid_{t=0}}(\prod\limits_{i=1}^s\psi_i(z_i+ta_i))) \cr
&={\omega}[X_{-1},[X_{-1}',Q]](\prod\limits_{i=1}^s\psi_i(z_i))\cr
&+{\omega}({{\rm d}\over {\rm d}t}_{\mid_{t=0}}(\prod\limits_{i=1}^s\Delta_i^{k_i}(z_i+ta_i')))({{\rm d}\over {\rm d}r}_{\mid_{r=0}}(\prod\limits_{i=1}^s\psi_i(z_i+ra_i)))\cr
&-{\omega}({{\rm d}\over {\rm d}r}_{\mid_{r=0}}(\prod\limits_{i=1}^s\Delta_i^{k_i}(z_i+ra_i')))({{\rm d}\over {\rm d}t}_{\mid_{t=0}}(\prod\limits_{i=1}^s\psi_i(z_i+ta_i)))\cr
&={\omega}[X_{-1},[X_{-1}',Q]](\prod\limits_{i=1}^s\psi_i(z_i)).}$$
Since  $[X_{-1},[X_{-1}',Q]]=[Y_{-1},[Y_{-1}',Q]]$,  then  
$$[{\rm d}\pi(X_{-1}),[{\rm d}\pi(X_{-1}'),{\cal M}(Q)]]=[{\rm d}\pi(Y_{-1}),[{\rm d}\pi(Y_{-1}'),{\cal M}(Q)]],$$
i.e.
$$[{\rm d}\pi(X_{-1}),{\cal M}(p_1)]=[{\rm d}\pi(Y_{-1}),{\cal M}(q_1)].$$

  \vskip 4 pt
 \noindent
Then, we  define ${\cal M}(p_0)$ by 
${\cal M}(p_0)={\cal M}([X_{-1},p_1])=[{\rm d}\pi(X_{-1}),{\cal M}(p_1)].$
\hfill
\eject
\noindent
{\rm (ii)} Let ${\cal D} : {\goth p} \rightarrow $End$({\cal O}_{\rm fin})$ be the linear map  defined as follows:
\vskip 2 pt
\noindent
- On ${\goth p}_{-2}$ and ${\goth p}_2$:
\vskip 1 pt

\centerline{${\cal D}(1)={\cal  D} , {\cal D}(Q)={\cal D}^{\sigma}$,}

\vskip 2 pt
\noindent
- On ${\goth p}_1$:
\vskip 1 pt

\centerline{${\cal D}(p_1)={\cal D}({\rm d}\kappa(X_{-1})Q)=[{\rm d}\pi(X_{-1}),{\cal D}(Q)]$,}

\vskip 2 pt
\noindent
- On ${\goth p}_{-1}$:  

\centerline{${\cal D}(p_{-1})={\cal D}({\rm d}\kappa(X_1)1)=[{\rm d}\pi(X_1),{\cal D}(1)]$,}

\vskip 2 pt
\noindent
- On ${\goth p}_0$ :   for  $p_0\in {\goth p}_0$ there is $(X_{-1},X_{-1}')\in {\goth k}_{-1}\times{\goth k}_{-1}$ such that $p_0=[X_{-1},[X_{-1}',Q]]]$.  For $\phi(w,z)=\prod\limits_{i=1}^sw_i^{m_i}\psi_i(z_i)$ and  ${\omega}=\prod\limits_{i=1}^sw_i^{m_i-k_i}$,  
$$({\cal D}(Q)\phi)(w,z)=({\cal D}^{\sigma}\phi)(w,z)=\prod\limits_{i=1}^sw_i^{m_i-k_i}(({\bf D}_{-m_i}^{(i)})^*)\psi_i)(z_i),$$
$$\eqalign{&([{\rm d}\pi(X_{-1}),[{\rm d}\pi(X_{-1}'),{\cal D}(Q)]]\phi)(w,z)\cr
 &=({\rm d}\pi(X_{-1}){\rm d}\pi(X_{-1}'){\cal D}(Q)\phi)(w,z)-({\rm d}\pi(X_{-1}){\cal D}(Q){\rm d}\pi(X_{-1}')\phi)(w,z)\cr
 &-({\rm d}\pi(X_{-1}'){\cal D}(Q){\rm d}\pi(X_{-1})\phi)(w,z)+({\cal D}(Q){\rm d}\pi(X_{-1}'){\rm d}\pi(X_{-1})\phi)(w,z)\cr
&={\omega}{\rm d}\pi(X_{-1})({{\rm d}\over {\rm d}t}_{\mid_{t=0}}((\prod\limits_{i=1}^s(({\bf D}_{-m_i}^{(i)})^*)\psi_i)(z_i+ta_i')\cr
&-{\omega}{\rm d}\pi(X_{-1})((\prod\limits_{i=1}^s(({\bf D}_{-m_i}^{(i)})^*)({{\rm d}\over {\rm d}t}_{\mid_{t=0}}(\prod\limits_{i=1}^s\psi_i(z_i+ta_i')))\cr
&-{\omega}{\rm d}\pi(X_{-1}')((\prod\limits_{i=1}^s(({\bf D}_{-m_i}^{(i)})^*)({{\rm d}\over {\rm d}t}_{\mid_{t=0}}(\prod\limits_{i=1}^s\psi_i(z_i+ta_i)))\cr
&+{\omega}\prod\limits_{i=1}^s(({\bf D}_{-m_i}^{(i)})^*)({{\rm d}\over {\rm d}r}_{\mid_{r=0}}({{\rm d}\over {\rm d}t}_{\mid_{t=0}}(\prod\limits_{i=1}^s\psi_i(z_i+ra_i'+ta_i)))).}$$
As in the case of ${\cal M}$ (but with  more technical calculations), one can show that if 
 $[X_{-1},[X_{-1}',Q]]=[Y_{-1},[Y_{-1}',Q]]$  with $(Y_{-1},Y_{-1}')\in {\goth k}_{-1}\times {\goth k}_{-1}$, then   $[{\rm d}\pi(X_{-1}),[{\rm d}\pi(X_{-1}'),{\cal D}(Q)]]=[{\rm d}\pi(Y_{-1}),[{\rm d}\pi(Y_{-1}'),{\cal D}(Q)]].$
\vskip 2 pt
\noindent
Then, we  define ${\cal D}(p_0)$ by ${\cal D}(p_0)={\cal D}([X_{-1},p_1])=[{\rm d}\pi(X_{-1}),{\cal D}(p_1)].$
\vskip 4 pt

\vskip 1 pt
\noindent
For $X\in {\goth k}, p\in {\goth p}$, one can show

\vskip 6 pt

\centerline{${\cal M}([X,p])=[{\rm d}\pi(X),{\cal M}(p)]$ and ${\cal D}([X,p])=[{\rm d}\pi(X),{\cal D}(p)]. \qquad  (*)$}

\vskip 6 pt

\noindent
Using property $(*)$, one can  deduce (iii). Finally,  using the irreducibility  of the representation ${\rm d}\kappa$ of ${\goth k}$ on ${\goth p}$, we deduce the  uniqueness of ${\cal M}$ and ${\cal D}$.

\vskip 1 pt

The proof of  property $(*)$ for   ${\cal M}$ and  for ${\cal D}$ consists in many cases,  based on the grading  decompositions of  ${\goth k}$ and  ${\goth p}$.
\hfill
\eject
\noindent
1)   For $p\in {\goth p}_2$, i.e.  up to a scalar, $p=Q$, one has:
\vskip 2 pt
\noindent
(i) for $X\in {\goth k}_{-1}$, the formula (*) is  true by  definition.
\vskip 2 pt
\noindent
(ii) for $X\in {\goth k}_0$, there is $\alpha(X)\in {\bboard C}$ such that $[X,Q]=\alpha(X)Q$. Then, 
\vskip 2 pt

\centerline{$({\cal M}(Q)d\pi(X)\phi)(w,z)=\prod\limits_{i=1}^sw_i^{k_i}\Delta_i^{k_i}(z_i){d\over dt}_{\mid_{t=0}}\phi(e^{-t\alpha(X)}w,$exp$(tX)z)=$}

\vskip 2  pt

\centerline{$(\prod\limits_{i=1}^sw_i^{m_i+k_i}\Delta_i^{k_i}(z_i))[-\alpha(X)(\sum\limits_{i=1}^s m_i)(\prod\limits_{i=1}^s\psi_i(z_i)+{{\rm d}\over {\rm d}t}_{\mid_{t=0}}(\prod\limits_{i=1}^s\psi_i($exp$(tX_i)z_i)]$}

\vskip 2 pt

 then $[{\rm d}\pi(X),{\cal M}(Q)]\phi)(w,z)={{\rm d}\over {\rm d}t}_{\mid_{t=0}}{\cal M}(Q)\phi(e^{-t\alpha(X)}w,$exp$(tX)z)-$

 \vskip 2  pt
 
 \centerline{$(\prod\limits_{i=1}^sw_i^{m_i+k_i}\Delta_i^{k_i}(z_i))[-\alpha(X)(\sum\limits_{i=1}^sm_i)(\prod\limits_{i=1}^s\psi_i(z_i)+{d\over dt}_{\mid_{t=0}}(\prod\limits_{i=1}^s\psi_i($exp$(tX_i)z_i)]$}

\vskip 2 pt

\centerline{$={{\rm d}\over {\rm d}t}_{\mid_{t=0}}(\prod\limits_{i=1}^sw_i^{m_i+k_i})e^{-t\alpha(X)(\sum\limits_{i=1}^sm_i)}(\prod\limits_{i=1}^s\Delta_i^{k_i}($exp$(tX_i)z_i)\psi_i($exp$(tX_i)z_i))$}

\vskip 2 pt

\centerline{$-(\prod\limits_{i=1}^sw_i^{m_i+k_i}\Delta_i^{k_i}(z_i))[-\alpha(X)(\sum\limits_{i=1}^sm_i)(\prod\limits_{i=1}^s\psi_i(z_i)+{{\rm d}\over {\rm d}t}_{\mid_{t=0}}(\prod\limits_{i=1}^s\psi_i($exp$(tX_i)z_i)]$}

\vskip 2 pt

\centerline{$={{\rm d}\over {\rm d}t}_{\mid_{t=0}}(\prod\limits_{i=1}^sw_i^{m_i+k_i})e^{-t\alpha(X)(\sum\limits_{i=1}^sm_i)+t\alpha(X)}(\prod\limits_{i=1}^s\psi_i($exp$(tX_i)z_i)))]$}

\centerline{$-(\prod\limits_{i=1}^sw_i^{m_i+k_i}\Delta_i^{k_i}(z_i))[-\alpha(X)(\sum\limits_{i=1}^sm_i)(\prod\limits_{i=1}^s\psi_i(z_i)+{{\rm d}\over {\rm d}t}_{\mid_{t=0}}(\prod\limits_{i=1}^s\psi_i($exp$(tX_i)z_i)]$}

\vskip 2 pt

\centerline{$=\alpha(X){\cal M}(Q)\phi(w,z)={\cal M}([X,Q])\phi(w,z).$}

\vskip 2 pt
\noindent
2)  For $p\in {\goth p}_1$, there is $X_{-1}\in {\goth k}_{-1}$ such that $p=[X_{-1},Q]$, then 
\vskip 4 pt
\noindent
(i) for $X\in {\goth k}_1$,   since $[X,X_{-1}]$ belongs to ${\goth k}_0$ and $[{\rm d}\pi(X),{\cal M}(Q)]=0$, then 

\centerline{${\cal M}([X,p])={\cal M}([X,[X_{-1},Q]])={\cal M}([[X,X_{-1}],Q])$}

\vskip 2 pt

\centerline{$=[{\rm d}\pi([X,X_{-1}]),{\cal M}(Q)]=[[{\rm d}\pi(X),{\rm d}\pi(X_{-1}],{\cal M}(Q)]$}

\vskip 2 pt

\centerline{$=[{\rm d}\pi(X),[{\rm d}\pi(X_{-1}),{\cal M}(Q)]]-[{\rm d}\pi(X_{-1}),[{\rm d}\pi(X),{\cal M}(Q)]]$}

\vskip 2 pt

\centerline{$=[{\rm d}\pi(X),{\cal M}(p)]$}

\vskip 2 pt

\noindent
(ii) for $X\in {\goth k}_0$,  since $[X,X_{-1}]$ belongs to ${\goth k}_{-1}$, $[X,Q]=\alpha(X)Q$ and $[{\rm d}\pi(X),{\cal M}(Q)]=\alpha(X){\cal M}(Q)$ with $\alpha(X)\in {\bboard C}$ then
\vskip 2 pt

\centerline{${\cal M}([X,p])={\cal M}([X,[X_{-1},Q]])={\cal M}([[X,X_{-1}],Q])+{\cal M}([X_{-1},[X,Q]])$}

\vskip 2 pt

\centerline{$=[{\rm d}\pi([X,X_{-1}]),{\cal M}(Q)]+\alpha(X)[{\rm d}\pi(X_{-1}),{\cal M}(Q)]$}

\vskip 2 pt

\centerline{$=[[{\rm d}\pi(X),{\rm d}\pi(X_{-1})],{\cal M}(Q)]+\alpha(X)[{\rm d}\pi(X_{-1}),{\cal M}(Q)]=$}

\vskip 2 pt

\centerline{$[{\rm d}\pi(X),[{\rm d}\pi(X_{-1}),{\cal M}(Q)]]-[{\rm d}\pi(X_{-1}),[{\rm d}\pi(X),{\cal M}(Q)]]$}

\vskip 2 pt

\centerline{$+\alpha(X)[{\rm d}\pi(X_{-1}),{\cal M}(Q)]$}

\vskip 2 pt

\centerline{$=[{\rm d}\pi(X),[{\rm d}\pi(X_{-1}),{\cal M}(Q)]]=[{\rm d}\pi(X),{\cal M}(p)].$}

\vskip 2 pt
\noindent
(iii) for $X\in {\goth k}_{-1}$,  $[X,p]=[X,[X_{-1},Q]]$ then by definition, 

\vskip  2 pt

\centerline{${\cal M}([X,p])=[{\rm d}\pi(X),[{\rm d}\pi(X_{-1}),{\cal M}(Q)]]=[{\rm d}\pi(X),{\cal M}(p)]$.}

\vskip 2 pt
\hfill
\eject
\noindent
3) For  $p\in {\goth p}_{-1}$, there is $p'\in {\goth p}_1$ such that $p=\kappa(\sigma)p'$. Then,
\vskip 2 pt
\noindent
(i) for $X\in {\goth k}_1$,  there is $X^{\sigma}\in {\goth k}_{-1}$ such that

\vskip 2 pt

\centerline{$[X,p]=[X,\kappa(\sigma)p']={\rm d}\kappa(X)\kappa(\sigma)p'=\kappa(\sigma){\rm d}\kappa(X^{\sigma})p'$,}

\vskip 1 pt

\centerline{${\cal M}([X,p])=\pi(\sigma){\cal M}({\rm d}\kappa(X^{\sigma})p')\pi(\sigma)^{-1}=\pi(\sigma)[{\rm d}\pi(X^{\sigma}),{\cal M}(p')]\pi(\sigma)^{-1}$}

\vskip 1 pt

\centerline{$=[\pi(\sigma){\rm d}\pi(X^{\sigma})\pi(\sigma)^{-1},\pi(\sigma){\cal M}(p')\pi(\sigma)^{-1}]=[{\rm d}\pi(X),{\cal M}(p)].$}

\vskip 1 pt
\noindent
(ii) for $X\in {\goth k}_0$, there is $X^{\sigma}\in {\goth k}_0$ such that

\vskip 1 pt

\centerline{$[X,p]=[X,\kappa(\sigma)p']={\rm d}\kappa(X)\kappa(\sigma)p'=\kappa(\sigma){\rm d}\kappa(X^{\sigma})p'$.}

\vskip 1 pt
\noindent
Then, ${\cal M}([X,p])=\pi(\sigma)[{\rm d}\pi(X^{\sigma}),{\cal M}(p')]\pi(\sigma)^{-1}=[{\rm d}\pi(X),{\cal M}(p)].$

\vskip 1 pt
\noindent
(iii) for $X\in {\goth k}_{-1}$, there is $X^{\sigma}\in {\goth k}_{1}$ such that

\vskip 1 pt

\centerline{$[X,p]=[X,\kappa(\sigma)p']={\rm d}\kappa(X)\kappa(\sigma)p'=\kappa(\sigma){\rm d}\kappa(X^{\sigma})p'$.}

\vskip 1 pt
\noindent
Then, 
${\cal M}([X,p])=\pi(\sigma)[{\rm d}\pi(X^{\sigma}),{\cal M}(p')]\pi(\sigma)^{-1}=[{\rm d}\pi(X),{\cal M}(p)].$

\vskip 1 pt
\noindent
4) For  $p\in {\goth p}_{-2}$, i. e. up to a scalar, $p=1=\kappa(\sigma)Q$,  one has 
\vskip 1 pt
\noindent
(i) for $X\in {\goth k}_1$,  there is $X^{\sigma}\in {\goth k}_{-1}$ such that 

\centerline{$[X,1]=[X,\kappa(\sigma)Q]=\kappa(\sigma){\rm d}\kappa(X^{\sigma})Q$.}

\vskip 1 pt
\noindent
Then, 
${\cal M}([X,1])=\pi(\sigma)[d\pi(X^{\sigma}),{\cal M}(Q)]\pi(\sigma)^{-1}=[{\rm d}\pi(X),{\cal M}(1)].$
\vskip 1 pt
\noindent
(ii) for $X\in {\goth k}_0$, there is $X^{\sigma}\in {\goth k}_0$ such that

\centerline{$[X,1]=[X,\kappa(\sigma)Q]=\kappa(\sigma){\rm d}\kappa(X^{\sigma})Q$.}

\noindent
Then,
${\cal M}([X,1])=\pi(\sigma)[{\rm d}\pi(X^{\sigma}),{\cal M}(Q)]\pi(\sigma)^{-1}=[{\rm d}\pi(X),{\cal M}(1)].$
\vskip 1 pt
\noindent
5)  For $p\in {\goth p}_0$,   $p=[X_{-1},[X_{-1}',Q]]$ with $(X_{-1},X_{-1}')\in {\goth k}_{-1}\times {\goth k}_{-1}$.
\vskip 1 pt
(i) for $X\in {\goth k}_1$,  
$$\eqalign{&{\cal M}([X,p])={\cal M}([X,[X_{-1},[X_{-1}',Q]]])\cr
&={\cal M}([X,X_{-1}],[X_{-1}',Q]])+{\cal M}([X_{-1},[X,[X_{-1}',Q]])\cr
&=[{\rm d}\pi([X,X_{-1}]),{\cal M}([X_{-1}',Q])]+[{\rm d}\pi(X_{-1}'),{\cal M}([X,[X_{-1}',Q]])]\cr
&=[{\rm d}\pi(X),[{\rm d}\pi(X_{-1}),[{\rm d}\pi(X_{-1}'),{\cal M}(Q)]]]=[{\rm d}\pi(X),{\cal M}(p)].}$$
\vskip 1 pt
\noindent
(ii) for $X\in {\goth k}_0$,  
$$\eqalign{&{\cal M}([X,p])={\cal M}([[X,X_{-1}],[X_{-1}',Q]])+{\cal M}([X_{-1},[X,[X_{-1}',Q]])\cr
&=[{\rm d}\pi([X,X_{-1}]),{\cal M}([X_{-1}',Q])]\cr
&+[{\rm d}\pi(X_{-1}),{\cal M}([X,X_{-1}',Q]])]\cr
&=[([{\rm d}\pi(X),{\rm d}\pi(X_{-1})]),([{\rm d}\pi(X_{-1}'),{\cal M}(Q)])]\cr
&+[{\rm d}\pi(X_{-1}),([{\rm d}\pi(X),{\rm d}\pi(X_{-1}'),{\cal M}(Q)]])]\cr
&=[{\rm d}\pi(X),{\cal M}(p)].}$$
\vskip 1 pt
\noindent
(iii) for $X\in {\goth k}_{-1}$, there is a unique $Y\in {\goth k}_1$ such that $[X,p]=[Y,1]$.  Then  
${\cal M}([X,p])=[{\rm d}\pi(Y),{\cal M}(1)]$.
Since  ${\cal M}(p)=[{\rm d}\pi(X_{-1},[{\rm d}\pi(X_{-1}'),{\cal M}(Q)]]$ and 
${\cal M}(Q)(\prod\limits_{i=1}^sw_i^{m_i}\psi_i(z_i))=\prod\limits_{i=1}^sw_i^{m_i+k_i}\psi_i(z_i),$
one can deduce that  $[{\rm d}\pi(Y),{\cal M}(1)]=[{\rm d}\pi(X),{\cal M}(p)]$.
\quad \qed

\noindent
Hence we get a map $\rho :{\goth g}={\goth k}\oplus {\goth p}\to {\rm End}\bigl({{\cal O}}_{\rm fin}\bigr)$, 
where  
$$\rho (X):={\rm d}\pi (X) \quad (X\in {\goth k}), \quad  \quad \rho (p):={\cal M}(p)-\delta\circ{\cal D}(p)\quad (p\in {\goth p}).$$
\noindent
\vskip 4 pt
\noindent
Since the subspaces  ${{\cal O}}_{q,\rm fin}(\widetilde\Xi)$ are stable under $\rho$, we denote by $\rho_q$ the restriction of $\rho$ to ${{\cal O}}_{q,\rm fin}(\widetilde\Xi)$.  We  obtain:
 \th Theorem 3.8|
Assume that the constant $\eta_q$ of Theorem 3.4  exists. Fix $(\delta _{(k_1m+q_1,\ldots,k_sm+q_s)})$ as in Theorems 3.4.
\vskip 1 pt
\noindent
{\rm (i)} $\rho_q $ is a representation of the Lie algebra $\goth g$ on the space ${{\cal O}}_{q,\rm fin}(\widetilde\Xi)$.
\vskip 1 pt
\noindent
{\rm (ii)} The representation $\rho _q$ is irreducible.

\vskip 1 pt
\finth

\noindent
\proof {\rm (i)}   Since $\pi$ is a representation of $\widetilde K=\prod\limits_{i=1}^sK_i$, for $X,X'\in {\goth k},$ $\rho_q([X,X'])=[\rho_q(X),\rho_q(X')].$
It follows from Proposition 3.3 that for $X\in {\goth k}, p\in {\goth p},$ $\rho_q([X,p])=[\rho_q(X),\rho_q(p)].$
It remains to show that for $p,p'\in {\goth p}$, $\rho_q([p,p'])=[\rho_q(p),\rho_q(p')].$
From Theorem 3.4, $[\rho_q(E),\rho_q(F)]=\rho_q([E,F])$. Then  this follows from  [BK95  , Lemma 3.6].  Consider the map 
$\tau_q : \bigwedge^2{\goth p} \rightarrow $ End$({\cal O}_{q,{\rm fin}}),$ defined by 
$\tau_q(p\wedge p')=[\rho_q(p),\rho_q(p')]-\rho_q([p,p']).$ We know that $\tau_q(E\wedge F)=0$.  We have to show that $\tau_q\equiv 0$ on the space $\bigwedge^2{\goth p}$. Since  the group  $\widetilde K$ acts on ${\goth p}$ by the representation $\tilde\kappa=\otimes\kappa_i$ and  that $d\tilde\kappa=d\kappa$. It follows that  $\widetilde K$   acts on the space 
$\bigwedge^2{\goth p} $ by the representation (still denoted $\tilde\kappa$) given by 
$\tilde\kappa(g)(p\wedge p')=(\tilde\kappa(g)p)\wedge(\tilde\kappa(g)p')$
and by differentiation, the Lie algebra ${\goth k}$ acts on the space 
$\bigwedge^2{\goth p} $ by the representation (still denoted $d\kappa$)  given by $d\kappa(X)(p\wedge p')=({\rm d}\kappa(X)p)\wedge p'+p\wedge(d\kappa(X)p').$ 
Since the representation d$\kappa$ of ${\goth k}$ on ${\goth p}$ is irreducible and since  $E$ and $F$ are  highest weight and lowest  weight vectors  with respect to  ${\goth k}_0+{\goth k}_{1}$ and ${\goth k}_0+{\goth k}_{-1}$ respectively and since
${\goth p}_1=$d$\kappa({\goth k}_{-1})E,  {\goth p}_{-1}=$d$\kappa({\goth k}_{1})F$ and ${\goth p}_0=$d$\kappa({\goth k}_{-1})($d$\kappa({\goth k}_{-1})E=$d$\kappa({\goth k}_{1})($d$\kappa({\goth k}_{1})F$,
 it follows that 
$E\wedge F$ is cyclic for the action of ${\goth k}$ on $\bigwedge^2{\goth p} $ ,  i.e.  $\bigwedge^2{\goth p}={\cal U}({\goth k})(E\wedge F)$. 
\vskip 2 pt
\noindent
Furthermore,  using  Proposition 3.3,  one can write for  every  $X\in {\goth k}, p\in {\goth p}$, 
$\rho_q($d$\kappa(X)p)=[{\rm d}\pi(X),\rho_q(p)].$ Then,  for every $X\in {\goth k}, p,p'\in {\goth p}$,
\vskip 2 pt

\centerline{$\tau_q({\rm d}\kappa(X)(p\wedge p'))=\tau_q({\rm d}\kappa(X)p\wedge p')+\tau_q(p\wedge {\rm d}\kappa(X)p')$}

\vskip 2 pt

\centerline{$=[\rho_q({\rm d}\kappa(X)p),\rho_q(p')]-\rho_q([{\rm d}\kappa(X)p,p'])$}

\vskip 2 pt

\centerline{$+[\rho_q(p),\rho_q({\rm d}\kappa(X)p')]-\rho_q([p,{\rm d}\kappa(X)p']$}

\vskip 2 pt

\centerline{$=[[{\rm d}\pi(X),\rho_q(p)],\rho_q(p')]- \rho_q([[X,p],p'])$}

\vskip 2 pt

\centerline{$+[\rho_q(p),[{\rm d}\pi(X),\rho_q(p')]]-\rho_q([p,[X,p']])$}

\vskip 2 pt

\centerline{$=[{\rm d}\pi(X),[\rho_q(p),\rho_q(p')]]-\rho_q([X,[p,p']])$}

\vskip 2 pt

\centerline{$=[{\rm d}\pi(X),[\rho_q(p),\rho_q(p')]]-[{\rm d}\pi(X),\rho_q([p,p'])]$}

\vskip 2 pt

\centerline{$=[{\rm d}\pi(X),\tau_q(p\wedge p')].$}

\vskip 2 pt
\noindent
We deduce that for every $X\in {\goth k}$, $\tau_q({\rm d}\kappa(X)(E\wedge F))=0$ 
and, since $E\wedge F$ is cyclic, it follows that  $\tau_q\equiv 0$.

\vskip 1 pt
\noindent
{\rm (ii)}  Let ${\cal V}\ne \{0\}$ be a $\rho_q({\goth g})$-invariant  subspace of ${\cal O}_{q,{\rm fin}}(\widetilde\Xi)$. Then ${\cal V}$ is $\rho({\goth k})$-invariant.  As ${\cal O}_{q,{\rm fin}}(\widetilde\Xi)
=\sum\limits_{m\in {\bboard N}}{\cal O}_{(k_1m+q_1,\ldots,k_sm+q_s)}(\widetilde\Xi)$ and  as the subspaces ${\cal O}_{(k_1m+q_1,\ldots,k_sm+q_s)}(\widetilde\Xi)$  are 
$\rho_q({\goth k})$-irreducible, then there exists  ${\cal J}\subset {\bboard N}$ (${\cal J}\ne \emptyset$)  such that 
${\cal V}=\sum\limits_{m\in {\cal J}}{\cal O}_{(k_1m+q_1,\ldots,k_sm+q_s)}(\widetilde\Xi)$. Observe now that if ${\cal V}$ contains ${\cal O}_{(k_1m+q_1,\ldots,k_sm+q_s)}(\widetilde\Xi)$ 
then it contains 
${\cal O}_{(k_1(m+1)+q_1,\ldots,k_s(m+1)+q_s)}(\widetilde\Xi)$. In fact,  denote by $\phi_{m,q}=\otimes_{i=1}^s\phi_{k_im+q_i}$ with $\phi_{k_im+q_i}\in
 {\cal O}_{k_im+q_i}(\Xi_i)$ given by $\phi_{k_im+q_i}(w_i,z_i)=w_i^{k_im+q_i}$. As ${\cal D}\phi_{m,q}=0$, it follows that 
 $\rho_q(F)\phi_{m,q}={\cal M}\phi_{m,q}=\phi_{m+1,q}$ and that  $\rho_q(F)\phi_{m,q}$ belongs to ${\cal O}_{(k_1(m+1)+q_1,\ldots,k_s(m+1)+q_s)}(\widetilde\Xi)$, then the space ${\cal O}_{(k_1(m+1)+q_1,\ldots,k_s(m+1)+q_s)}(\widetilde\Xi)$ is included in $ {\cal V}.$ Let $m_0$ be 
the minimum of the $m$  such that ${\cal O}_{(k_1m+q_1,\ldots,k_sm+q_s)}(\widetilde\Xi) \subset {\cal V}$, then
${\cal V}=\sum\limits_{m=m_0}^{\infty}{\cal O}_{(k_1m+q_1,\ldots,k_sm+q_s)}(\widetilde\Xi).$
The function $\phi=\otimes_{i=1}^s\phi_i$ with $\phi_i(w_i,z_i)=\Delta_i^{k_im+q_i}w_i^{k_im+q_i}$ belongs to  ${\cal O}_{(k_1m+q_1,\ldots,k_sm+q_s)}(\widetilde\Xi)$ and 
for  $(w,z)=((w_1,z_1),\ldots,(w_s,z_s))$,
$$\rho_q(F)\phi(w,z)=\prod\limits_{i=1}^s\Delta_i^{(k_im+q_i)}(z_i)w_i^{k_im+q_i+p_i}$$
$$-\delta_{(k_1(m-1)+q_1,\ldots,k_s(m-1)+q_s)}\prod\limits_{i=1}^s\Delta_i^{k_i}\Bigl({\partial \over \partial z_i}\Bigr)\Delta_i^{k_im+q_i}(z_i)w_i^{k_im+q_i-p_i}.$$
\noindent
By the Bernstein identity (Proposition 2.6)
$$\Delta_i^{k_i}\Bigl({\partial \over \partial z_i}\Bigr) \Delta_i^{k_i(m+{q_i\over k_i}) } =B_i(m+{q_i\over k_i})\Delta_i^{k_i(m+{q_i\over k_i})-k_i} $$
\noindent
with 
$$B_i(m+{q_i\over k_i} )=b_i(k_im+q_i)b_i(k_im+q_i -1)\ldots b_i(k_im+q_i-k_i+1),$$
$b_i$ is the Bernstein polynomial relative to the determinant polynomial $\Delta _i$.
\vskip 2 pt
\noindent
Since $B_i(m+{q_i\over k_i}) >0$ for $m>0$, then,   if ${\cal O}_{(k_1m+q_1,\ldots,k_sm+q_s)}(\widetilde\Xi)\subset {\cal V}$, then ${\cal O}_{(k_1(m-1)+q_1,\ldots,k_s(m-1)+q_s)}(\widetilde\Xi) \subset {\cal V}$. Therefore $m_0=0$ and ${\cal V}={\cal O}_{q,{\rm fin}}(\widetilde\Xi)$.

\hfill \qed
\bigskip
Let us now consider,  in the following table, all the cases  where $\rho_q$ is a representation  of the Lie algebra $\goth g$.. In particular,  there is a  constant $\eta_q$ such that  for all $i$,
$$\eta_q={q_i\over k_i}+{n_i\over k_ir_i}.$$
\hfill \eject

\centerline{\bf Table 2} \rm

\vskip 10 pt
\noindent
(1)  $\widetilde{\rm SL}(3,{\bboard R})$ :  $V={\bboard C},  Q(z)=z^4$, $q\in \{0,1,2\}$,   $\eta_q={q\over 4}+{1\over 4}$.
\vskip 4pt
\noindent
(2)   $\widetilde{\rm SL}(p+2,{\bboard R})$ : $V={\bboard C}^p,  Q(z)=(z_1^2+\ldots+z_p^2)^2$,  $q\in \{0,1\}$, $\eta_q={q\over 2}+{p\over 4}$.
\vskip 4 pt
\noindent
(3)  $\widetilde{\rm SO}(3,3)$ : $V={\bboard C}\oplus{\bboard C}, Q(z)=z_1^2z_2^2, q_1=q_2=:\tilde q\in\{0,1\}, \eta_q={\tilde q\over 2}+{1\over 2}$.
\vskip 4pt
\noindent
(4) $\widetilde{\rm SO}(4,4)$ : $V=\oplus^4{\bboard C}$, $Q(z)=z_1z_2z_3z_4$, $q_i=\tilde q=0$,  $\eta_q=1$. 
\vskip 4 pt
\noindent
(5) $\widetilde{\rm SO}(p+2,3)$ :  $V={\bboard C}^p\oplus{\bboard C}, Q(z)=(z_1^2+\ldots+z_p^2){z'}^2$, $p$ odd

\vskip 4 pt

\centerline{$q_1=0$,  $q_2={p-1\over 2}$, $\eta_q={p\over 2}$.}

\vskip 4pt
\noindent
(6)  $\widetilde{\rm SO}(p+2,4)$ :  $V={\bboard C}^p\oplus {\bboard C}\oplus{\bboard C}$,  $Q(z)=(z_1^2+\ldots+z_p^2)z'z''$,  

\vskip 4 pt

\centerline{$p$  even, $q_1=0$,  $q_2=q_3={p\over 2}-1$,  $\eta_q={p\over 2}$.}

\vskip 4pt
\noindent
(7)  $\widetilde{\rm SO}(p_1+2,p_2+2)$ : $(V,Q(z))$ is the pair

\vskip 4 pt

\centerline{$({\bboard C}^{p_1}\oplus{\bboard C}^{p_2},(z_1^2+\ldots+z_{p_1}^2)({z_1'}^2+\ldots+{z_{p_2}'}^2)$,}

\vskip 4 pt

 $p_1-p_2\geq 0$ even,    $q_1=0$,  $q_2={p_1-p_2\over 2}$,   $\eta_q={p_1\over 2}$.
\vskip 4pt
\noindent
(8) $E_{6(6)}, E_{7(7)},E_{8(8)}$ : $(V,Q(z))$ is respectively the pair

\vskip 4 pt

\centerline{$({\rm Sym}(4,{\bboard C}), \det z)$ ;  $(M(4,{\bboard C}), \det z)$ ;  $({\rm Skew}(8,{\bboard C}),{\rm Pfaff}(z))$,}

 $q=0$,  $\eta_q=1+{3d\over 2}$, with $d=1, 2$ or $4$ respectively.
\vskip 4 pt
\noindent
(9) $F_{4(4)}, E_{6(2)},E_{7(-5)}, E_{8(-24)}$ :  $(V,Q(z))$ is respectively the pair

\vskip 4 pt

\centerline{$({\rm Sym}(3,{\bboard C})\oplus {\bboard C},\det(z_1)\cdot z_2)$}

\vskip 4 pt

\centerline{($M(3,{\bboard C})\oplus {\bboard C}, \det(z_1)\cdot z_2)$}

\vskip 4 pt

\centerline{(${\rm Skew}(6,{\bboard C})\oplus {\bboard C},{\rm Pfaff}(z_1)\cdot z_2$)}

\vskip 4 pt

\centerline{ 
$({\rm Herm}(3,{\bboard C})\oplus {\bboard C},\det(z_1)\cdot z_2)$,}

 $q_1=0$,    $q_2={n_1\over 3}-1$,  $\eta_0=1+d_1$, with $d_1=1, 2, 4, 8$.
\vskip 4 pt
\noindent
(10) $G_{2(2)} : V={\bboard C}\oplus{\bboard C}, Q=z_1^3z_2, q_1=2, q_2={q_1-2\over 3}=0, \eta_q={q_1\over 3}+{1\over 3}=1$.
\vskip 20 pt

\centerline{$(p,p_1,p_2\geq 2$)}

\vskip 1 pt

\vfill \eject
\vskip 1 pt
\noindent
\section 4. Irreducible unitary representations of  the   corresponding real  Lie group|
We consider, for a sequence $(c_{(k_1m+q_1,\ldots,k_sm+q_s)})$ of positive numbers, an inner product on 
${{\cal O}}_{q,\rm fin}(\widetilde\Xi)$  such that 
\vskip 4 pt

\centerline{$\Vert\phi\Vert^2=\sum\limits_{m\in {\bboard N}}{1\over c_{(k_1m+q_1,\ldots,k_sm+q_s)}}\Vert
\psi_{(k_1m+q_1,\ldots,k_sm+q_s)}\Vert_{(k_1m+q_1,\ldots,k_sm+q_s)}^2,$}

\vskip 2 pt
\noindent
for
\vskip 2 pt

\centerline{$\phi(w,z)=\sum\limits_{m\in{\bboard N}}\psi_{(k_1m+q_1,\ldots,k_sm+q_s)}(z_1,\ldots,z_s)w_1^{k_1m+q_1}\ldots w_s^{k_sm+q_s}.$}

\vskip 2 pt
\noindent
This inner product is invariant under $\widetilde K_{\bboard R}=\prod\limits_{i=1}^s(K_i)_{\bboard R}$.  We assume that the constant $\eta_0$  of Theorem 3.4. and the constants $\delta_{k_1m+q_1,\ldots,k_sm+q_s)}$ exist, and we  will determine the sequence
$(c_{(k_1m+q_1,\ldots,k_sm+q_s)})$ such that this inner product  is invariant under  the representation $\rho _q$ restricted to ${\goth g}_{\bboard R}$. We denote
by ${\cal F}_q(\widetilde\Xi)$ the Hilbert space completion of ${{\cal O}}_{q,\rm fin}(\widetilde\Xi)$ with
respect to this  inner product. We will assume $c_{(q_1,\ldots,q_s)}=1$.
\vskip 10 pt
Recall  from Proposition 2.6, the Bernstein identity
$$\Delta_i^{k_i}\Bigl({\partial \over \partial z_i}\Bigr)\Delta_i(z_i)^{k_im+q_i}=B_i(m+{q_i\over k_i})\Delta_i(z_i)^{k_im+q_i-k_i}.$$
where $$B_i(\alpha )=b_i(k_i\alpha )b_i(k_i\alpha -1)\ldots b_i(k_i\alpha -k_i+1),$$
$b_i$ is the Bernstein polynomial relative to the determinant polynomial $\Delta _i$, given by
$b_i(\alpha)=\alpha(\alpha+{d_i\over 2})\ldots(\alpha+(r_i-1){d_i\over 2})$
where $d_i$ is the degree of the simple summand $V_i$, (${n_i\over r_i}=1+(r_i-1){d_i\over 2}$).
The Bernstein polynomial $B(\alpha)=\prod\limits_{i=1}^sB_i(\alpha)$ is of degree 4, and vanishes at  $0$. The roots  of $B$ are: 

\vskip 2 pt

\centerline{${x\over k_i}-y{d_i\over 2k_i}$ with $1\leq i\leq s,  x, y \in{ \bboard N}, \quad 0\leq x\leq k_i-1, \quad 0\leq y\leq  r_i-1$. }

\vskip 6pt
\noindent
 We  have two cases:
\bigskip
First  case:  $B(1-\eta_q)=0$. This  is the case  if  $\exists i$, $q_i=0$. Let us denote by $a_0$ and $b_0$ the two other roots of $B(\alpha)$,  such  that $B(\alpha)=A\alpha(\alpha+\eta_q-1)(\alpha-a_0)(\alpha-b_0). \quad \quad $.
\vskip 2 pt
Second case:  $B(1-\eta_q)\ne 0$. This  is the case  if  $\forall i,  q_i\ne 0$. Let us denote by  $a_0', b_0',c_0'$ the three other zeros of $B$, in such a way that $B(\alpha )=A\alpha (\alpha -a_0')(\alpha -b_0')(\alpha -c_0').$
\vskip 1 pt

The remaining roots of the Bernstein polynomial  are given in  the following tables,  table 3 for the case 1, and table 4 for the case 2.

\vfill \eject

\centerline{\bf Table 3} \rm 
\vskip 10 pt

\def\tvi{\vrule height 12pt depth 1pt width 0pt}

\def\traithorizontal{\noalign{\hrule}}
$$\eqalignno{&\vbox{\offinterlineskip \halign{
\tvi #& \quad #&\quad # & \quad # & \quad # & \quad #  & \quad #  \cr
\qquad&$q$& $\eta_q $&$1-\eta_q$&$a_0$&$b_0$ \cr
\traithorizontal 
 \qquad(1)& $q=0$\qquad& ${1\over 4}$& ${3\over 4}$&${2\over 4} $ &${1\over 4}$ \cr
\qquad (2)&$q=0$& ${p\over 4}$& $1-{p\over 4}.$ & ${1\over 2} $ & ${1\over 2}-{p\over 4}$\cr
\qquad (3)&$q_1=q_2=0$& ${1\over 2}$ & ${1\over 2}$ & $0  $ & ${1\over 2} $ \cr
\qquad (4)& $q_1=q_2=q_3=q_4=0$& $1$ & $0$ & $0$ & $0$  \cr
\qquad (5)& $q_1=0$, $q_2=p-1$ & ${p\over 2}$& $1-{p\over 2}$ &${1\over 2}$  & $0$\cr
\qquad (6)& $q_1=0$, $p$ even, $q_2=q_3={p\over 2}-1$& ${p\over 2}$& $1-{p\over 2}$ \qquad &  $0$ & $0$ \cr
\qquad (7)& $q_1=0$, $p_1+p_2$ even, $q_2={p_1-p_2\over 2}$\qquad& ${p_1\over 2}$\qquad   & $1-{p_1\over 2}$ & $0$  & $1-{p_2\over 2}$\cr
\qquad (8)& $1+{3d\over 2}$ & $-{3d\over 2}$& $-2{d\over 2}$ & $-{d\over 2}$ \cr
\qquad (9)& $q_1=0$,  $q_2=d_1$& $1+d_1$& $-d_1$& $-{d_1\over 2}$ & $0$ \cr
\qquad (10)&$q_1=2$, $q_2=0$&$1$&$-q_2=0$&${1\over 3}$&${2\over 3} $\cr}}\qquad\qquad\qquad\cr}$$
\bigskip 

\vskip 40 pt

\centerline{\bf Table 4}\rm 

\vskip 10 pt

\def\tvi{\vrule height 12pt depth 1pt width 0pt}

\def\traithorizontal{\noalign{\hrule}}
$$\eqalignno{&\vbox{\offinterlineskip \halign{
\tvi # & \quad #&\quad# & \quad # & \quad # & \quad #  & \quad #  \cr
&$q$& $\eta_q $ & $c_0' $ & $a_0'$ & $b_0'$ \cr
\traithorizontal 
(1) & $q\in \{1,2\}$\qquad \qquad\qquad &${q\over 4}+{1\over 4}$ & ${3\over 4}$  & ${2\over 4} $ &${1\over 4}$ \cr
(2) &$q=1$, $p\ne 2$ \qquad\qquad& ${1\over 2}+{p\over 4}$& $1-{p\over 4}.$ & ${1\over 2} $ & ${1\over 2}-{p\over 4}$\cr
(3) &$\tilde q=1 $,  \quad $q_1=q_2=1$ \quad \qquad\qquad& ${\tilde q\over 2}+{1\over 2}=1$ & ${1\over 2}$ & $0  $ & ${1\over 2} $ \cr}}\qquad \qquad\qquad \cr}$$

\bigskip 
\vfill \eject
\th Theorem 4.1|
\vskip 1 pt
\noindent
Case 1 : Assume $B(1-\eta_q)=0$. Then 
\vskip 1 pt
\noindent
{\rm (i)} The inner product of ${\cal F}_q(\widetilde\Xi)$ is ${\goth g}_{\bboard R}$-invariant 
iff
$$c_{(k_1m+q_1,\ldots,k_sm+q_s)}={(\eta_q +1)_m\over (\eta_q+a_0)_m(\eta_q +b_0)_m}{1\over m!}.$$
{\rm (ii)} The  reproducing kernel of ${\cal F}_q(\widetilde\Xi)$ is given by 
$${\cal K}(\xi,\xi')={}_1F_2\bigl(\eta_q+1 ;\eta_q +a_0,\eta_q +b_0; \prod\limits_{i=1}^sH_i(z_i,z_i')(w_i\overline{w_i'})\bigr){\cal K}_q(\xi,\xi'),$$
\vskip 1 pt
\noindent
Case 2 : Assume $B(1-\eta_q)\ne 0$. Then 
\vskip 1 pt
\noindent
{\rm (i)} The inner product of ${\cal F}_q(\widetilde\Xi)$ is ${\goth g}_{\bboard R}$-invariant 
iff
$$c_{(k_1m+q_1,\ldots,k_sm+q_s)}={(\eta_q+1 )_m(1)_m\over(\eta_q+a_0')_m (\eta_q +b_0')_m(\eta_q +c_0')_m}{1\over m!}$$
{\rm (ii)} The  reproducing kernel of ${\cal F}_q(\widetilde\Xi)$ is given by 
$$
{\cal K}(\xi,\xi')=$$
$${}_2F_3\bigl(\eta_q+1,1; \eta_q+a_0',\eta_q +b_0',\eta_q +c_0';\prod\limits_{i=1}^sH_i(z_i,z_i')(w_i\overline{w_i'})\bigr){\cal K}_q(\xi,\xi')$$
where $\xi=((w_1,z_1),\ldots,(w_s,z_s))  , \xi'=((w_1',z_1'),\ldots,(w_s',z_s')),$
and 
\vskip 2 pt

\centerline{${\cal K}_q(\xi,\xi')=\prod\limits_{i=1}^sH_i^{q_i}(z_i,z_i')(w_i\overline{w_i'})^{q_i}$}

\vskip 2 pt
\noindent
is the reproducing kernel of the space ${{\cal O}}_{({q_1},\ldots,{q_s})}(\widetilde\Xi)$.
\vskip 1 pt
\noindent
(considering the values of $\eta_q+1, \eta_q+a_0', \eta_q+b_0', \eta_q+c_0'$ occurring in table 4, the function  ${}_2F_3$ becomes ${}_1F_2$).
 
 \finth

\vskip 1 pt
\proof
(i) Recall that ${\goth p}_{\bboard R}=\{p\in {\goth p}\mid \beta (p)=p\},$
where $\beta $ is the conjugation of $\goth p$, we introduced at the end of Section 1.
Recall also that  $\beta (\kappa (g)p)=\kappa \bigl(\alpha (g)\bigr)\beta (p).$
\vskip 2 pt
\noindent
The inner product of ${\cal F}_q(\widetilde\Xi)$ is ${\goth g}_{\bboard R}$-invariant if and only if, for every $p\in {\goth p}$,
$$\rho_q(p)^*=-\rho_q \bigl(\beta (p)\bigr).$$
But this is equivalent to the single condition $\rho_q(E)^*=-\rho _q(F)$
i.e, for 
$\phi\in {{\cal O}}_{(k_1(m+1)+q_1,\ldots,k_s(m+1)+q_s)}(\widetilde\Xi), \quad \phi'\in {{\cal O}}_{(k_1m+q_1,\ldots,k_sm+q_s)}(\widetilde\Xi)$,
$${1\over c_{(k_1(m+1)+q_1,\ldots,k_s(m+1)+q_s)}}(\phi\mid{\cal M}^{\sigma}\phi')
={\delta_{(k_1m+q_1,\ldots,k_sm+q_s)}\over c_{(k_1m+q_1,\ldots,k_sm+q_s)}}({\cal D}\phi\mid\phi').$$
\hfill
\eject
Recall that  the norm of $\psi=\otimes_{i=1}^s\psi_i \in \tilde{{\cal O}}_{(k_1m+q_1,\ldots,k_sm+q_s)}(\widetilde\Xi)$ can be written
$$\Vert \psi\Vert_{(k_1m+q_1,\ldots,k_sm+q_s)}^2$$
$$={1\over a_{(k_1m+q_1,\ldots,k_sm+q_s)}}\prod\limits_{i=1}^s\int_{V_i}\vert\psi_i(z_i)\vert^2H_i(z_i)^{-k_im-q_i-2{n_i\over r_i} }m(dz_i).$$
Then, the  required condition of invariance becomes 
$$\eqalignno{
&{1\over c_{(m+1)}a_{(m+1)}}\int_{\prod\limits_{i=1}^sV_i}\prod\limits_{i=1}^s\psi_i(z_i)
\Delta_i(z_i)^{k_i}
\overline{\psi_i'(z)}H_i(z_i)^{-k_i(m+1)-q_i-2{n_i\over r_i} }m(dz_i) \cr
&={\delta_{(m)}\over c_{(m+1)}a_{(m+1)}}
\int_{\prod\limits_{i=1}^sV_i}\prod\limits_{i=1}^s(\Delta_i^{k_i}\Bigl({\partial \over \partial z_i}\Bigr)\psi_i)(z_i)\overline{\psi_i'(z_i)}H_i(z_i)^{-k_im-q_i-2{n_i\over r_i} }m(dz_i).\cr}$$
where we denote 
$$x_{(m)}=x_{(k_1m+q_1,\ldots,k_sm+q_s)}.$$
By  integrating by parts:
$$\eqalign{
&\int_{\prod\limits_{i=1}^sV_i}\prod\limits_{i=1}^s(\Delta_i^{k_i}\Bigl({\partial\over \partial z_i}\Bigr)\psi_i)(z)\overline
{\psi'(z)}H_i(z_i)^{-k_im-q_i-2{n_i\over r_i} }m(dz_i) \cr
&=\int_{\prod\limits_{i=1}^sV_i}\prod\limits_{i=1}^s\psi_i(z_i)\overline{\psi_i'(z_i)}
\biggl(\Bigl(\Delta_i^{k_i}({\partial\over \partial z_i})\Bigr)H_i(z_i)^{-k_im-q_i-2{n_i\over r_i} }\Bigr)m(dz_i),\cr}$$
and, by the relation
$$\Delta_i^{k_i}\Bigl({\partial \over \partial z_i}\Bigr)H_i(z)^{-k_im-q_i-2{n_i\over r_i} } 
$$
$$=B_i(-m-{q_i\over k_i}-2{n_i\over k_i r_i} )\overline{\Delta_i^{k_i}(z_i)}H_i(z_i)^{-k_i(m+1)-q_i-2{n_i\over r_i} },$$
the invariance condition can be written 
$${c_{(m)}\over c_{(m+1)}}=
{a_{(m+1)}\over a_{(m)}}\delta _{(m)}\prod\limits_{i=1}^sB_i(-m-{q_i\over k_i}-2{n_i\over k_ir_i} ).$$

\vskip 2 pt
\noindent
From Proposition 2.7  it follows that
$${a_{(m+1)}\over a_{(m)}}={\prod\limits_{i=1}^sB_i(-m-{q_i\over k_i}-{n_i\over k_ir_i} )\over \prod\limits_{i=1}^sB_i(-m-{q_i\over k_i}-2{n_i\over  k_ir_i})}.$$

\vskip 2 pt
\noindent
Then we obtain 
$${c_{(m)}\over c_{(m+1)}}=\delta _{[m]}\prod\limits_{i=1}^sB_i(-m-{q_i\over k_i}-{n_i\over k_ir_i} ),$$
$$\eqalign{
&{c_{(m+1)}\over c_{(m)}}
={(m+\eta_q)(m+\eta_q +1)\over B(-m-\eta_0)}\cr
&={(m+\eta_q+1)\over(m+\eta_q+a_0')(m+\eta_q +b_0')(m+\eta_q +c_0')}.\cr}$$
Since $c_{(q_1,\ldots,q_s)}=1$,  we get 
\vskip 6 pt
\noindent
in case 1, i.e. if $c_0'=1-\eta_q, a_0'=a_0, b_0'=b_0$, then 
$$c_{(m)}={(\eta_q +1)_m\over (\eta_q +a_0 )_m(\eta_q +b_0)_m}{1\over m!}$$
\noindent
and in case 2,   we obtain :
$$c_{(m)}={(\eta_q +1)_m(1)_m\over(\eta_q+a_0')_m (\eta_q +b_0' )_m(\eta_q +c_0')_m}{1\over m!}$$
(ii) By Theorem 2.5 the reproducing kernel of ${\cal F}_q(\widetilde\Xi)$ is in case 1 given by
$$\eqalign{
&{\cal K}(\xi,\xi')
=\sum _{m\in{\bboard N}}
c_{(k_1m+q_1,\ldots,k_sm+q_s)} \prod\limits_{i=1}^sH_i(z_i,z_i')^{k_im+q_i}
(w_i\overline {w_i'})^{k_im+q_i} \cr
&={}_1F_2\bigl(\eta_q+1;\eta_q +a_0,\eta_q +b_0; \prod\limits_{i=1}^sH_i(z_i,z_i')(w_i\overline{w_i'})\bigr){\cal K}_q(\xi,\xi'),\cr}$$
\hfill 
\vskip 2 pt
\noindent
and in case 2 it is given by
$$\eqalign{
&{\cal K}(\xi,\xi')
=\sum _{m\in{\bboard N}}
c_{(k_1m+q_1,\ldots,k_sm+q_s)} \prod\limits_{i=1}^sH_i(z_i,z_i')^{k_im+q_i}
(w_i\overline {w_i'})^{k_im+q_i} =\cr
&{}_2F_3\bigl(1,\eta_q+1;\eta_q+a_0',\eta_q +b_0',\eta_q +c_0'; \prod\limits_{i=1}^sH_i(z_i,z_i')(w_i\overline{w_i'})\bigr){\cal K}_q(\xi,\xi'),\cr}$$

\centerline{($\xi=((w_1,z_1),\ldots,(w_s,z_s))$, $\xi'=((w_1',z_1'),\ldots,(w_s',z_s'))$).}
\hfill \qed 
\vskip 2 pt

\hfill
\eject
The  formulas  obtained in Theorem 4.1, for the reproducing kernel  of the Hilbert space ${\cal F}_q(\widetilde\Xi)$,  by considering  the cases  in table 3 and  table 4,  agree exactly  with  the formulas given by Theorem 8.1 and table 6.9 in [B98]  for the minimal representations (one has to replace $r_0$ by $\eta_q$, $a$ by $\eta_q+a_0$, $b$ by $\eta_q+b_0$) and in [B97] for the $SL(3,{\bboard R})$-case .

\vskip 4 pt
\noindent
In the special case of $\widetilde{\rm SL}(p+2,{\bboard R})$ (with $p\geq 3$), the reproducing kernels of   the  two  representations corresponding to $q=0$ and $1$  are respectively
\vskip 4 pt

\centerline{${}_1F_2({5\over 4}; {3\over 4}, {1\over 2} ; \phi_p(1+z\overline{z'})(w\overline{w'})\bigr),$}

\vskip 4 pt

\centerline{${}_1F_2({p\over 4}+{3\over 2}; {3\over 2}, {p\over 4} ; \phi_p(1+z\overline{z'})(w\overline{w'})\bigr)$}

\vskip 4 pt
\noindent
where $\xi=(w,z)\in{\bboard C}\times {\bboard C}^p)$ and 
${\cal K}_q(\xi,\xi')=[\phi_p(1+z\overline{z'})(w\overline{w'})]^q.$
\vskip 4 pt
\noindent
In the special case of $\widetilde{\rm SO}(3,3)\simeq \widetilde{\rm SL}(4,{\bboard R})$, the reproducing kernels of   the  two  representations corresponding to $q_1=q_2=0$ and $q_1=q_2=1$  are respectively given by
\vskip 4 pt

\centerline{${}_1F_2({3\over 2}; {1\over 2}, 1 ; (1+z_1\overline{z_1'})
(1+z_2\overline{z_2'})(w_1\overline{w_1'}(w_2\overline{w_2'})\bigr),$}

\vskip 4 pt

\centerline{${}_1F_2(2; {3\over 2}, {3\over 2} ; (1+z_1\overline{z_1'})
(1+z_2\overline{z_2'})(w_1\overline{w_1'}(w_2\overline{w_2'}]\bigr)$}

\vskip 4 pt
\noindent
where 
$\xi=(w,z)=((w_1,z_1),(w_2,z_2)), \xi'=(w',z')=((w_1',z_1'),(w_2',z_2'))\in ({\bboard C}\times {\bboard C})^2$
 and 
${\cal K}_{(q_1,q_2)}(\xi,\xi')=[(1+z_1\overline{z_1'})(w_1\overline{w_1'})]^{q_1}
[(1+z_2\overline{z_2'})(w_2\overline{w_2'})]]^{q_2}.$

\vskip 2 pt

\noindent
One can observe  that the case  ((1), $q=0$) agrees with the case ((2), $q=0$, $p=1$),  the case ((1), $q=2$) agrees with the case ((2), $q=1$, $p=1$), 
the case  ((3), $q_1=q_2=0$) agrees with the case ((2), $q=0$, $p=2$) and the  case  ((3), $q_1=q_2=1$) agrees with the case ((2), $q=1$, $p=2$). 
\hfill
\eject
\bigskip
In the following, we will see that the Hilbert space ${\cal F}_q(\widetilde\Xi)$ is a pseudo-weighted Bergman space. It means
that the norm of $\phi\in {\cal F}_q(\widetilde\Xi)$  is given by an integral of $|\phi |^2$ with respect to a weight taking both positive and negative values.
The weight involves a Meijer $G$-function:

$$G(u)=G_{2,4}^{4,0}\Bigl(u\big|\matrix{\eta_q&\eta_q-1 & \cr \beta_1 &\beta_2 & \beta_3 & \beta_4 \cr}\Bigr).$$

 \bigskip
 \th Theorem 4.2|
For $\phi \in {\cal F}_q(\widetilde\Xi)$,
$$\|\phi \|^2=\int _{\prod\limits_{i=1}^s({\bboard C}\times V_i)}
\phi (w,z)|^2p(w,z)\prod\limits_{i=1}^sm(dw_i)m_{i,0}(dz_i),$$
with 
$$p(w,z)=CG\bigl(\prod\limits_{i=1}^s|w_i|^2H_i(z_i)\bigr) \prod\limits_{i=1}^sH_i(z_i).$$
\finth
\vskip 4 pt

\noindent
In fact,  from the proof of Theorem 4.1 it follows that
$${c_{(m)} a_{(m)}\over c_{(m+1)}a_{(+1)}}=\delta _{(m)}
\prod\limits_{i=1}^sB_i(-m-{q_i\over k_i}-2{n_i\over  k_ir_i} )={\prod\limits_{i=1}^sB_i(-m-\eta_q-{n_i\over k_ir_i})\over(m+\eta_q)(m+\eta_q+1)}.$$
\vskip 2 pt
\noindent
Then
$$\eqalignno{&{1\over c_{(m+1)}a_{(m+1)}}=
{\prod\limits_{i=1}^sB_i(-m-\eta_q-{n_i\over  k_ir_i} )
\over(m+\eta_q)(m+\eta_q+1)}{1\over (c\cdot a)_{(k_1m+q_1,\ldots,k_sm+q_s)}}\cr}$$
$$\eqalignno{&={(m+\eta_q+\alpha_1')(m+\eta_q+\alpha_2')(m+\eta_q+\alpha_3')(m+\eta_q+\alpha_4')\over (m+\eta_q)(m+\eta_q+1)}{1\over c_{(m)}a_{(m)}}\cr}$$
where 
$$\tilde B(\alpha)=\prod\limits_{i=1}^sB_i(\alpha-{n_i\over  k_ir_i})=A(\alpha-\alpha_1')(\alpha-\alpha_2')(\alpha-\alpha_3')(\alpha-\alpha_4').$$
  Then
  \hfill
  \eject
 
$${1\over c_{(m)}a_{(m)}}={(\eta_q+\alpha_1')_m(\eta_q+\alpha_2')_m(\eta_q+\alpha_3')_m(\eta_q+\alpha_4')_m\over(\eta_q)_m(\eta_q+1)_m}$$
$$={\Gamma(\eta_q+\alpha_1'+m)\Gamma(\eta_q+\alpha_2'+m)\Gamma(\eta_q+\alpha_3'+m)\Gamma(\eta_q+\alpha_4'+m)\over\Gamma(\eta_q+m)\Gamma(\eta_q+1+m)}$$
$$={\Gamma(\eta_q+\alpha_1'+m)\Gamma(\eta_q+\alpha_2'+m)\Gamma(\eta_q+\alpha_3'+m)\Gamma(\eta_q+\alpha_4'+m)\over\Gamma(\eta_q+m)\Gamma(\eta_q+1+m)}$$

$$=C\int_0^{\infty}G(u)u^mdu$$
where 
$$G(u)=G_{2,4}^{4,0}\Bigl(u\big|\matrix{\eta_q&\eta_q-1 & \cr \beta_1 &\beta_2 & \beta_3 & \beta_4 \cr}\Bigr).$$

$$\beta_1=\eta_q+\alpha_1'-1, \quad \beta_2=\eta_q+\alpha_2'-1, \quad \beta_3=\eta_q+\alpha_3'-1, \quad \beta_4=\eta_q+\alpha_4'-1.$$
\vskip 1 pt
Observe that  in the particular case $q=0$,  $\tilde B(\alpha)=B(\alpha-\eta_q)$ and,  since  the Bernstein polynomial $B$   vanishes   at  $0$ and $1-\eta_q$,  then 

\vskip 2 pt

\centerline{$\alpha_1'=\eta_q, \quad  \alpha_2'=1,  \quad \alpha_3'=\eta_q+a_0, \quad  \alpha_4'=\eta_q+b_0$}

then

\centerline{$\beta_1=2\eta_q-1, \quad \beta_2=\eta_q, \quad \beta_3=2\eta_q+a_0-1$ and $ \beta_4=2\eta_q+b_0-1$}

\vskip 2 pt
\noindent
and  the $G$-function  is given by 
$$G(u)=G_{1,3}^{3,0}\Bigl(u\big|\matrix{\eta_q-1 & \cr2\eta_q-1& 2\eta_q+a_0-1&2\eta_q+b_0-1 \cr}\Bigr), $$
and these formulas agree with those obtained in [AF12].
\bigskip
\hfill \eject 
\centerline{\bf Table 5} \rm

\vskip 10 pt

$$\eqalignno{&\vbox{\offinterlineskip \halign{
\tvi # & # &  # &  # & # &\quad #\cr
&$$& $\alpha_1'$ &$\alpha_2'$& $\alpha_3' $ & $\alpha_4'$  \cr
\traithorizontal 
& (1)  &${1\over 4}$ & ${2\over 4}$ & ${3\over 4}$ & $1$ \cr
& (2) &${p\over 4}$ & ${1\over 2}$ & ${p\over 4}+{1\over 2}$ & $1$ \cr
& (3)  &${1\over 2}$& $1$ & ${1\over 2}$ & $1$ \cr
& (4)  &$1$ & $1$ & $1$ & $1$ \cr
&  (5) & ${p\over 2}$ & $1$ & ${1\over 2}$ & $1$ \cr
& (6) & ${p\over 2}$ & $1$ & $1$ & $1$\cr
 & (7)  & ${p_1\over 2}$ & $1 $ & ${p_2\over 2}$ & $1$ \cr
&  (8) & $1+{3d\over 2}$ & $1+d$ & $1+{d\over 2}$ & $1$ \cr
 & (9)  & $1+d_1$ & $1+{d_1\over 2}$ & $1$ & $1$\cr
&  (10) & ${1\over 3}$ & ${2\over 3}$ & $1$  & $1$ \cr}}\cr}$$
\bigskip

\centerline{\bf Table 6} \rm

\vskip 10 pt

$$\eqalignno{&\vbox{\offinterlineskip \halign{
\tvi#& # & # &  # &  # & # & #&\quad #\cr
&$$& $\eta_q-1$ &$\eta_q$& $\beta _1 $ & $\beta _2$ & $\beta _3$&$\beta_4$ \cr
\traithorizontal 
 & (1) & ${q\over 4}-{3\over 4}, q\ne 3$ & ${q\over 4}+{1\over 4}$ & ${q\over 4}-{1\over 2}$ & ${q\over 4}-{1\over 4}$ & ${q\over 4}$& ${q\over 4}+{1\over 4}$\cr
 &(2)& ${q\over 2}+{p\over 4}-1,q\in \{0,1\}$ & ${q\over 2}+{p\over 4}$ & ${q\over 2}+{p\over 2}-1$ & ${q\over 2}+{p\over 4}-{1\over 2}$&${q\over 2}+{p\over 2}-{1\over 2}$&${q\over 2}+{p\over 4}$ \cr
 &(3) &${\tilde q\over 2}-{1\over 2}, \tilde q\in\{0,1\}$ & ${\tilde q\over 2}+{1\over 2}$ & ${\tilde q\over 2}$ & ${\tilde q\over 2}$ &${\tilde q\over 2}+{1\over 2}$&${\tilde q\over 2}+{1\over 2}$\cr
 &(4) &$0$ & $1$ & $1$ & $1$ &$1$&$1$\cr
 &(5)& ${p\over 2}-1$ & ${p\over 2}$ & $p-1$ & ${p\over 2}$&${p\over 2}-{1\over 2}$&${p\over 2}$ \cr
 &(6)& ${p\over 2}-1$ & ${p\over 2}$ & $p-1$ & ${p\over 2}$&${p\over 2}$&${p\over 2}$ \cr
 &(7) &${p_1\over 2}-1$ & ${p_1\over 2}$ & $p_1-1$ & ${p_1\over 2}$ &${p_1+p_2\over 2}-1$&${p_1\over 2}$\cr
 &(8)& ${3d\over 2}$ & ${3d\over 2}+1$ & $3d+1$ & ${5d\over 2}+1$ & $2d+1$& ${3d\over 2}+1$\cr
 &(9) &$d_1$ & $d_1+1$ & $2d_1+1$  & ${3d_1\over 2}+1$  & $d_1+1$ & $d_1+1$\cr
 &  (10) & $0$ & $1$ & ${1\over 3}$&${2\over 3} $ & $1$&$1$ \cr}}\cr}$$
\bigskip
Observe  that one of  the $\alpha_i'$  always equals $1$, and then  one of the $\beta_i$ always equals $\eta_q$, therefore the  involved $G$-function is always $G_{1,3}^{3,0}$.

\bigskip
\hfill \eject
Let ${G_{\bboard R}}$ be the simply
connected Lie group with Lie algebra ${\goth g}_{\bboard R}.$

\vskip 1 pt
\th Theorem 4.3| 
{\rm (i)} There is a unique unitary irreducible
representation $\pi_q$  of $G_{\bboard R}$ on ${\cal F}_q(\widetilde\Xi)$ such that
d$\pi_q=\rho_q $.
\vskip 1 pt
{\rm (ii)} The representation $\pi_q$ is spherical if and only if $q=0$. In that case,  the space ${{\cal O}}_{(0,\ldots,0)}(\widetilde\Xi)$ reduces to the constant functions,  the $K$-fixed vectors.
\finth
The proof of (i)  is the same than the proof of Theorem 6.3 in [AF12]. The (ii) is rather obvious.
\vskip 6 pt

\noindent
\section 5. The ${\goth sl}(3,{\bboard R})$-case|
\vskip 4 pt
In the special case $s=1, r=1, k=4$,  the polynomials $\Delta$ and  $Q$  are given by $\Delta(z)=z$ and $Q(z)=z^4$. We denote by $E=Q\in {\cal V}:={\goth p}_2$ and $ F=1\in {\cal V}^{\sigma}:={\goth p}_{-2}$. Observe that the ${\goth k}_i$ and the ${\goth p}_i$ are one dimensional. For $H\in {\goth k}_0$ such that exp$(tH)$ is the dilation  : $z\in {\bboard C} \mapsto e^{-t}z \in{\bboard C}$, we have $[H, E]=2E$ and $[H, F]=-2F$, 
and  we  consider the Lie algebra structure on ${\goth g}:={\goth k}\oplus {\goth p}$ such that $[E,F]=H$. In this case ${\goth g}$ is isomorphic to ${\goth sl}_3({\bboard C})$ and the real form ${\goth g}_{{\bboard R}}$ is isomorphic to ${\goth sl}(3,{\bboard R}).$ The structure group of $V={\bboard C}$ is $L={\rm Str}(V,\Delta)={\bboard C}^*$ acting by dilations $l_{\lambda}$,  and,   since for $\lambda\in {\bboard C}^*, \Delta(\lambda\cdot z)=\lambda\Delta(z)$,  the character $\chi$ doesn't exist , then  $K=\{(g,\mu) \mid g\in {\rm Conf}(V,\Delta)$, $\mu(z)^2=z\}$, with ${\rm Conf}(V,\Delta)={\rm SL}(2,{\bboard C})/\{\pm I\}$. The orbit $\Xi=\Xi^{(4)}$ of  $Q$  under $K$ acting by $\kappa^{(4)}$ and the orbit $\widetilde\Xi$ of  $\Delta$  under $K$ acting by $\kappa$ are  given by
$$\Xi=\{w^4\cdot(v-z)^4\mid (w,z)\in {\bboard C}\times {\bboard C}\},$$
\noindent
$$\widetilde\Xi=\{w\cdot(v-z)\mid (w,z)\in {\bboard C}\times {\bboard C}\}$$
\noindent
in such a way that the map $\widetilde\Xi \rightarrow \Xi, \xi \mapsto \xi^4$ is a covering and $\Xi$ is diffeomorphic to $\widetilde\Xi/Z_4$, where $Z_4$ is the group of the 4-roots of unity.
\vskip 4 pt
\noindent
The space  ${\cal O}_{q,{\rm fin}}(\widetilde\Xi)$ is given by 
${\cal O}_{q,{\rm fin}}(\widetilde\Xi)=\sum\limits_{m\in {\bboard N}}{\cal O}_{4m+q}(\widetilde\Xi)$, 

\centerline{${\cal O}_{4m+q}(\widetilde\Xi)=\{f\in {\cal O}(\widetilde\Xi) \mid f(\lambda \xi)=\lambda^{4m+q}f(\xi)\}$}

\vskip 2 pt

\centerline{$\simeq\{\phi \in  {\cal O}({\bboard C}^2) \mid \phi(w,z)=w^{4m+q}\psi(z)\}$.}

\vskip 2 pt
\noindent
 The norm  on the space  ${{\cal O}}_{4m+q}({\widetilde\Xi})$ is  given for $\phi(w,z)=w^{4m+q}\psi(z)$ by 
\vskip 2 pt

\centerline{$\Vert \phi\Vert_{4m+q}^2={1\over a_{4m+q}}\int_{{\bboard C}}\vert \psi(z)\vert(1+\vert z\vert^2)^{-(4m+q+2)}dz,$}

\vskip 2 pt

\centerline{$a_{4m+q}=\int_{{\bboard C}}(1+\vert z\vert^2)^{-(4m+q+2)}dz=\pi{1\over 4m+q+1}.$}

\vskip 4 pt
\noindent
In particular, for 
$$\Phi_{4m+q}(w,z):=w^{4m+q}$$
\noindent
one has 
$$\Vert\Phi_{4m+q}\Vert_{4m+q}=1.$$
\hfill
\eject
\noindent
The reproducing kernel  of the space ($\widetilde{{\cal O}}_{4m+q}({\bboard C}), \Vert\cdot\Vert_{4m+q})$  is given by 
$$\tilde{\cal K}(z,z')=(1+z\bar z')^{4m+q}.$$
\vskip 4 pt
\noindent
The representation $\pi$  of $K$ on ${{\cal O}}(\widetilde\Xi)={\cal O}({\bboard C}^2)$  is given by
$$(\pi(g)\phi)(w,z)=\phi(\mu(g^{-1},z)w,g^{-1}\cdot z).$$
It follows  that  for $\lambda=e^{-t}$,  $(\pi(l_{\lambda})\phi)(w,z)=\phi(e^{-2t}w,e^{t}\cdot z)$, 
$$(d\pi(H)\phi)(w,z)={d\over dt}_{\mid_{t=0}}\phi(e^{-2t}w,e^{t}\cdot z)=-2\phi(w,z)+({\cal E}\phi)(w,z).$$

\vskip 4 pt
\noindent
The operators ${\cal M}$, $ {\cal D}$, ${\cal M}^{\sigma}$, ${\cal D}^{\sigma}$ are given by
$$({\cal M}\phi)(w,z)=w^4\phi(w,z),  ({\cal D}\phi)(w,z)=w^{-4}{\partial^4\over\partial z^4}\phi(w,z).$$
$${\cal M}^{\sigma}(w^{4m+q}\psi(z))$$
$$=w^{4(m+1)+q}z^4\phi(w,z) ,  ({\cal D}^{\sigma}\phi)(w,z)=w^{4(m-1)+q}({\bf D}_{-4m-q})^*\psi(z)$$
($({\bf D}_{-4m-q})^*$ is the adjoint of the operator 
${\bf D}_{-4m-q}=z^{4-4m-q}{\partial^4\over\partial z^4}z^{4m+q}.$)
\vskip 8 pt
\noindent
The representation $\rho_q$ is given  by 
$$(\rho_q(E)(w^{4m+q}\psi(z))=$$
$$w^{4(m+1)+q}\psi(z)-\delta_{4(m-1)+q}w^{4(m-1)+q}{\partial^4\over\partial z^4}\psi(z),$$
$$\rho_q(F)(w^{4m+q}\psi(z))=$$
$$w^{4(m+1)+q}z^4\psi(z)-\delta_{4(m-1)+q}w^{(4(m-1)+q}({\bf D}_{-4m-q})^*\psi(z),$$
$$\rho_q(H)(w^{4m+q}\psi(z))=d\pi(H)(w^{4m+q}\psi(z))=w^{4m+q}(-2\psi+\tilde{\cal E}\psi)(z)$$
with
$$\tilde{\cal E}\psi)(z)={d\over dt}_{\mid_{t=0}}\psi(e^tz)$$
and
$$\delta_{4m+q}={1\over 4^4(m+{q\over 4}+{1\over 4}) (m+{q\over 4}+{1\over 4}+1)}.$$
\bigskip
The invariance condition $\rho_q(F)^*=-\rho_q(E)$ is equivalent to the   following:
$${1\over a_{4(m+1)+q}c_{4(m+1)+q}}\int_{{\bboard C}}z^4\psi(z)\overline{\psi'(z)}(1+\vert z\vert^2)^{-(4(m+1)+q+2)}dz=$$
$$-{\delta_{4m+q}\over a_{4m+q}c_{4m+q}}\int_{{\bboard C}}\psi(z){\partial^4\over\partial z^4}\overline{\psi'(z)}(1+\vert z\vert^2)^{-(4m+q+2)}dz$$
which is  equivalent to
$${1\over a_{4(m+1)+q}c_{4(m+1)+q}}\int_{{\bboard C}}z^4\psi(z)\overline{\psi'(z)}(1+\vert z\vert^2)^{-(4m+q+6)}dz=$$
$$
-{\delta_{4m+q}\over a_{4m+q}c_{4m+q}}\int_{{\bboard C}}\psi(z)\overline{\psi'(z)}{\partial^4\over\partial z^4}(1+\vert z\vert^2)^{-(4m+q+2)}dz.$$
Since $${\partial^4\over\partial\bar z^4}(1+\vert z\vert^2)^{-(4m+q+2)}$$ 
$$=(4m+q+2)(4m+q+3)(4m+q+4)(4m+q+5)z^4(1+\vert z\vert^2)^{-(4(m+1)+q+2)},$$
$$\delta_{4m+q}={1\over 4^4(m+{q\over 4}+{1\over 4})(m+{q\over 4}+{1\over 4}+1)},$$
and $a_m={\pi\over 4m+q+1}$, it follows, after an integration by parts,  that 
$${1\over a_{4(m+1)+q}c_{4(m+1)+q}}$$
$$=
{(4m+q+2)(4m+q+3)(4m+q+4)(4m+q+5)\over 4^4(m+{q\over 4}+{1\over 4})(m+{q\over 4}+{1\over 4}+1)}{1\over a_{4m+q}c_{4m+q}}$$
i.e. 
$$ (4m+q+5){1\over c_{4(m+1)+q}}$$
$$\eqalignno{={(4m+q+2)(4m+q+3)(4m+q+4)(4m+q+5)(4m+q+1)\over 4^4(m+{q\over 4}+{1\over 4})(m+{q\over 4}+{1\over 4}+1)}{1\over c_{4m+q}}}$$
i.e.
$${1\over c_{4(m+1)+q}}={(m+{q\over 4}+{2\over 4})(m+{q\over 4}+{3\over 4})(m+{q\over 4}+1)\over m+1+{q\over 4}+{1\over 4}} {1\over c_{4m+q}}$$
and for $c_q=1$, one obtains 

\vskip 4 pt

\centerline{ ($q=0$) : \quad $c_{4m}={({5\over 4})_m\over ({2\over 4})_m({3\over 4})_m}{1\over m!}$\quad ; \quad 
($q=1$) :  $c_{4m+1}={({6\over 4})_m\over ({3\over 4})_m({5\over 4})_m}{1\over m!}$,}

 \vskip 2 pt
 
 \centerline{($q=2$) :  $c_{4m+2}={({7\over 4})_m\over ({5\over 4})_m({6\over 4})_m}{1\over m!}$.}

\vskip 4 pt
The representation $d\pi+\rho_q$ integrates to a unitary representation of $\widetilde{\rm SL}(3,{\bboard R})$ in 
 a Hilbert space  ${\cal F}_q(\widetilde\Xi)$, the completion of ${\cal O}_{q,{\rm fin}}(\widetilde\Xi)$ with respect to the norm  given for 
 $\phi=\sum\limits_{m\in {\bboard N}}\phi_{4m+q}$,  by 
 $\Vert \phi\Vert^2=\sum\limits_{m\in {\bboard N}}{1\over c_{4m+q}}\Vert\phi_{4m+q}\Vert _{4m+q}^2,$
 and the reproducing kernels  of ${\cal F}_q(\widetilde\Xi)$ for $q=0, 1, 2 $ are respectively   given by:

$${}_1F_2({5\over 4}; {3\over 4}, {1\over 2} ; (1+z\overline{z'})(w_i\overline{w_i'})\bigr){\cal K}_0(\xi,\xi')\quad (q=0)$$
$${}_1F_2({3\over 2}; {5\over 4}, {3\over 4} ; (1+z\overline{z'})(w_i\overline{w_i'})\bigr){\cal K}_1(\xi,\xi')(q=1)$$
$${}_1F_2({7\over 4}; {6\over 4}, {5\over 4} ; (1+z\overline{z'})(w\overline{w'})\bigr){\cal K}_2(\xi,\xi')\quad (q=2),$$

where $\xi=(w,z), \xi'=(w',z')\in {\bboard C}^2$ and ${\cal K}_q(\xi,\xi')=[(1+z\overline{z'})(w\overline{w'})]^q$.
\vskip 4 pt

 Furthermore the Hilbert space ${\cal F}_q(\widetilde\Xi)$ is a weighted Bergman space and its norm is given by
$$\|\phi \|^2=\int _{{\bboard C}^2}
|\phi (w,z)|^2p(w,z)m(dw)(1+\vert z\vert^2)^{-2}(dz),$$
with $p(w,z)=G\bigl(\vert w\vert^2(1+\vert z\vert^2)\bigr)(1+\vert z\vert^2))$ and  the $G$-function is given by

$$G(u)=G_{1,3}^{3,0}\Bigl(u\big|\matrix{-{3\over 4}& \cr -{2\over 4} && 0& -{1\over 4} \cr}\Bigr)\quad (q=0),$$

$$G(u)=G_{1,3}^{3,0}\Bigl(u\big|\matrix{-{1\over 2}& \cr -{1\over 4} &0 & {1\over 4}  \cr}\Bigr)\quad (q=1),$$

$$G(u)=G_{1,3}^{3,0}\Bigl(u\big|\matrix{-{1\over 4} & \cr 0 &{1\over 4}&-{1\over 2}  \cr}\Bigr)\quad (q=2).$$

\bigskip
\bigskip
\bigskip
\bigskip
\hfill \eject
\centerline{\bf References}
\vskip 6pt
\noindent
\rm
\vskip 4 pt
\noindent
[A11] D. Achab (2011), \it Construction process for simple Lie algebras, \rm
J.  of Algebra, \bf 325,\rm 186-204.
\vskip 4 pt
\noindent
[AF12] D. Achab and J. Faraut (2012), \it Analysis of the Brylinski-Kostant model for minimal representations, Canad.J. of Math.,\bf 64, \rm721-754.
\vskip 4 pt
\noindent
[B97] R. Brylinski (1997), \it Quantization of the 
4-dimensional nilpotent orbit of $SL(3,{\bboard R})$, \rm 
Canad.  Jou.l of Math.,\bf 49,\rm 916-943.
\vskip 4 pt
\noindent
[B98]  R. Brylinski (1998), \it Geometric quantization of real minimal nilpotent orbits, \rm Symp. Geom.,
Diff. Geom.  and Appl., \rm 9, 5-58.
\vskip 4pt
\noindent
[BK94] R. Brylinski and B. Kostant (1994), \it 
Minimal representations, geometric quantization and unitarity, \rm 
Proc. Nat. Acad., \bf 91, \rm 6026-6029.
\vskip 4 pt
\noindent
[DS99-I] A. Dvorsky and S. Sahi (1999), \it Explicit Hilbert spaces for certain unipotent representations II, \rm 
Invent. Math., \bf 138, no. 1, \rm 203-224.
\vskip 4 pt
\noindent
[DS99-II] A. Dvorsky and S. Sahi (1999), \it Explicit Hilbert spaces for certain unipotent representations III, \rm 
J. Funct. Anal.\bf 201, no. 2, \rm 430-456.
\vskip 4 pt
\noindent
[FK94]  J. Faraut and A. Kor\'anyi (1994), \it Analysis on symmetric cones, \rm Oxford University Press.
\vskip 4 pt
\noindent
[FG96]  J. Faraut and S. Gindikin (1996), \it Pseudo-Hermitian symmetric spaces of
tube type, in  
{\rm Topics in Geometry (S. Gindikin ed.)}., \rm Progress in non linear differential
equations and their applications, {\bf 20}, \rm 123-154, Birkh\"auser.
\vskip 4 pt
\noindent
[G08] R. Goodman (2008), \it Harmonic analysis on compact symmetric spaces : \rm the
legacy of Elie Cartan and Hermann Weyl in Groups and analysis, 
London Math. Soc. Lecture Note,\bf 354, \rm 1-23.
\vskip 4 pt
\noindent
[HKMO12]  (2012) J. Hilgert, T. Kobayashi,  J. M\"ollers, B. Orsted (2012), \it  Fock model and Segal-Bargmann transform  for minimal representations of Hermitian Lie groups, \rm Journal of Functional Analysis,\bf  263, \rm 3492-3563.
\vskip 4 pt
\noindent
[M78] K. McCrimmon (1978), \it Jordan algebras and their applications, \rm 
Bull. A.M.S., \bf 84, \rm 612-627.
\vskip 4 pt
\noindent
[PW86]R.B. Paris and A.D. Wood (1986), \it Asymptotics of high order differential
equations, \rm Pitman Research Notes in Math Series, vol.
\bf 129,\rm Longman Scientific and Technical, Harlow.
\vskip 4 pt
\noindent
[S92] S. Sahi (1992), \it Explicit Hilbert spaces for certain unipotent representations, \rm 
Invent. Math.,\bf 110, no. 2, \rm 409-418.
\vskip 4 pt
\noindent
[V91] D. A. Vogan (1991), \it Associated varieties and unipotent representations|,  \rm in: W. Barker and P. Sally eds., Harmonic Analysis on Reductive Groups, Birkhausser, 315-388.

\end